\documentclass[11pt]{amsart}

\renewcommand{\bar}{\overline}

\newcommand{\pa}{\partial}
\renewcommand{\phi}{\varphi}

\newcounter{hours}\newcounter{minutes}

\newcommand{\ka}{K\"ahler }

\newcommand{\pz}{\partial_z}
\newcommand{\pzb}{\partial_{\bar z}}

\usepackage{amsmath,amsthm,amscd}
\usepackage{calc}

\title
[]{Canonical Metrics on the Moduli Space of Riemann Surfaces I}
\author{Kefeng Liu}
\author{Xiaofeng Sun}
\author{Shing-Tung Yau}
\date{January 31, 2004}

\thanks{The authors are supported by NSF. }

\newtheorem{theorem}{Theorem}[section]

\newtheorem{lemma}{Lemma}[section]
\newtheorem{cor}{Corollary}[section]
\newtheorem{prop}{Proposition}[section]

\newtheorem{definition}{Definition}[section]

\theoremstyle{remark}
\newtheorem{rem}{Remark}[section]

\begin{document}
\maketitle

\numberwithin{equation}{section}

\tableofcontents

\newcommand{\M}{{\mathcal M}}
\section{Introductions}

One of the main purpose of this paper is to compare those
well-known canonical and complete metrics on the Teichm\"uller and
the moduli spaces of Riemann surfaces. We use as bridge two new
metrics, the Ricci metric and the perturbed Ricci metric. We will
prove that these metrics are equivalent to those classical
complete metrics. For this purpose we study in detail the
asymptotic behaviors and the signs of the curvatures of these new
metrics. In particular we prove that the perturbed Ricci metric is
a complete K\"ahler metric with bounded negative holomorphic
sectional curvature and bounded bisectional and Ricci curvature.

The study of the Teichm\"uller spaces and moduli spaces of Riemann
surfaces has a long history. It has been intensively studied by
many mathematicians in complex analysis, differential geometry,
topology and algebraic geometry for the past 60 years. They have
also appeared in theoretical physics such as string theory. The
moduli space can be viewed as the quotient of the corresponding
Teichm\"uller space by the modular group. There are several
classical metrics on these spaces: the Weil-Petersson metric, the
Teichm\"uller metric, the Kobayashi metric, the Bergman metric,
the Caratheodory metric and the K\"ahler-Einstein metric. These
metrics have been studied over the years and have found many
important applications in various areas of mathematics. Each of
these metrics has its own advantages and disadvantages in studying
different problems.

The Weil-Petersson metric is a K\"ahler metric as first proved by
Ahlfors, both of its holomorphic sectional curvature and Ricci
curvature have negative upper bounds as conjectured by Royden and
proved by Wolpert. These properties have found many applications
by Wolpert, and they were also used in solving problems from
algebraic geometry by combining with the Schwarz lemma of Yau
(\cite{liu1}, \cite{Yau1}). But as first proved by Masur it is not
a complete metric which prevents the understanding of some aspects
of the geometry of the moduli spaces. Siu and Schumacher extended
some results to higher dimensional cases. The works of Masur and
Wolpert, Siu and Schumacher will play important roles in our
study.

The Teichm\"uller metric, the Kobayashi metric and the
Caratheodory metric are only Finsler metrics. They are very
effective in studying the hyperbolic property of the moduli space.
Royden proved that the Teichm\"uller metric is equal to the
Kobayashi metric from which he deduced the important corollary
that the isometry group of the Teichm\"uller space is exactly the
modular group. Recently C. McMullen introduced a new complete
K\"ahler metric on the moduli space by perturbing the
Weil-Petersson metric \cite{mc}. By using this metric he was able
to prove that the moduli space is K\"ahler hyperbolic, and also to
derive several topological consequences. The McMullen metric has
bounded geometry, but we lose control on the signs of its
curvatures.

In the early 80s Cheng-Yau \cite{cy1} and Mok-Yau \cite{mokyau1}
proved the existence of the K\"ahler-Einstein metrics on the
Teichm\"uller space. Since the K\"ahler-Einstein metric is
canonical, it also descends to a complete K\"ahler metric on the
moduli space. More than 20 years ago Yau \cite{yau2} conjectured
the equivalence of the K\"ahler-Einstein metric to the
Teichm\"uller metric. We will prove this conjecture in this paper.
Since the McMullen metric is equivalent to the Teichm\"uller
metric, so we have also proved the equivalence of the
K\"ahler-Einstein metric and the McMullen metric.

The method of our proof is to study in detail another complete
K\"ahler metric, the
 metric induced by the negative Ricci curvature of the
 Weil-Petersson metric which we call the Ricci metric. We first study its
 asymptotic behavior near the boundary of the moduli space, we prove that it is asymptotically
 equivalent to the Poincar\'e metric, and asymptotically its holomorphic sectional curvature
has negative upper and lower bound in the degeneration directions.
But its curvatures in the non-degeneration directions near the
boundary and in the interior of the moduli space can not be
controlled well. To solve this
 problem, we introduce another new complete K\"ahler metric which we call the perturbed
Ricci metric, it is obtained by
 adding a multiple of the Weil-Petersson metric. We
 compute the holomorophic sectional curvature and the Ricci
 curvature of this new metric. We show that they are all bounded below and above, and the
 holomorphic sectional curvature has negative upper and lower bounds. By
applying
 the Schwarz lemma of Yau we can prove the equivalence of this new
 metric to the K\"ahler-Einstein metric. The equivalence of the perturbed Ricci
 metric to the McMullen metric is proved by a careful estimate of
 the asymptotic behavior of these two metrics.

To state our main results in detail, let us introduce some
definitions and notations. Here for convenience we will use the
same notation for a K\"ahler metric and its K\"ahler form. First
two metrics $\omega_{\tau_1}$ and $\omega_{\tau_2}$ are called
equivalent, if they are quasi-isometric to each other in the sense
that

$$ C^{-1} \omega_{\tau_2} \leq \omega_{\tau_1} \leq C\omega_{\tau_2} $$
for some positive constant $C$. We will write this as
$\omega_{\tau_1}\sim \omega_{\tau_2}$.

Our first result is the following asymptotic behavior of the Ricci
metric near the boundary divisor of the moduli space. Let
$\mathcal T_g$ denote the Teichm\"uller space and $\M_g$ be the
moduli space of Riemann surfaces of genus $g$ where $g \geq 2$.
$\M_g$ is a complex orbifold of dimension $3g-3$ as a quotient of
$\mathcal T_g$ by the modular group. Let $n=3g-3$. Let
$\omega_{WP}$ denote the Weil-Petersson metric and
$\omega_\tau=-Ric(\omega_{WP})$ be the Ricci metric. It is easy to
show that there is an asymptotic Poincar\'e metric on $\mathcal
M_g$. See Section \ref{sec4} for the construction.

\begin{theorem}
The Ricci metric is equivalent to the asymptotic Poincar\'e
metric.
\end{theorem}

This theorem is proved in Section 4. Our second result is the
following estimates of the holomorphic sectional curvature of the
Ricci metric. Note our convention of the sign of the curvature may
be different from some literature.

\begin{theorem}
Let $X_0\in\bar{\mathcal{M}_g}\setminus\mathcal{M}_g$ be a
codimension $m$ point and let
$(t_1,\cdots,t_m,s_{m+1},\cdots,s_n)$ be the pinching coordinates
at $X_0$ where $t_1,\cdots,t_m$ correspond to the degeneration
directions. Then the holomorphic sectional curvature of the Ricci
metric is negative in the degeneration directions and is bounded
in the non-degeneration directions. Precisely, there is a
$\delta>0$ such that if $|(t,s)|<\delta$, then
\begin{eqnarray*}
\widetilde R_{i\bar ii\bar i}=
\frac{3u_i^4}{8\pi^4|t_i|^4}(1+O(u_0)) >0 \text{\hskip 0.2in if
\hskip 0.1in} i \leq m
\end{eqnarray*}
and
\begin{eqnarray*}
\widetilde R_{i\bar ii\bar i}=O(1) \text{\hskip 0.2in if \hskip
0.1in} i \geq m+1.
\end{eqnarray*}
Furthermore, on $\mathcal M_g$ the holomorphic sectional
curvature, the bisectional curvature and the Ricci curvature of
the Ricci metric are bounded from above and below.
\end{theorem}

This is Theorem 4.4 of Section 4 of this paper. One of the main
purposes of our work was to find a natural complete metric whose
holomorphic sectional curvature is negative. To do this, we
introduce the perturbed Ricci metric. In Section 5 we will prove
the following theorem:

\begin{theorem}\label{perholocurv}
For suitable choice of positive constant $C$, the perturbed Ricci
metric
$$\omega_{\widetilde\tau}=\omega_\tau+C\omega_{WP}$$ is complete
and its holomorphic sectional curvatures are negative and bounded
from above and below by negative constants. Furthermore, the Ricci
curvature of the perturbed Ricci metric is bounded from above and
below.
\end{theorem}

Note that the perturbed Ricci metric is equivalent to the Ricci
metric, since its asymptotic behavior is dominated by the Ricci
metric. Now we denote the K\"ahler-Einstein metric of
Cheng-Mok-Yau by $\omega_{KE}$ which is another complete K\"ahler
metric on the moduli space. By applying the Schwarz lemma of Yau
we derive our fourth result in Section 6:

\begin{theorem} We have the equivalence of the following three
complete K\"ahler metrics on the moduli spaces of curves:

$$\omega_{KE}\sim\omega_\tau\sim\omega_{\tilde{\tau}}.$$

\end{theorem}

Our final result in this paper proved in Section 6 is the
equivalence of the Ricci metric and the perturbed Ricci metric to
the McMullen metric. Let us denote the McMullen metric by
$\omega_M$.

\begin{theorem} We have the equivalence of the following metrics: the McMullen metric,
the Ricci metric and the perturbed Ricci metric:

$$\omega_M\sim \omega_\tau \sim \omega_{\tilde{\tau}}.$$
\end{theorem}

As a corollary we know that these metrics are also equivalent to
the Teichm\"uller metric, the Kobayashi metric, and the
K\"ahler-Einstein metric. This proved the conjecture of Yau
\cite{yau2}. In the second part of this work, we will study the
Bergman metric and the Caratheodory metric. We believe that these
two metrics are also equivalent to the above metrics. We will also
study the goodness of the Ricci metric in the sense of Mumford,
discuss the bounded geometry of the K\"ahler-Einstein metric and
the perturbed Ricci metric, and study the stability of the tangent
bundle of the moduli space of curves.

This paper is organized as follows. In Section 2 we set up some
notations and introduce the Weil-Petersson metric and its
curvatures. In Section 3 we introduce various operators needed for
our computations, we compute and simplify the curvature of the
Ricci metric by using these operators and their various special
properties. This section consists of long and complicated
computations. Section 4 consists of several subtle estimates of
the Ricci metric and its curvatures near the boundary of the
moduli space. In Section 5 we introduce the perturbed Ricci
metric, compute its curvature and study its asymptotic behavior
near the boundary of the moduli space. These results are then used
in Section 6 to prove the equivalence of the several well-known
classical complete K\"ahler metrics as stated above. In the
appendix we add some details of the computations for the
convenience of the readers.

{\bf Acknowledgements} {\em The second author would like to thank
H. Cao, P. Li, Z. Lu, R. Schoen and R. Wentworth for their help
and encouragement.}

\section{The Weil-Petersson metric}

The purpose of this section is to set up notations for our
computations. We will introduce the Weil-Petersson metric and
recall some of its basic properties. Let $\M_g$ be the moduli
space of Riemann surfaces of genus $g$ where $g \geq 2$. $\M_g$ is
a complex orbifold of dimension $3g-3$. Let $n=3g-3$. Let $\frak
X$ be the total space and $\pi:\frak X \to \M_g$ be the projection
map. There is a natural metric, called the Weil-Petersson metric
which is defined on the orbiford $\M_g$ as follows:

Let $s_1, \cdots, s_n$ be holomorphic local coordinates near a
regular point $s \in \M_g$ and assume that $z$ is a holomorphic
local  coordinate on the fiber $X_s=\pi^{-1}(s)$. For the
holomorphic vector fields $\frac{\pa}{\pa s_1},\cdots,
\frac{\pa}{\pa s_n}$, there are vector fileds $v_1,\cdots, v_n$ on
$\frak X$ such that
\begin{enumerate}
\item  $\pi_*(v_i)=\frac{\pa}{\pa s_i}$ for $i=1,\cdots,n$; \item
$\bar\pa v_i$ are harmonic $TX_s$-valued $(0,1)$ forms for
$i=1,\cdots,n$.
\end{enumerate}
The vector fields $v_1,\cdots,v_n$ are called the harmonic lift of
the vectors $\frac{\pa}{\pa s_1},\cdots, \frac{\pa}{\pa s_n}$. The
existence of such harmonic vector fields was pointed out by
Siu~\cite{siu1}. In his work ~\cite{sc1} Schumacher gave an
explicit construction of such lift which we now describe.

Since $g \geq 2$, we can assume that each fiber is equipped with
the K\"ahler-Einstein, or the Poincare metric,
$\lambda=\frac{\sqrt{-1}}{2}\lambda(z,s)dz\wedge d\bar z$. The
K\"ahler-Einstein condition gives the following equation:
\begin{eqnarray}\label{ke100}
\pz\pzb\log\lambda=\lambda.
\end{eqnarray}
For the rest of this paper we denote $\frac{\partial}{\partial
s_i}$ by $\partial_i$  and $\frac{\pa}{\pa z}$ by $\pa_z$. Let
\begin{eqnarray}\label{smalla}
a_i=-\lambda^{-1}\partial_i\pzb\log\lambda
\end{eqnarray}
and let
\begin{eqnarray}\label{biga}
A_i=\pzb a_i.
\end{eqnarray}
Then we have the following

\begin{lemma}\label{har}
The harmonic horizontal lift of $\partial_i$ is
\[
v_i=\partial_i+a_i\pz.
\]
In particular
\[
B_i=A_i\pz\otimes d\bar z \in H^{1}(X_s,T_{X_s})
\]
is harmonic. Further more, the lift $\pa_i\mapsto B_i$ gives the
Kodaira-Spencer map  $T_s\M_g\rightarrow H^1(X_s,T_{X_s})$.
\end{lemma}

Now we have the well-known definition of the Weil-Petersson
metric:

\begin{definition}
The Weil-Petersson metric on $\M_g$ is defined to be
\begin{eqnarray}\label{wp}
h_{i\bar j}(s)=\int_{X_s}B_i \cdot \bar{B_j} \ dv=
\int_{X_s}A_i\bar{A_j}\ dv,
\end{eqnarray}
where $dv=\frac{\sqrt{-1}}{2}\lambda dz\wedge d\bar z$ is the
volume form on the fiber $X_s$.
\end{definition}

It is known that the curvature tensor of the Weil-Petersson metric
can be represented by
\[
R_{i\bar j k\bar l}= \int_{X_s}\left\{(B_i\cdot \bar
B_j)(\Box+1)^{-1} (B_k\cdot\bar B_l) + (B_i\cdot \bar
B_l)(\Box+1)^{-1} (B_k\cdot\bar B_j)\right\}dv,
\]
where $\Box$ is the complex Laplacian defined by
\[
\Box=-\lambda^{-1}\frac{\pa^2}{\pa z\pa \bar z}.
\]

By the expression of the curvature operator, we know
 that the curvature operator
is nonpositive. Furthermore, the Ricci curvature of the metric is
negative.

However, the Weil-Petersson metric is incomplete. In \cite{tr1}
Trapani proved the negative Ricci curvature of the Weil-Petersson
metric is a complete K\"ahler metric on the moduli space. We call
this metric the Ricci metric. It is interesting to understand the
curvature of the Ricci metric, at least asymptotically. To
estimate it, we first derive an integral formula of its curvature.

\section{Ricci metric and its curvature}\label{sec3}
In this section we establish an integral formula \eqref{finalcurv}
of the curvature of the Ricci metric. The importance of this
formula is that the functions being integrated only involve
derivatives in the fiber direction which we are able to control.
Thus we can use this formula to estimate the asymptotics of the
curvature of the Ricci metric in next section.

The main tool we use is the harmonic lift of Siu and Schumacher
described in the previous section. These lifts together with
formula \eqref{vol0} enable us to transfer derivatives in the
moduli direction into derivatives in the fiber direction.

We use the same notations as in the previous section. We first
introduce several operators which will be used for the
computations and simplifications of the curvatures of the Ricci
metric.

 Define an
$(1,1)$ form on the total space $\frak X$ by
\[
g=\frac{\sqrt{-1}}{2}\pa\bar\pa\log\lambda
=\frac{\sqrt{-1}}{2}(g_{i\bar j} ds_i\wedge d\bar s_j-\lambda a_i
ds_i\wedge d\bar z- \lambda\bar a_i dz\wedge d\bar s_i+\lambda
dz\wedge d\bar z).
\]
The form $g$ is not necessarily positive.
 Introduce
 \[
e_{i\bar j}=\frac{2}{\sqrt{-1}}g(v_i, \bar v_j)=g_{i\bar
j}-\lambda a_i\bar{a_j}
\]
 be a global function. Let us write
$f_{i\bar j}=A_i\bar{A_j}$.  Schumacher proved the following
result:
\begin{lemma}\label{inner}
By using the same notations as above, we have
\begin{eqnarray}\label{dot}
(\Box+1)e_{i\bar j}=f_{i\bar j}.
\end{eqnarray}
\end{lemma}

Since $e_{i\bar j}$ and $f_{i\bar j}$ are the building blocks of
the Ricci metric, it is interesting to study its property under
the action of the vector fields $v_i$'s.

\begin{lemma}
With the same notations as above, we have
\[
v_k(e_{i\bar j})=v_i(e_{k\bar j}).
\]
\end{lemma}

{\bf Proof.} Since  $dg=0$, we have the following
\begin{align*}
&0=dg(v_i, v_k, \bar v_j)=v_i(e_{k\bar j})-v_k(e_{i\bar j})
+\bar v_j g(v_i, v_k)\\
& \qquad -g(v_i, [v_k, \bar v_j])+g(v_k, [v_i, \bar v_j]) -g(\bar
v_j, [v_i, v_k]).
\end{align*}
The Lie bracket of $v_j$ with $\bar v_j$ or $v_k$ are vector
fields tangent to $X_s$, which are perpendicular to the horizontal
vector fields $v_i$ with respect to the form $g$. Thus the last
three terms of the above equations are zero. On the other hand,
$g(v_i, v_k)=0$. The lemma thus follows from the above equation.

\qed

We also need to define the following operator
\[
P: \ C^\infty(X_s)\rightarrow \Gamma(\Lambda^{1,0} (T^{0,1}
X_s)),\,\, f\mapsto \pz(\lambda^{-1}\pz f).
\]
The dual operator $P^*$ can be written as follows
\[
P^*:\ \Gamma(\Lambda^{0,1} (T^{1,0} X_s))\rightarrow
C^\infty(X_s),\,\, B\mapsto
\lambda^{-1}\pz(\lambda^{-1}\pz(\lambda B)).
\]

The operator $P$ is actually a composition of the Maass operators.
We recall the definitions from \cite{wol1}. Let $X$ be a Riemann
surface and let $\kappa$ be its canonical bundle. For any integer
$p$, let $S(p)$ be the space of smooth sections of
$(\kappa\otimes\bar\kappa^{-1})^{\frac{p}{2}}$. Fix a conformal
metric $ds^2=\rho^2(z)|dz|^2$.

\begin{definition}\label{maass}
The Maass operators $K_p$ and $L_p$ are defined to be the metric
derivatives  $K_p:\ S(p)\to S(p+1)$ and $L_p:\ S(p)\to S(p-1)$
given by
\[
K_p(\sigma)=\rho^{p-1}\pz(\rho^{-p}\sigma)
\]
and
\[
L_p(\sigma)=\rho^{-p-1}\pzb(\rho^{p}\sigma)
\]
where $\sigma \in S(p)$.
\end{definition}
Clearly we have $P=K_1K_0$. Also each element $\sigma \in S(p)$
has a well-defined absolute value $|\sigma|$ which is independent
of the choice of the local coordinate. We define the $C^k$ norm of
$\sigma$ as in \cite{wol1}:
\begin{definition}\label{norm}
Let $Q$ be an operator which is a composition of operators
$K_\ast$ and $L_\ast$. Denote by $|Q|$ the number of such factors.
For any $\sigma \in S(p)$, define
\[
\Vert\sigma\Vert_{0}=\sup_{X}|\sigma|
\]
and
\[
\Vert\sigma\Vert_k=\sum_{|Q|\leq k}\Vert Q\sigma\Vert_0.
\]
We can also localize the norm on a subset of $X$. Let $\Omega
\subset X$ be a domain. We can define
\[
\Vert\sigma\Vert_{0,\Omega}=\sup_{\Omega}|\sigma|
\]
and
\[
\Vert\sigma\Vert_{k,\Omega} =\sum_{|Q|\leq k}\Vert
Q\sigma\Vert_{0,\Omega}.
\]
\end{definition}
Both of the above definitions depend on the choice of conformal
metric on $X$. In the following, we always use the
K\"ahler-Einstein metric on the surface unless otherwise stated.

Since the Weil-Petersson metric is defined by using the integral
along the fibers, the following formula is very useful:

\begin{equation}\label{vol0}
\partial_i \int_{X_s} \eta=\int_{X_s}L_{v_i}\eta
\end{equation}
where $\eta$ is a relative $(1,1)$ form on $\frak X$.

The Lie derivative defined here is slightly different from the
ordinary definition. Let $\phi_t$ be the one parameter group
generated by the vector field $v_i$. Then $\phi_t$ can be viewed
as a diffeomorphism between two fibers $X_s\rightarrow X_{s'}$.
Then we define
\[
L_{v_i}\eta=\underset{t\rightarrow 0}{\lim}\frac 1t(\phi^*_t
(\sigma)-\sigma)
\]
for any one form $\sigma$. On the other hand, let $\xi$ be a
vector field on the fiber $X_s$. Then we define
\[
L_{v_i}\xi=\underset{t\rightarrow 0}{\lim} \frac 1t
((\phi_{-t})_*\xi-\xi).
\]

We have the following
\begin{prop}
By using the above notations, we have
\[
L_{v_i}\sigma=i(v_i) d_1\sigma+d_1 i(v_i)\sigma,
\]
where $d_1$ is the differential operator along the fiber, and
\[
L_{v_i}\xi=[v_i,\xi].
\]
\end{prop}

In the following, we denote $L_{v_i}$ by $L_i$.

\begin{lemma}\label{us1}
By using the above notations, we have
\begin{enumerate}
\item $L_{i} dv=0$;  \item $L_{\bar l}(B_i)=-\bar P(
e_{i\bar l})-f_{i\bar l}\pzb\otimes d\bar z + f_{i\bar
l}\pz\otimes dz$; \item $L_k(\bar{B_j})=-P (e_{k\bar j})-f_{k\bar
j}\pz\otimes dz+ f_{k\bar j}\pzb\otimes d\bar z$; \item
$L_k(B_i)=(v_k(A_i)-A_i\pz a_k) \pz\otimes d\bar z$; \item
 $L_{\bar l}(\bar{A_j})=
(\bar{v_l}(\bar{A_l})-\bar{A_l}\pzb\bar{a_l}) \pzb\otimes dz$.
\end{enumerate}

\end{lemma}

{\bf Proof.} The first formula was proved by Schumacher in
\cite{sc1}. To check the other formulas, we note that the third
and fifth formulas follow from the second and fourth, which we
will prove, by taking conjugation. We first have
\begin{align*}
\begin{split}
\pz a_k &=\pz(-\lambda^{-1}\partial_k\pzb\log\lambda)=
\lambda^{-2}\pz\lambda\partial_k\pzb\log\lambda
-\lambda^{-1}\pz\partial_k\pzb\log\lambda \\
&=-\lambda^{-1}\pz\lambda a_k-\lambda^{-1}\partial_k\pz\pzb
\log\lambda=-\lambda^{-1}\pz\lambda
a_k-\lambda^{-1}\partial_k\lambda.
\end{split}
\end{align*}
We also have
\begin{align*}
\begin{split}
\partial_{\bar l}a_i &=\partial_{\bar l}
(-\lambda^{-1}\partial_i\pzb\log\lambda))=
\lambda^{-2}\partial_{\bar l}\lambda\partial_i\pzb\log\lambda
-\lambda^{-1}\pzb\partial_i\partial_{\bar l}\log\lambda \\
&=-\lambda^{-1}\partial_{\bar l}\lambda a_i-\lambda^{-1} \pzb
g_{i\bar l}= -\lambda^{-1}\partial_{\bar l}\lambda
a_i-\lambda^{-1}
\pzb (e_{i\bar l}+\lambda a_i\bar{a_l})\\
&=-\lambda^{-1}\partial_{\bar l}\lambda a_i-\lambda^{-1} \pzb
e_{i\bar l}-\lambda^{-1}\pzb\lambda a_i\bar{a_l}
-A_i\bar{a_l}-a_i\pzb\bar{a_l}\\
&=-(\lambda^{-1}\partial_{\bar l}\lambda+
\lambda^{-1}\pzb\lambda\bar{a_l}+\pzb\bar{a_l})a_i
-\lambda^{-1}\pzb e_{i\bar l}-A_i\bar{a_l}\\
&=-\lambda^{-1}\pzb e_{i\bar l}-A_i\bar{a_l}.
\end{split}
\end{align*}

For the second formula we have
\begin{align*}
\begin{split}
L_{\bar l}(B_i)&=\bar{v_l}(A_i)\pz\otimes d\bar z
+A_i(-\pz\bar{a_l}\pzb)\otimes d\bar z
+A_i\pz\otimes(\pz\bar{a_l}dz+\pzb\bar{a_l}d\bar z)\\
&=(\bar{v_l}(A_i)+A_i\pzb\bar{a_l})\pz\otimes d\bar z -f_{i\bar
l}\pzb\otimes d\bar z + f_{i\bar l}\pz\otimes dz.
\end{split}
\end{align*}
So we only need to check that $\bar{v_l}(A_i)+A_i\pzb\bar{a_l}=
-\pzb(\lambda^{-1}\pzb e_{i\bar l})$. To prove this, we have
\begin{align*}
\begin{split}
\bar{v_l}(A_i)+A_i\pzb\bar{a_l}&=\bar{a_l}\pzb
A_i+\partial_{\bar{l}}A_i+A_i\pzb\bar{a_l}
=\pzb(A_i\bar{a_l})+\pzb\partial_{\bar l}a_i\\
&=\pzb(A_i\bar{a_l})-\pzb(\lambda^{-1}\pzb e_{i\bar l})
-\pzb(A_i\bar{a_l})=-\pzb(\lambda^{-1}\pzb e_{i\bar l}).
\end{split}
\end{align*}
This proved the second formula. For the fourth one, we have
\begin{eqnarray*}
L_k(B_i)=v_k(A_i)\pz\otimes d\bar z+A_i(-\pz a_k \pz) \otimes
d\bar z=(v_k(A_i)-A_i\pz a_k)\pz\otimes d\bar z.
\end{eqnarray*}
This finishes the proof.

\qed

An interesting and useful fact is that the Lie derivative of $B_i$
in the direction of $v_k$ is still harmonic. This result is true
only for the moduli space of Riemann surfaces. In the general case
of moduli space of \ka-Einstein manifolds, we only have
$\bar{\partial}^{\ast}L_k B_i=0$.
\begin{lemma}\label{us2}
$L_k(B_i) \in H^{1}(X_s,T{X_s})$ is harmonic.
\end{lemma}
{\bf Proof.} From Lemma \ref{us1} we know that $L_k(B_i)
=(v_k(A_i)-A_i\pz a_k) \pz\otimes d\bar z\in
H^{0,1}(X_s,T_{X_s})$. So it is clear that
$\bar\partial(L_k(B_i))=0$. To prove
$\bar\partial^{\ast}(L_k(B_i))=0$ we only need to check that
\[
\pz(\lambda(v_k(A_i)-A_i\pz a_k))=0.
\]
From the computation in the above lemma, we have
\begin{align*}
\begin{split}
v_k(A_i)-A_i\pz a_k=&\lambda a_ia_k-
\pzb(\lambda^{-1}\partial_k\partial_i\pzb\log\lambda)\\
=&\lambda a_ia_k+\lambda^{-2}\pzb\lambda
\partial_k\partial_i\pzb\log\lambda
-\lambda^{-1}\partial_k\partial_i\pzb\pzb\log\lambda
\end{split}
\end{align*}
which implies
\begin{align}\label{f10}
\begin{split}
\pz(\lambda(v_k(A_i)-A_i\pz a_k))=& \pz(\lambda^2
a_ia_k+\lambda^{-1}\pzb\lambda
\partial_k\partial_i\pzb\log\lambda
-\partial_k\partial_i\pzb\pzb\log\lambda)\\
=&\pz(\lambda^2 a_ia_k)+\pz(\lambda^{-1}\pzb\lambda)
\partial_k\partial_i\pzb\log\lambda+\lambda^{-1}\pzb\lambda
\pz(\partial_k\partial_i\pzb\log\lambda)\\
&-\partial_k\partial_i\pzb\pz\pzb\log\lambda\\
=&\pz(\lambda^2 a_ia_k)+\lambda\partial_k\partial_i\pzb\log\lambda
+\lambda^{-1}\pzb\lambda\partial_k\partial_i\lambda
-\partial_k\partial_i\pzb\lambda.
\end{split}
\end{align}
Now we analyze the second term in \eqref{f10}. We have
\begin{align}\label{f20}
\begin{split}
\lambda\partial_k\partial_i\pzb\log\lambda=&
\lambda\partial_k\partial_i\frac{\pzb\lambda}{\lambda}
=\lambda\partial_k\frac{\lambda\partial_i\pzb\lambda
-\partial_i\lambda\pzb\lambda}{\lambda^2}\\
=&\lambda\frac{\lambda^2(\partial_k\lambda
\partial_i\pzb\lambda+\lambda
\partial_k\partial_i\pzb\lambda-\partial_k\partial_i
\lambda\pzb\lambda-\partial_i\lambda\partial_k\pzb\lambda)}
{\lambda^4}\\
&-\lambda\frac{2\lambda\partial_k\lambda
(\lambda\partial_i\pzb\lambda
-\partial_i\lambda\pzb\lambda)}{\lambda^4}\\
=&-\lambda^{-1}\partial_k\lambda
\partial_i\pzb\lambda
+\partial_k\partial_i\pzb\lambda -\lambda^{-1}\partial_k\partial_i
\lambda\pzb\lambda -\lambda^{-1}\partial_i\lambda
\partial_k\pzb\lambda\\
&+2\lambda^{-2}\partial_i\lambda
\partial_k\lambda\pzb\lambda\\
=&-\partial_i\lambda(\lambda^{-1}\partial_k\pzb\lambda
-\lambda^{-2}\partial_k\lambda\pzb\lambda)
-\partial_k\lambda(\lambda^{-1}\partial_i\pzb\lambda
-\lambda^{-2}\partial_i\lambda\pzb\lambda)\\
&+\partial_k\partial_i\pzb\lambda
-\lambda^{-1}\partial_k\partial_i
\lambda\pzb\lambda\\
=&-\partial_i\lambda\partial_k\pzb\log\lambda
-\partial_k\lambda\partial_i\pzb\log\lambda
+\partial_k\partial_i\pzb\lambda -\lambda^{-1}\partial_k\partial_i
\lambda\pzb\lambda\\
=&\lambda\partial_i\lambda a_k+\lambda\partial_k\lambda a_i
+\partial_k\partial_i\pzb\lambda -\lambda^{-1}\partial_k\partial_i
\lambda\pzb\lambda.
\end{split}
\end{align}
By combining \eqref{f10} and \eqref{f20} we have
\begin{align*}
\begin{split}
\pz(\lambda(v_k(A_i)-A_i\pz a_k))=&\pz(\lambda^2 a_ia_k)
+\lambda\partial_i\lambda a_k+\lambda\partial_k\lambda a_i\\
=&2\lambda\pz\lambda a_ia_k+\lambda^2 \pz a_ia_k +\lambda^2 a_i\pz
a_k
+\lambda\partial_i\lambda a_k+\lambda\partial_k\lambda a_i\\
=&\lambda^2 a_k(\lambda^{-1}\pz\lambda a_i+\pz a_i
+\lambda^{-1}\partial_i\lambda)\\
&+\lambda^2 a_i(\lambda^{-1}\pz\lambda a_k+\pz a_k
+\lambda^{-1}\partial_k\lambda)\\
=&0.
\end{split}
\end{align*}
This proves that $\bar\partial^{\ast}(L_k(B_i))=0$.

\qed

The above lemma is very helpful in computing the curvature when we
use normal coordinates of the Weil-Petersson metric. We have
\begin{cor}\label{normalgood}
Let $s_1,\cdots,s_n$ be normal coordinates at $s \in \mathcal M_g$
with respect to the Weil-Petersson metric. Then at $s$ we have,
for all $i,k$,
\[
L_k B_i=0.
\]
\end{cor}
{\bf Proof.} From Lemma \ref{us2} we know that $L_k B_i$ is
harmonic. Since $B_1,\cdots,B_n$ is a basis of $T_s \mathcal M_g$,
we have
\[
L_k B_i=h^{p\bar q}(\int_{X_s}L_k B_i\cdot \bar{B_q}\
dv)B_p=h^{p\bar q}\partial_k h_{i\bar q}B_p=0.
\]

\qed

The commutator of $v_k$ and $\bar{v_l}$ will be used later. We
give a formula here which is essentially due to Schumacher.
\begin{lemma}\label{commu0}
$[\bar{v_l},v_k]=-\lambda^{-1}\pzb e_{k\bar l}\pz +\lambda^{-1}\pz
e_{k\bar l}\pzb$.
\end{lemma}
{\bf Proof.} From a direct computation we have
\[
[\bar{v_l},v_k]=\bar{v_l}(a_k)\pz-v_k(\bar{a_l})\pzb.
\]
By using Lemma \ref{us1} we have
\[
\bar{v_l}(a_k)=\bar{a_l}\pzb a_k+\partial_{\bar l}a_k
=-\lambda^{-1}\pzb e_{k\bar l}
\]
and
\[
v_k(\bar{a_l})=a_k\pz\bar{a_l}+\partial_k \bar{a_l}
=\lambda^{-1}\pz e_{k\bar l}.
\]
These finish the proof.

\qed

\begin{rem}
In the rest of this paper, we will use the following notation
for curvature:\\
Let $(M,g)$ be a K\"ahler manifold. Then the curvature tensor is
given by
\begin{eqnarray}\label{notation1}
R_{i\bar j k\bar l}=\frac{\partial^2 g_{i\bar j}}{\partial z_k
\partial \bar z_l}-g^{p\bar q}\frac{\partial g_{i\bar q}}
{\partial z_k}\frac{\partial g_{p\bar j}} {\partial \bar z_l}.
\end{eqnarray}
In this situation, the Ricci curvature is given by
\[
R_{i\bar j}=-g^{k\bar l}R_{i\bar j k\bar l}.
\]
\end{rem}

In \cite{siu1} and \cite{sc1}, Siu and Schumacher proved the
following curvature formula for the Weil-Petersson metric. This
formula was also proved by Wolpert in \cite{wol3}. Here we give a
short proof here.
\begin{theorem}\label{wpcurv}
The curvature of Weil-Petersson metric is given by
\begin{eqnarray}\label{WPcurv1}
R_{i\bar j k\bar l}=\int_{X_s}(e_{i\bar j}f_{k\bar l}+ e_{i\bar
l}f_{k\bar j})\ dv.
\end{eqnarray}
\end{theorem}
{\bf Proof.} We have
\begin{align}\label{wpcurv10}
\begin{split}
R_{i\bar j k\bar l}=&\partial_{\bar l}\partial_k h_{i\bar j}
-h^{p\bar q}\partial_k h_{i\bar q}\partial_{\bar l}
h_{p\bar j}\\
=&\partial_{\bar l}\int_{X_s}L_k B_i\cdot \bar{B_j}\ dv -h^{p\bar
q}\int_{X_s}L_k B_i\cdot \bar{B_q}\ dv \int_{X_s}B_p\cdot
L_{\bar l}\bar{B_j}\ dv\\
=&\int_{X_s}(L_{\bar l}L_k B_i\cdot \bar{B_j} +L_k B_i\cdot
L_{\bar l}\bar{B_j})\ dv -h^{p\bar q}\int_{X_s}L_k B_i\cdot
\bar{B_q}\ dv \int_{X_s}B_p\cdot L_{\bar l}\bar{B_j}\ dv.
\end{split}
\end{align}
Since $B_1,\cdots,B_n$ is a basis of $T_s M_g$, we have
\[
h^{p\bar q}\int_{X_s}L_k B_i\cdot \bar{B_q}\ dv \int_{X_s}B_p\cdot
L_{\bar l}\bar{B_j}\ dv =\int_{X_s}L_k B_i\cdot L_{\bar
l}\bar{B_j}\ dv.
\]
By combining this formula with \eqref{wpcurv10} we have
\begin{align}\label{wpcurv20}
\begin{split}
R_{i\bar j k\bar l}=& \int_{X_s}L_{\bar l}L_k B_i\cdot \bar{B_j}\
dv =\int_{X_s}L_k L_{\bar l}B_i\cdot \bar{B_j}\ dv
+\int_{X_s}L_{[\bar{v_l},v_k]}B_i\cdot \bar{B_j}\ dv\\
=&\partial_k\int_{X_s}L_{\bar l}B_i\cdot \bar{B_j}\ dv
-\int_{X_s}L_{\bar l}B_i\cdot L_k\bar{B_j}\ dv
+\int_{X_s}L_{[\bar{v_l},v_k]}B_i\cdot \bar{B_j}\ dv\\
=&-\int_{X_s}L_{\bar l}B_i\cdot L_k\bar{B_j}\ dv
+\int_{X_s}L_{[\bar{v_l},v_k]}B_i\cdot \bar{B_j}\ dv
\end{split}
\end{align}
since $\int_{X_s}L_{\bar l}B_i\cdot \bar{B_j}\ dv=0$. Now we
compute $\int_{X_s}L_{[\bar{v_l},v_k]}B_i\cdot \bar{B_j}\ dv$. Let
$\pi_{\bar 1}^1(L_{[\bar{v_l},v_k]}B_i)$ be the projection of
$L_{[\bar{v_l},v_k]}B_i$ onto $H^{0,1}(X_s,T_{X_s})$ which gives
the $\pz\otimes d\bar z$ part of $L_{[\bar{v_l},v_k]}B_i$. Since
$B_i$ is harmonic, we know $\pz(\lambda A_i)=0$ which implies $\pz
A_i= -\lambda^{-1}\pz\lambda A_i$. By Lemma \ref{commu0} we have
\begin{align}\label{wpcurv30}
\begin{split}
\pi_{\bar 1}^1(L_{[\bar{v_l},v_k]}B_i)=& (-\lambda^{-1}\pzb
e_{k\bar l}\pz A_i +A_i\pz(\lambda^{-1}\pzb e_{k\bar l})
+\pzb(\lambda^{-1}A_i\pz e_{k\bar l}))\pz\otimes d\bar z\\
=&(\lambda^{-2}\pz\lambda A_i\pzb e_{k\bar l}
-\lambda^{-2}\pz\lambda A_i\pzb e_{k\bar l} -A_i\Box e_{k\bar l}
+\pzb(\lambda^{-1}A_i\pz e_{k\bar l}))\pz\otimes d\bar z\\
=&(-A_i\Box e_{k\bar l} +\pzb(\lambda^{-1}A_i\pz e_{k\bar
l}))\pz\otimes d\bar z.
\end{split}
\end{align}
This implies
\begin{align}\label{wpcurv40}
\begin{split}
\int_{X_s}L_{[\bar{v_l},v_k]}B_i\cdot \bar{B_j}\ dv=&
\int_{X_s}\pi_{\bar 1}^1(L_{[\bar{v_l},v_k]}B_i)\cdot
\bar{B_j}\ dv\\
=&\int_{X_s}(-A_i\Box e_{k\bar l}
+\pzb(\lambda^{-1}A_i\pz e_{k\bar l}))\bar{A_j}\ dv\\
=&-\int_{X_s}f_{i\bar j}\Box e_{k\bar l}\ dv
+\int_{X_s}\pzb(\lambda^{-1}A_i\pz e_{k\bar l})\bar{A_j}\ dv\\
=&-\int_{X_s}f_{i\bar j}\Box e_{k\bar l}\ dv
-\int_{X_s}\lambda^{-2}A_i\pz e_{k\bar l}
\pzb(\lambda \bar{A_j})\ dv\\
=&-\int_{X_s}f_{i\bar j}\Box e_{k\bar l}\ dv.
\end{split}
\end{align}
To compute $\int_{X_s}L_{\bar l}B_i\cdot L_k\bar{B_j}\ dv$, by
using Lemma \ref{us1} we obtain
\begin{align}\label{wpcurv50}
\begin{split}
&\int_{X_s}L_{\bar l}B_i\cdot L_k\bar{B_j}\ dv=
\int_{X_s}(\pzb(\lambda^{-1}\pzb e_{i\bar l}) \pz(\lambda^{-1}\pz
e_{k\bar j})-2f_{k\bar j}f_{i\bar l})\
dv\\
=&\int_{X_s}(\lambda^{-2}\pz e_{k\bar j}
\pz(\lambda\pzb(\lambda^{-1}\pzb e_{i\bar l}))\ dv
-2\int_{X_s}f_{k\bar j}f_{i\bar l}\ dv\\
=&-\int_{X_s}(\lambda^{-2}\pz\lambda\pz e_{k\bar j}
\pzb(\lambda^{-1}\pzb e_{i\bar l})+\lambda^{-1} \pz e_{k\bar j}\pz
\pzb(\lambda^{-1}\pzb e_{i\bar l}))\ dv
-2\int_{X_s}f_{k\bar j}f_{i\bar l}\ dv\\
=&\int_{X_s}(\lambda^{-2}\pzb e_{i\bar l}\pzb
(\lambda^{-1}\pz\lambda\pz e_{k\bar j}) +\lambda^{-1}\pz\pzb
e_{k\bar j}\pz (\lambda^{-1}\pzb e_{i\bar l}))\ dv
-2\int_{X_s}f_{k\bar j}f_{i\bar l}\ dv\\
=&\int_{X_s}(\lambda^{-2}\pzb e_{i\bar l} (\lambda\pz e_{k\bar
j}-\pz\lambda\Box e_{k\bar j}) -\Box e_{k\bar j}
(-\lambda^{-2}\pz\lambda\pzb e_{i\bar l} -\Box e_{i\bar l}))\ dv
-2\int_{X_s}f_{k\bar j}f_{i\bar l}\ dv\\
=&\int_{X_s}(\lambda^{-1}\pzb e_{i\bar l}\pz e_{k\bar j}) +\Box
e_{k\bar j}\Box e_{i\bar l})\ dv
-2\int_{X_s}f_{k\bar j}f_{i\bar l}\ dv\\
=&\int_{X_s}(\Box e_{k\bar j} e_{i\bar l} +\Box e_{k\bar j}\Box
e_{i\bar l})\ dv
-2\int_{X_s}f_{k\bar j}f_{i\bar l}\ dv\\
=&\int_{X_s}(\Box e_{k\bar j} f_{i\bar l} -2f_{k\bar j}f_{i\bar
l})\ dv =-\int_{X_s}(f_{k\bar j}f_{i\bar l}+ e_{k\bar j}f_{i\bar
l})\ dv.
\end{split}
\end{align}
By combining \eqref{wpcurv20}, \eqref{wpcurv40} and
\eqref{wpcurv50} with the identity $f_{k\bar j}f_{i\bar l}
=A_i\bar{A_j}A_k\bar{A_l}=f_{i\bar j}f_{k\bar l}$, we have
\begin{align}\label{wpcurv60}
\begin{split}
R_{i\bar j k\bar l}=&\int_{X_s}(f_{k\bar j}f_{i\bar l}+ e_{k\bar
j}f_{i\bar l}- f_{i\bar j}\Box e_{k\bar l})\ dv
=\int_{X_s}(f_{i\bar j}e_{k\bar l}+f_{i\bar l}e_{k\bar j})
\ dv\\
=&\int_{X_s}(e_{i\bar j}f_{k\bar l}+ e_{i\bar l}f_{k\bar j})\ dv.
\end{split}
\end{align}
Here we have used the fact the $(\Box+1)$ is a self-adjoint
operator. This finished the proof.

\qed

It is well-known that the Ricci curvature of the Weil-Petersson
metric is negative which implies that the negative Ricci curvature
of the Weil-Petersson metric defines a \ka metric on the moduli
space $\M_g$.
\begin{definition}
The Ricci metric $\tau_{i\bar j}$ on the moduli space $\M_g$ is
the negative Ricci curvature of the Weil-Petersson metric. That is
\begin{eqnarray}\label{riccidef}
\tau_{i\bar j}=-R_{i\bar j}=h^{\alpha\bar\beta} R_{i\bar
j\alpha\bar\beta}.
\end{eqnarray}
\end{definition}

Now we define a new operator which acts on functions on the
fibers.
\begin{definition}
For each $1 \leq k \leq n$ and for any smooth function $f$ on the
fibers, we define the commutator operator $\xi_k$ which acts on a
function $f$ by
\begin{eqnarray}\label{oper}
\xi_k(f)=\bar\partial^{\ast}(i(B_k)\partial f)
=-\lambda^{-1}\pz(A_k\pz f).
\end{eqnarray}
\end{definition}
The reason we call $\xi_k$ the commutator operator is that $\xi_k$
is the commutator of $(\Box+1)$ and $v_k$ and the following lemma.
\begin{lemma}\label{commu}
As operators acting on functions, we have \\
(1) $(\Box+1)v_k-v_k(\Box+1)=\Box
v_k-v_k\Box=\xi_k$;\\
(2) $(\Box+1)\bar{v_l}-\bar{v_l}(\Box+1)=\Box\bar{v_l}
-\bar{v_l}\Box=\bar{\xi_l}$;\\
(3) $\xi_k(f)=-A_k\pz(\lambda^{-1}\pz f)=A_kP(f)
=-A_kK_1K_0(f)$. \\
Furthermore, we have
\begin{eqnarray}\label{Tinv}
(\Box+1)v_k(e_{i\bar j})=\xi_k(e_{i\bar j}) +\xi_i(e_{k\bar
j})+L_k B_i \cdot \bar{B_j}.
\end{eqnarray}
\end{lemma}
{\bf Proof.} To prove (1), we have
\begin{align*}
\begin{split}
(\Box+1)v_k-v_k(\Box+1)=&\Box v_k+v_k-v_k\Box-v_k
=\Box v_k-v_k\Box\\
=&-\lambda^{-1}\pz\pzb(a_k\pz+\partial_k)-
(a_k\pz+\partial_k)(-\lambda^{-1}\pz\pzb)\\
=&-\lambda^{-1}\pz(A_k\pz+a_k\pz\pzb+\partial_k\pzb)\\
&+a_k\pz(\lambda^{-1})\pz\pzb+\lambda^{-1}a_k\pz\pz\pzb
+\partial_k(\lambda^{-1})\pz\pzb
+\lambda^{-1}\partial_k\pz\pzb\\
=&-\lambda^{-1}\pz(A_k\pz)-\lambda^{-1}\pz a_k\pz\pzb
-\lambda^{-1}a_k\pz\pz\pzb-\lambda^{-1}\partial_k\pz\pzb\\
&-\lambda^{-2}\pz\lambda a_k\pz\pzb+\lambda^{-1}a_k\pz\pz\pzb
-\lambda^{-2}\partial_k\lambda\pz\pzb
+\lambda^{-1}\partial_k\pz\pzb\\
=&\xi_k-\lambda^{-1}(\pz a_k+\lambda^{-1}\pz\lambda a_k
+\lambda^{-1}\partial_k\lambda)\pz\pzb=\xi_k
\end{split}
\end{align*}
where we have used Lemma \ref{us1} in the last equality of the
above formula. By taking conjugation we can prove (2) by using
(1). To prove (3), we use the harmonicity of $B_k$. Since
$\bar\partial^{\ast}B_k=0$ we have $\pz(\lambda A_k)=0$. So
\begin{eqnarray*}
\xi_k(f)=-\lambda^{-1}\pz(A_k\pz f) =-\lambda^{-1}\pz(\lambda A_k
\lambda^{-1}\pz f) =-\lambda^{-1}\lambda A_k\pz(\lambda^{-1}\pz f)
=-A_k\pz(\lambda^{-1}\pz f).
\end{eqnarray*}
To prove the last part, by using part 1 of this lemma, we have
\begin{align*}
\begin{split}
(\Box+1)v_k(e_{i\bar j})=&v_k((\Box+1)(e_{i\bar j}))
+\xi_k(e_{i\bar j})=v_k(f_{i\bar j})+\xi_k(e_{i\bar j})\\
=&L_k B_i \cdot \bar{B_j}+B_i \cdot L_k\bar{B_j} +\xi_k(e_{i\bar
j}) =L_k B_i \cdot \bar{B_j}-A_i\pz(\lambda^{-1}\pz
e_{k\bar j})+\xi_k(e_{i\bar j})\\
=&L_k B_i \cdot \bar{B_j}+\xi_i(e_{k\bar j})+\xi_k(e_{i\bar j}).
\end{split}
\end{align*}
This finishes the proof.

\qed
\begin{rem}
From Corollary \ref{normalgood} and the above lemma, when we use
the normal coordinates on the moduli space, we have the clean
formula $(\Box+1)v_k(e_{i\bar j})=\xi_i(e_{k\bar
j})+\xi_k(e_{i\bar j})$.
\end{rem}

The main result in this section is to prove the curvature formula
of the Ricci metric. The terms produced here are very symmetric
with respect to indices. For convenience, we introduce the
symmetrization operator.
\begin{definition}
Let $U$ be any quantity which depends on indices $i,k,\alpha,\bar
j,\bar l, \bar\beta$. The symmetrization operator $\sigma_1$ is
defined by taking the summation of all orders of the triple
$(i,k,\alpha)$. That is
\begin{align*}
\begin{split}
\sigma_1(U(i,k,\alpha,\bar j,\bar l, \bar\beta))=&
U(i,k,\alpha,\bar j,\bar l, \bar\beta)+ U(i,\alpha,k,\bar j,\bar
l, \bar\beta)+ U(k,i,\alpha,\bar j,\bar l, \bar\beta)+
U(k,\alpha,i,\bar j,\bar l, \bar\beta)\\
&+U(\alpha,i,k,\bar j,\bar l, \bar\beta)+ U(\alpha,k,i,\bar j,\bar
l, \bar\beta).
\end{split}
\end{align*}
Similarly, $\sigma_2$ is the symmetrization operator of $\bar j$
and $\bar \beta$ and $\widetilde{\sigma_1}$ is the symmetrization
operator of $\bar j$, $\bar l$ and $\bar \beta$.
\end{definition}

Now we are ready to compute the curvature of the Ricci metric. For
the first order derivative we have
\begin{theorem}\label{1stderiv}
\begin{align}\label{1st}
\begin{split}
\partial_k \tau_{i\bar j}=h^{\alpha\bar\beta}
\left\{\sigma_1\int_{X_s} (\xi_k(e_{i\bar j})e_{\alpha\bar\beta})\
dv \right\} +\tau_{p\bar j}\Gamma_{ik}^p
\end{split}
\end{align}
where $\Gamma_{ik}^p$ is the Christoffell symbol of the
Weil-Petersson metric.
\end{theorem}
{\bf Proof.} From Lemma \ref{inner} we know that $(\Box+1)e_{i\bar
j}=f_{i\bar j}$. By using Lemma \ref{commu} and Theorem
\ref{wpcurv} we have
\begin{align}\label{1st10}
\begin{split}
\partial_k R_{i\bar j\alpha\bar\beta}=&
\partial_k \int_{X_s}(e_{i\bar j}f_{\alpha\bar\beta}+
e_{i\bar\beta}f_{\alpha\bar j})\ dv\\
=& \int_{X_s}(v_k(e_{i\bar j})f_{\alpha\bar\beta}+ e_{i\bar
j}v_k(f_{\alpha\bar\beta})+ v_k(e_{i\bar\beta})f_{\alpha\bar j}
+e_{i\bar\beta}v_k(f_{\alpha\bar j}))\ dv\\
=&\int_{X_s}((\Box+1)v_k(e_{i\bar j})e_{\alpha\bar\beta}+ e_{i\bar
j}v_k(f_{\alpha\bar\beta})+
(\Box+1)v_k(e_{i\bar\beta})e_{\alpha\bar j}
+e_{i\bar\beta}v_k(f_{\alpha\bar j}))\ dv\\
=&\int_{X_s}(v_k(f_{i\bar j})e_{\alpha\bar\beta}+ e_{i\bar
j}v_k(f_{\alpha\bar\beta})+ v_k(f_{i\bar\beta})e_{\alpha\bar j}
+e_{i\bar\beta}v_k(f_{\alpha\bar j}))\ dv\\
+&\int_{X_s}(\xi_k(e_{i\bar j})e_{\alpha\bar\beta}
+\xi_k(e_{i\bar\beta})e_{\alpha\bar j})\ dv\\
=&\int_{X_s}((L_k B_i \cdot \bar{B_j})e_{\alpha\bar\beta} +(L_k
B_\alpha \cdot \bar{B_\beta})e_{i\bar j} +(L_k B_i \cdot
\bar{B_\beta})e_{\alpha\bar j}
+(L_k B_\alpha \cdot \bar{B_j})e_{i\bar\beta})\ dv\\
+&\int_{X_s}((B_i \cdot L_k \bar{B_j})e_{\alpha\bar\beta}
+(B_\alpha \cdot L_k \bar{B_\beta})e_{i\bar j} +(B_i \cdot L_k
\bar{B_\beta})e_{\alpha\bar j}
+(B_\alpha \cdot L_k \bar{B_j})e_{i\bar\beta})\ dv\\
+&\int_{X_s}(\xi_k(e_{i\bar j})e_{\alpha\bar\beta}
+\xi_k(e_{i\bar\beta})e_{\alpha\bar j})\ dv.
\end{split}
\end{align}
Now we simplify the right hand side of \eqref{1st10}. Since
$B_1,\cdots,B_n$ is a basis of $T_s M_g$, we know that the first
line of the  right hand side of \eqref{1st10} is
\begin{align}\label{1st20}
\begin{split}
&\int_{X_s}((L_k B_i \cdot \bar{B_j})e_{\alpha\bar\beta} +(L_k
B_\alpha \cdot \bar{B_\beta})e_{i\bar j} +(L_k B_i \cdot
\bar{B_\beta})e_{\alpha\bar j}
+(L_k B_\alpha \cdot \bar{B_j})e_{i\bar\beta})\ dv\\
=&\int_{X_s}(L_k B_i \cdot (\bar{B_j}e_{\alpha\bar\beta}
+\bar{B_\beta}e_{\alpha\bar j}) +L_k B_\alpha
\cdot(\bar{B_j}e_{i\bar\beta}
+\bar{B_\beta}e_{i\bar j}))\ dv\\
=&h^{p\bar q}\int_{X_s}(L_k B_i \cdot \bar{B_q})\ dv
\int_{X_s}(B_p \cdot (\bar{B_j}e_{\alpha\bar\beta}
+\bar{B_\beta}e_{\alpha\bar j})\ dv\\
&+h^{p\bar q}\int_{X_s}(L_k B_\alpha \cdot \bar{B_q})\ dv
\int_{X_s}(B_p \cdot (\bar{B_j}e_{i\bar\beta}
+\bar{B_\beta}e_{i\bar j})\ dv\\
=&h^{p\bar q}\partial_k h_{i\bar q}R_{p\bar j\alpha\bar\beta}
+h^{p\bar q}\partial_k h_{\alpha\bar q}R_{i\bar j p\bar\beta}
=\Gamma_{ik}^p R_{p\bar j\alpha\bar\beta} +\Gamma_{\alpha k}^p
R_{i\bar j p\bar\beta}.
\end{split}
\end{align}
We deal with the second line of the right hand side of
\eqref{1st10} by using Lemma \ref{us1} and Lemma \ref{commu} to
get
\begin{eqnarray}\label{1st30}
B_i \cdot L_k \bar{B_j}=-A_i\pz(\lambda^{-1}\pz e_{k\bar j})
=\xi_i(e_{k\bar j}).
\end{eqnarray}
This implies
\begin{align}\label{1st40}
\begin{split}
&\int_{X_s}((B_i \cdot L_k \bar{B_j})e_{\alpha\bar\beta}
+(B_\alpha \cdot L_k \bar{B_\beta})e_{i\bar j} +(B_i \cdot L_k
\bar{B_\beta})e_{\alpha\bar j}
+(B_\alpha \cdot L_k \bar{B_j})e_{i\bar\beta})\ dv\\
=&\int_{X_s}(\xi_i(e_{k\bar j})e_{\alpha\bar\beta}
+\xi_\alpha(e_{k\bar\beta})e_{i\bar j}
+\xi_i(e_{k\bar\beta})e_{\alpha\bar j} +\xi_\alpha(e_{k\bar
j})e_{i\bar\beta})\ dv.
\end{split}
\end{align}
We also have
\begin{eqnarray}\label{1st50}
\partial_k \tau_{i\bar j}=h^{\alpha\bar\beta}\partial_k
R_{i\bar j\alpha\bar\beta}+\partial_k h^{\alpha\bar\beta} R_{i\bar
j\alpha\bar\beta}=h^{\alpha\bar\beta} (\partial_k R_{i\bar
j\alpha\bar\beta}- R_{i\bar j p\bar\beta}\Gamma_{k\alpha}^p).
\end{eqnarray}
By combining \eqref{1st10}, \eqref{1st20}, \eqref{1st40} and
\eqref{1st50}, together with the fact that $\xi_i$ is a real
symmetric operator and the definition of $\tau_{i\bar j}$, we have
proved this theorem.

\qed

To compute the second order derivative, we need to compute the
commutator of $\xi_k$ and $\bar{v_l}$. We have
\begin{lemma}\label{commu10}
For any smooth function $f \in C^\infty(X_s)$,
\begin{eqnarray}\label{klcomm}
\bar{v_l}(\xi_k f)-\xi_k(\bar{v_l}f)= \bar P(e_{k\bar
l})P(f)-2f_{k\bar l}\Box f +\lambda^{-1}\pz f_{k\bar l}\pzb f.
\end{eqnarray}
\end{lemma}
{\bf Proof.} We will fix local holomorphic coordinates and compute
locally. First we know that the commutator of $\bar{v_l}$ and
$\pz$ is
\begin{eqnarray}\label{kl10}
\bar{v_l}\pz-\pz\bar{v_l}=-\pz\bar{a_l}\pzb=-\bar{A_l}\pzb.
\end{eqnarray}
Similarly, the commutator of $\bar{v_l}$ and $\lambda^{-1}\pz$ is
\begin{eqnarray}\label{kl20}
\bar{v_l}(\lambda^{-1}\pz)-\lambda^{-1}\pz\bar{v_l}
=\bar{v_l}(\lambda^{-1})\pz+\lambda^{-1}
(\bar{v_l}\pz-\pz\bar{v_l})
=\lambda^{-1}\pzb\bar{a_l}\pz-\lambda^{-1}\bar{A_l}\pzb.
\end{eqnarray}
The above two formulae imply
\begin{align}\label{kl30}
\begin{split}
\bar{v_l}P-P\bar{v_l}=&-\bar{v_l}(\pz(\lambda^{-1}\pz))
+\pz(\lambda^{-1}\pz)\bar{v_l}\\
=&(\bar{A_l}\pzb-\pz\bar{v_l})(\lambda^{-1}\pz)
+\pz(\bar{v_l}(\lambda^{-1}\pz)
-\lambda^{-1}\pzb\bar{a_l}\pz+\lambda^{-1}\bar{A_l}\pzb)\\
=&\bar{A_l}\pzb(\lambda^{-1}\pz)-\pz
(\lambda^{-1}\pzb\bar{a_l}\pz)+\pz
(\lambda^{-1}\bar{A_l}\pzb)\\
=&-\lambda^{-2}\pzb\lambda\bar{A_l}\pz
+\lambda^{-1}\bar{A_l}\pz\pzb
+\lambda^{-2}\pz\lambda\pzb\bar{a_l}\pz
-\lambda^{-1}\pzb\bar{A_l}\pz
-\lambda^{-1}\pzb\bar{a_l}\pz\pz\\
&-\lambda^{-2}\pz\lambda\bar{A_l}\pzb
+\lambda^{-1}\pz\bar{A_l}\pzb +\lambda^{-1}\bar{A_l}\pz\pzb.
\end{split}
\end{align}
By using the harmonicity, we have $\pzb(\lambda\bar{A_l})=0$ which
implies $\pzb\bar{A_l}=-\lambda^{-1}\pzb\lambda \bar{A_l}$. By
plugging this into formula \eqref{kl30} we have
\begin{align}\label{kl40}
\begin{split}
\bar{v_l}P-P\bar{v_l}=&-2\bar{A_l}\Box
+\lambda^{-2}\pz\lambda\pzb\bar{a_l}\pz
-\lambda^{-1}\pzb\bar{a_l}\pz\pz
-\lambda^{-2}\pz\lambda\bar{A_l}\pzb
+\lambda^{-1}\pz\bar{A_l}\pzb\\
=&-2\bar{A_l}\Box+\pzb\bar{a_l}P
-\lambda^{-2}\pz\lambda\bar{A_l}\pzb
+\lambda^{-1}\pz\bar{A_l}\pzb.
\end{split}
\end{align}
Now, since $\xi_k=A_kP$, we have
\begin{align}\label{kl50}
\begin{split}
\bar{v_l}(\xi_k f)-\xi_k(\bar{v_l}f)=&
\bar{v_l}(A_k)P(f)+A_k(\bar{v_l}P(f)-P\bar{v_l}(f))\\
=&(\bar{v_l}(A_k)+A_k\pzb\bar{a_l})P(f) -2f_{k\bar l}\Box f
-\lambda^{-2}\pz\lambda A_k\bar{A_l}\pzb
+\lambda^{-1}A_k\pz\bar{A_l}\pzb.
\end{split}
\end{align}
From the proof of lemma \ref{us1} we know
$\bar{v_l}(A_k)+A_k\pzb\bar{a_l}=\bar{P}(e_{k\bar l})$. By using
the harmonicity we have $-\lambda^{-1}\pz\lambda A_k =\pz A_k$. So
from \eqref{kl50} we have
\begin{align}\label{kl60}
\begin{split}
\bar{v_l}(\xi_k f)-\xi_k(\bar{v_l}f)=& \bar{P}(e_{k\bar
l})P(f)-2f_{k\bar l}\Box f +\lambda^{-1}\pz A_k\bar{A_l}\pzb f
+\lambda^{-1}A_k\pz\bar{A_l}\pzb f\\
=&\bar{P}(e_{k\bar l})P(f)-2f_{k\bar l}\Box f +\lambda^{-1}\pz
f_{k\bar l}\pzb f.
\end{split}
\end{align}
This finishes the proof.

\qed

From the above lemma, it is convenient to define the commutator of
$\xi_k$ and $\bar{v_l}$ as an operator.
\begin{definition}\label{commu30}
For each $k,l$, we define the operator $Q_{k\bar l}$ which acts on
a function to produce another function by
\begin{eqnarray}\label{commu50}
Q_{k\bar l}(f)=\bar{P}(e_{k\bar l})P(f)-2f_{k\bar l}\Box f
+\lambda^{-1}\pz f_{k\bar l}\pzb f.
\end{eqnarray}
\end{definition}
Now we are ready to compute the curvature tensor of the Ricci
metric. The formula consists of four types of terms.
\begin{theorem}\label{riccicurv}
Let $s_1,\cdots,s_n$ be local holomorphic coordinates at $s \in
M_g$. Then at $s$, we have
\begin{align}\label{finalcurv}
\begin{split}
\widetilde{R}_{i\bar j k\bar l}=&h^{\alpha\bar\beta}
\left\{\sigma_1\sigma_2\int_{X_s}
\left\{(\Box+1)^{-1}(\xi_k(e_{i\bar j}))
\bar{\xi}_l(e_{\alpha\bar\beta})+ (\Box+1)^{-1} (\xi_k(e_{i\bar
j})) \bar{\xi}_\beta(e_{\alpha\bar l})
\right\}\ dv\right\}\\
&+h^{\alpha\bar\beta} \left\{\sigma_1\int_{X_s}Q_{k\bar
l}(e_{i\bar j}) e_{\alpha\bar\beta}\ dv
\right\}\\
&-\tau^{p\bar q}h^{\alpha\bar\beta}h^{\gamma\bar\delta}
\left\{\sigma_1\int_{X_s}\xi_k(e_{i\bar q}) e_{\alpha\bar\beta}\
dv\right\}\left\{ \widetilde\sigma_1\int_{X_s}\bar{\xi}_l(e_{p\bar
j})
e_{\gamma\bar\delta})\ dv\right\}\\
&+\tau_{p\bar j}h^{p\bar q}R_{i\bar q k\bar l}.
\end{split}
\end{align}
\end{theorem}
{\bf Proof.} By Lemma \ref{us2} we know that $L_kB_i$ is harmonic.
Since $B_1,\cdots,B_n$ is a basis of harmonic Beltrami
differentials, from the proof of Theorem \ref{wpcurv} we have
\begin{eqnarray}\label{linear}
L_kB_i=\Gamma_{ik}^s B_s.
\end{eqnarray}
We first compute $\partial_{\bar l} \int_{X_s}\xi_k(e_{i\bar
j})e_{\alpha\bar\beta}\ dv$. By Lemma \ref{commu} and Lemma
\ref{commu10} we have
\begin{align}\label{key100}
\begin{split}
\partial_{\bar l}
\int_{X_s}\xi_k(e_{i\bar j})e_{\alpha\bar\beta}\ dv
=&\int_{X_s}(\bar v_l(\xi_k(e_{i\bar j}))
e_{\alpha\bar\beta}+\xi_k(e_{i\bar j})\bar v_l
(e_{\alpha\bar\beta}))\ dv\\
=&\int_{X_s}(\xi_k(\bar v_l(e_{i\bar j}))
e_{\alpha\bar\beta}+\xi_k(e_{i\bar j})\bar v_l
(e_{\alpha\bar\beta})
+Q_{k\bar l}(e_{i\bar j})e_{\alpha\bar\beta})\ dv\\
=&\int_{X_s}(\xi_k(e_{\alpha\bar\beta})\bar v_l (e_{i\bar
j})+\xi_k(e_{i\bar j})\bar v_l (e_{\alpha\bar\beta})
+Q_{k\bar l}(e_{i\bar j})e_{\alpha\bar\beta})\ dv\\
=&\int_{X_s}(\Box+1)^{-1}(\xi_k(e_{\alpha\bar\beta}))
(\Box+1)(\bar v_l(e_{i\bar j}))\ dv\\
&+\int_{X_s}(\Box+1)^{-1}(\xi_k(e_{i\bar j})) (\Box+1)(\bar
v_l(e_{\alpha\bar\beta}))\ dv
+\int_{X_s}Q_{k\bar l}(e_{i\bar j})e_{\alpha\bar\beta}\ dv\\
=&\int_{X_s}(\Box+1)^{-1}(\xi_k(e_{\alpha\bar\beta}))
(\bar\xi_l(e_{i\bar j})+\bar v_l(f_{i\bar j}))\ dv\\
&+\int_{X_s}(\Box+1)^{-1}(\xi_k(e_{i\bar j}))
(\bar\xi_l(e_{\alpha\bar\beta})
+\bar v_l(f_{\alpha\bar\beta}))\ dv\\
&+\int_{X_s}Q_{k\bar l}(e_{i\bar j})e_{\alpha\bar\beta}\ dv\\
=&\int_{X_s}((\Box+1)^{-1}(\xi_k(e_{\alpha\bar\beta}))
\bar\xi_l(e_{i\bar j})+(\Box+1)^{-1}(\xi_k(e_{i\bar j}))
\bar\xi_l(e_{\alpha\bar\beta}))\ dv\\
&+\int_{X_s}(\Box+1)^{-1}(\xi_k(e_{\alpha\bar\beta}))
(\bar\xi_j(e_{i\bar l})+A_i\cdot L_{\bar l}\bar{A_j})\ dv\\
&+\int_{X_s}(\Box+1)^{-1}(\xi_k(e_{i\bar j}))
(\bar\xi_\beta(e_{\alpha\bar l})
+A_\alpha\cdot L_{\bar l}\bar{A_\beta})\ dv\\
&+\int_{X_s}Q_{k\bar l}(e_{i\bar j})e_{\alpha\bar\beta}\ dv.
\end{split}
\end{align}
Now by using \eqref{linear} we have
\begin{align}\label{key110}
\begin{split}
&\int_{X_s}((\Box+1)^{-1}(\xi_k(e_{\alpha\bar\beta})) (A_i\cdot
L_{\bar l}\bar{A_j})+ (\Box+1)^{-1}(\xi_k(e_{i\bar j}))
(A_\alpha\cdot L_{\bar l}\bar{A_\beta}))\ dv\\
=&\int_{X_s}((\Box+1)^{-1}(\xi_k(e_{\alpha\bar\beta}))
(\bar{\Gamma_{jl}^t}A_i\cdot \bar{A_t})+
(\Box+1)^{-1}(\xi_k(e_{i\bar j}))
(\bar{\Gamma_{\beta l}^t}A_\alpha\cdot \bar{A_t}))\ dv\\
=&\bar{\Gamma_{jl}^t}\int_{X_s}\xi_k(e_{\alpha\bar\beta})
(\Box+1)^{-1}(A_i\cdot \bar{A_t})\ dv +\bar{\Gamma_{\beta
l}^t}\int_{X_s}\xi_k(e_{i\bar j})
(\Box+1)^{-1}(A_\alpha\cdot \bar{A_t})\ dv\\
=&\bar{\Gamma_{jl}^t}\int_{X_s}\xi_k(e_{\alpha\bar\beta}) e_{i\bar
t}\ dv +\bar{\Gamma_{\beta l}^t}\int_{X_s}\xi_k(e_{i\bar j})
e_{\alpha\bar t}\ dv.
\end{split}
\end{align}
By combining \eqref{key100} and \eqref{key110} we have
\begin{align}\label{key120}
\begin{split}
\partial_{\bar l}
\int_{X_s}\xi_k(e_{i\bar j})e_{\alpha\bar\beta}\ dv
=&\int_{X_s}(\Box+1)^{-1}(\xi_k(e_{i\bar j}))
(\bar\xi_l(e_{\alpha\bar\beta})
+\bar\xi_\beta(e_{\alpha\bar l}))\ dv\\
&+\int_{X_s}(\Box+1)^{-1}(\xi_k(e_{\alpha\bar\beta}))
(\bar\xi_l(e_{i\bar j})
+\bar\xi_j(e_{i\bar l}))\ dv\\
&+\bar{\Gamma_{jl}^t}\int_{X_s}\xi_k(e_{\alpha\bar\beta}) e_{i\bar
t}\ dv +\bar{\Gamma_{\beta l}^t}\int_{X_s}\xi_k(e_{i\bar j})
e_{\alpha\bar t}\ dv\\
&+\int_{X_s}Q_{k\bar l}(e_{i\bar j})e_{\alpha\bar\beta}\ dv.
\end{split}
\end{align}
We also have
\begin{align}\label{key130}
\begin{split}
\partial_{\bar l}\Gamma_{ik}^p=&
\partial_{\bar l}(h^{p\bar q}\partial_k h_{i\bar q})
=-h^{p\bar\beta}h^{\alpha\bar q}\partial_{\bar l}
h_{\alpha\bar\beta}\partial_k h_{i\bar q}+h^{p\bar q}
\partial_{\bar l}\partial_k h_{i\bar q}\\
=&h^{p\bar q}(\partial_{\bar l}\partial_k h_{i\bar q}
-h^{\alpha\bar\beta}\partial_{\bar l} h_{\alpha\bar q}\partial_k
h_{i\bar\beta}) =h^{p\bar q}R_{i\bar q k\bar l}.
\end{split}
\end{align}
From Theorem \ref{1stderiv}, formula \eqref{key120} and
\eqref{key130} we derive
\begin{align}\label{key140}
\begin{split}
\partial_{\bar l}\partial_k \tau_{i\bar j}=&
(\partial_{\bar l}h^{\alpha\bar\beta}) \left\{\sigma_1\int_{X_s}
\xi_k(e_{i\bar j})e_{\alpha\bar\beta}\ dv \right\}
+h^{\alpha\bar\beta}\left\{\sigma_1
\partial_{\bar l}\int_{X_s}
\xi_k(e_{i\bar j})e_{\alpha\bar\beta}\ dv \right\}\\
&+h^{\gamma\bar\delta}\left\{\widetilde\sigma_1\int_{X_s}
\bar\xi_l(e_{p\bar j})e_{\gamma\bar\delta}\ dv \right\}
\Gamma_{ik}^p+\tau_{p\bar q}\Gamma_{ik}^p \bar{\Gamma_{jl}^q}
+\tau_{p\bar j}h^{p\bar q}R_{i\bar q k\bar l}\\
=&-h^{\alpha\bar t}\bar{\Gamma_{lt}^\beta}
\left\{\sigma_1\int_{X_s}
\xi_k(e_{i\bar j})e_{\alpha\bar\beta}\ dv \right\}\\
&+h^{\alpha\bar\beta}\left\{\sigma_1\sigma_2\int_{X_s}
(\Box+1)^{-1}(\xi_k(e_{i\bar j})) (\bar\xi_l(e_{\alpha\bar\beta})
+\bar\xi_\beta(e_{\alpha\bar l}))\ dv\right\}\\
&+h^{\alpha\bar\beta} \left\{\sigma_1\int_{X_s}Q_{k\bar
l}(e_{i\bar j}) e_{\alpha\bar\beta}\ dv\right\}
+h^{\alpha\bar\beta}\bar{\Gamma_{jl}^t} \left\{\sigma_1\int_{X_s}
\xi_k(e_{i\bar t})e_{\alpha\bar\beta}\ dv\right\}\\
&+h^{\alpha\bar\beta}\bar{\Gamma_{\beta l}^t}
\left\{\sigma_1\int_{X_s} \xi_k(e_{i\bar j})e_{\alpha\bar t}\
dv\right\}
+h^{\gamma\bar\delta}\left\{\widetilde\sigma_1\int_{X_s}
\bar\xi_l(e_{p\bar j})e_{\gamma\bar\delta}\ dv \right\}
\Gamma_{ik}^p\\
&+\tau_{p\bar q}\Gamma_{ik}^p \bar{\Gamma_{jl}^q} +\tau_{p\bar
j}h^{p\bar q}R_{i\bar q k\bar l}.
\end{split}
\end{align}
Now from the above formula, by using Theorem \ref{1stderiv} we can
easily check the formula \eqref{finalcurv}.

\qed

The curvature formula of the Ricci metric would be simpler if we
have used the normal coordinates. However, when we estimate the
asymptotic behavior of the curvature, it is hard to describe the
normal coordinates near the boundary points. Thus we will use this
general formula directly in our computations. The estimates are
quite subtle.

\section{The asymptotics of the Ricci metric and its curvatures}\label{sec4}
From formula \eqref{WPcurv1} we can easily see the sign of the
curvature of the Weil-Petersson metric directly. However, the sign
of the curvature of the Ricci metric cannot be derived from
formula \eqref{finalcurv}. In this section, we estimate the
asymptotics of the Ricci metric and its curvatures. We first
describe the local pinching coordinates near the boundary of the
moduli space due to the plumbing construction of Wolpert. Then we
use Masur's construction of the holomorphic quadratic
differentials to estimate the harmonic Beltrami differentials.
Finally, we construct $\widetilde e_{i\bar j}$ which is an
approximation of $e_{i\bar j}$. By doing this we avoid the
estimates of the Green function of $\Box+1$ on the Riemann
surfaces.

Let $\mathcal M_g$ be the moduli space of Riemann surfaces of
genus $g \geq 2$ and let $\bar{\mathcal M}_g$ be its
Deligne-Mumford compactification \cite{dm1}. Each point $y \in
\bar{\mathcal M}_g \setminus \mathcal M_g$ corresponds to a stable
nodal surface $X_y$. A point $p \in X_y$ is a node if there is a
neighborhood of $p$ which is isometric to the germ $\{ (u,v)\mid
uv=0,\ |u|,|v|<1 \} \subset \mathbb{C}^2$.

We first recall the rs-coordinate on a Riemann surface defined by
Wolpert in \cite{wol1}. There are two cases: the puncture case and
the short geodesic case. For the puncture case, we have a nodal
surface $X$ and a node $p\in X$. Let $a,b$ be two punctures which
are glued together to form $p$.
\begin{definition}
A local coordinate chart $(U,u)$ near $a$ is called rs-coordinate
if $u(a)=0$ where $u$ maps $U$ to the punctured disc $0<|u|<c$
with $c>0$, and the restriction to $U$ of the K\"ahler-Einstein
metric on $X$ can be written as $\frac{1}{2|u|^2(\log |u|)^2}
|du|^2$. The rs-coordinate $(V,v)$ near $b$ is defined in a
similar way.
\end{definition}
For the short geodesic case, we have a closed surface $X$, a
closed geodesic $\gamma \subset X$ with length $l <c_\ast$ where
$c_\ast$ is the collar constant.
\begin{definition}
A local coordinate chart $(U,z)$ is called rs-coordinate at
$\gamma$ if $\gamma \subset U$ where $z$ maps $U$ to the annulus
$c^{-1}|t|^{\frac{1}{2}}<|z| <c|t|^{\frac{1}{2}}$, and the
K\"ahler-Einstein metric on $X$ can be written as
$\frac{1}{2}(\frac{\pi}{\log |t|}\frac{1}{|z|}\csc \frac{\pi\log
|z|}{\log |t|})^2 |dz|^2$.
\end{definition}
\begin{rem}
We put the factor $\frac{1}{2}$ in the above two definitions to
normalize metrics such that \eqref{ke100} hold.
\end{rem}

By Keen's collar theorem \cite{ke1}, we have the following lemma:
\begin{lemma}\label{gcollar}
Let $X$ be a closed surface and let $\gamma$ be a closed geodesic
on $X$ such that the length $l$ of $\gamma$ satisfies $l <c_\ast$.
Then there is a collar $\Omega$ on $X$ with holomorphic coordinate
$z$ defined on $\Omega$ such that
\begin{enumerate}
\item $z$ maps $\Omega$ to the annulus
$\frac{1}{c}e^{-\frac{2\pi^2}{l}}<|z|<c$ for $c>0$; \item the
K\"ahler-Einstein metric on $X$ restricted to $\Omega$ is given by
\begin{eqnarray}\label{precmetric}
(\frac{1}{2}u^2 r^{-2}\csc^2\tau) |dz|^2
\end{eqnarray}
where $u=\frac{l}{2\pi}$, $r=|z|$ and $\tau=u\log r$; \item the
geodesic $\gamma$ is given by the equation $|z|=
e^{-\frac{\pi^2}{l}}$.
\end{enumerate}
We call such a collar $\Omega$ a genuine collar.
\end{lemma}
We notice that the constant $c$ in the above lemma has a lower
bound such that the area of $\Omega$ is bounded from below. Also,
the coordinate $z$ in the above lemma is rs-coordinate. In the
following, we will keep using the above notations $u$, $r$ and
$\tau$.

Now we describe the local manifold cover of $\bar{\mathcal M}_g$
near the boundary. We take the construction of Wolpert
\cite{wol1}. Let $X_{0,0}$ be a nodal surface corresponding to a
codimension $m$ boundary point. $X_{0,0}$ have $m$ nodes
$p_1,\cdots,p_m$. $X_0=X_{0,0}\setminus \{ p_1,\cdots,p_m \}$ is a
union of punctured Riemann surfaces. Fix the rs-coordinate charts
$(U_i,\eta_i)$ and $(V_i,\zeta_i)$ at $p_i$ for $i=1,\cdots,m$
such that all the $U_i$ and $V_i$ are mutually disjoint. Now pick
an open set $U_0 \subset X_0$ such that the intersection of each
connected component of $X_0$ and $U_0$ is a nonempty relatively
compact set and the intersection $U_0 \cap (U_i\cup V_i)$ is empty
for all $i$. Now pick Beltrami differentials
$\nu_{m+1},\cdots,\nu_{n}$ which are supported in $U_0$ and span
the tangent space at $X_0$ of the deformation space of $X_0$. For
$s=(s_{m+1},\cdots,s_n)$, let $\nu(s)=\sum_{i=m+1}^n s_i\nu_i$. We
assume $|s|=(\sum |s_i|^2)^{\frac{1}{2}}$  small enough such that
$|\nu(s)|<1$. The nodal surface $X_{0,s}$ is obtained by solving
the Beltrami equation $\bar\partial w=\nu(s)\partial w$. Since
$\nu(s)$ is supported in $U_0$, $(U_i,\eta_i)$ and $(V_i,\zeta_i)$
are still holomorphic coordinates on $X_{0,s}$. Note that they are
no longer rs-coordinates. By the theory of Alhfors and Bers
\cite{ab1} and Wolpert \cite{wol1} we can assume that there are
constants $\delta,c>0$ such that when $|s|<\delta$, $\eta_i$ and
$\zeta_i$ are holomorphic coordinates on $X_{0,s}$ with
$0<|\eta_i|<c$ and $0<|\zeta_i|<c$. Now we assume
$t=(t_1,\cdots,t_m)$ has small norm. We do the plumbing
construction on $X_{0,s}$ to obtain $X_{t,s}$. We remove from
$X_{0,s}$ the discs $0<|\eta_i|\leq \frac{|t_i|}{c}$ and
$0<|\zeta_i|\leq \frac{|t_i|}{c}$ for each $i=1,\cdots,m$, and
identify $\frac{|t_i|}{c}<|\eta_i|< c$ with
$\frac{|t_i|}{c}<|\zeta_i|< c$ by the rule $\eta_i \zeta_i=t_i$.
This defines the surface $X_{t,s}$. The tuple
$(t_1,\cdots,t_m,s_{m+1},\cdots,s_n)$ are the local pinching
coordinates for the manifold cover of $\bar{\mathcal M}_g$. We
call the coordinates $\eta_i$ (or $\zeta_i$) the plumbing
coordinates on $X_{t,s}$ and the collar defined by
$\frac{|t_i|}{c}<|\eta_i|< c$ the plumbing collar.
\begin{rem}
From the estimate of Wolpert \cite{wol2}, \cite{wol1} on the
length of short geodesic, we have $u_i=\frac{l_i}{2\pi}\sim
-\frac{\pi}{\log|t_i|}$.
\end{rem}

We also need the following version of the Schauder estimate proved
by Wolpert \cite{wol1}.
\begin{theorem}\label{laesti}
Let $X$ be a closed Riemann surface equipped with the unique
K\"ahler-Einstein metric. Let $f$ and $g$ be smooth functions on
$X$ such that $(\Box+1)g=f$. Then for any integer $k\geq 0$, there
is a constant $c_k$ such that $\Vert g\Vert _{k+1} \leq c_k \Vert
f\Vert _{k}$ where the norm is defined by \eqref{norm}.
\end{theorem}

Now we estimate the asymptotics of the Ricci metric in the
pinching coordinates. We will use the following notations. Let
$(t,s)=(t_1,\cdots,t_m,s_{m+1},\cdots,s_n)$ be the pinching
coordinates near $X_{0,0}$. For $|(t,s)|<\delta$, let $\Omega^j_c$
be the $j$-th genuine collar on $X_{t,s}$ which contains a short
geodesic $\gamma_j$ with length $l_j$. Let $u_j=\frac{l_j}{2\pi}$,
$u_0=\sum_{j=1}^m u_j+\sum_{j=m+1}^n |s_j|$, $r_j=|z_j|$ and
$\tau_j=u_j\log r_j$ where $z_j$ is the properly normalized
rs-coordinate on $\Omega^j_c$ such that
\[
\Omega^j_c=\{ z_j\mid c^{-1}e^{-\frac{2\pi^2}{l_j}}<|z_j|<c \}.
\]
From the above argument, we know that the K\"ahler-Einstein metric
$\lambda$ on $X_{t,s}$ restrict to the collar $\Omega^j_c$ is
given by
\begin{eqnarray}\label{jmetric}
\lambda=\frac{1}{2}u_j^2 r_j^{-2}\csc^2\tau_j .
\end{eqnarray}
For convenience, we let $\Omega_c=\cup_{j=1}^m \Omega^j_c$ and
$R_c=X_{t,s}\setminus \Omega_c$. In the following, we may change
the constant $c$ finitely many times, clearly this will not affect
the estimates.

To estimate the curvature of the Ricci metric, the first step is
to find all the harmonic Beltrami differentials $B_1,\cdots,B_n$
which correspond to the tangent vectors $\frac{\partial}{\partial
t_1},\cdots, \frac{\partial}{\partial s_n}$. In \cite{ma1}, Masur
constructed $3g-3$ regular holomorphic quadratic differentials
$\psi_1,\cdots,\psi_n$ on the plumbing collars by using the
plumbing coordinate $\eta_j$. These quadratic differentials
correspond to the cotangent vectors $dt_1,\cdots,ds_n$.

However, it is more convenient to estimate the curvature if we use
the rs-coordinate on $X_{t,s}$ since we have the accurate form of
the K\"ahler-Einstein metric $\lambda$ in this coordinate. In
\cite{tr1}, Trapani used the graft metric constructed by Wolpert
\cite{wol1} to estimate the difference between the plumbing
coordinate and rs-coordinate and gave the holomorphic quadratic
differentials constructed by Masur in the rs-coordinate. We
collect Trapani's results (Lemma 6.2-6.5, \cite{tr1}) in the
following theorem:
\begin{theorem}\label{imp}
Let $(t,s)$ be the pinching coordinates on $\bar{\mathcal M}_g$
near $X_{0,0}$ which corresponds to a codimension $m$ boundary
point of $\bar{\mathcal M}_g$. Then there exist constants
$M,\delta>0$ and $1>c>0$ such that if $|(t,s)|<\delta$, then the
$j$-th plumbing collar on $X_{t,s}$ contains the genuine collar
$\Omega^j_c$. Furthermore, one can choose rs-coordinate $z_j$ on
the collar $\Omega_c^j$ properly such that the holomorphic
quadratic differentials $\psi_1,\cdots,\psi_n$ corresponding to
the cotangent vectors $dt_1,\cdots,ds_n$ have the form
$\psi_i=\phi_i(z_j)dz_j^2$ on the genuine collar $\Omega^j_c$ for
$1\leq j \leq m$, where
\begin{enumerate}
\item $\phi_i(z_j)=\frac{1}{z_j^2}(q_i^j(z_j)+\beta_i^j)$ if
$i\geq m+1$; \item
$\phi_i(z_j)=(-\frac{t_j}{\pi})\frac{1}{z_j^2}(q_j(z_j)+\beta_j)$
if $i=j$; \item $\phi_i(z_j)=(-\frac{t_i}{\pi})
\frac{1}{z_j^2}(q_i^j(z_j)+\beta_i^j)$ if $1\leq i \leq m$ and
$i\ne j$.
\end{enumerate}
Here $\beta_i^j$ and $\beta_j$ are functions of $(t,s)$, $q_i^j$
and $q_j$ are functions of $(t,s,z_j)$ given by
\[
q_i^j(z_j)=\sum_{k<0}\alpha_{ik}^j(t,s)t_j^{-k}z_j^k
+\sum_{k>0}\alpha_{ik}^j(t,s)z_j^k
\]
and
\[
q_j(z_j)=\sum_{k<0}\alpha_{jk}(t,s)t_j^{-k}z_j^k
+\sum_{k>0}\alpha_{jk}(t,s)z_j^k
\]
such that
\begin{enumerate}
\item $\sum_{k<0}|\alpha_{ik}^j|c^{-k}\leq M$ and
$\sum_{k>0}|\alpha_{ik}^j|c^{k}\leq M$ if $i\ne j$; \item
$\sum_{k<0}|\alpha_{jk}|c^{-k}\leq M$ and
$\sum_{k>0}|\alpha_{jk}|c^{k}\leq M$; \item
$|\beta_i^j|=O(|t_j|^{\frac{1}{2}-\epsilon})$ with
$\epsilon<\frac{1}{2}$ if $i\ne j$; \item $|\beta_j|=(1+O(u_0))$.
\end{enumerate}
\end{theorem}

An immediate consequence of the above theorem is the following
refined version of Masur's estimates of the Weil-Petersson metric.
In the following, we will fix $(t,s)$ with small norm and let
$X=X_{t,s}$.
\begin{cor}\label{wpasymp}
Let $(t,s)$ be the pinching coordinates. Then
\begin{enumerate}
\item $h^{i\bar i}=2u_i^{-3}|t_i|^2(1+O(u_0))$ and $h_{i\bar i}
=\frac{1}{2}\frac{u_i^{3}}{|t_i|^2}(1+O(u_0))$ for $1\leq i\leq
m$; \item $h^{i\bar j}=O(|t_it_j|)$ and $h_{i\bar
j}=O(\frac{u_i^3u_j^3}{|t_it_j|})$, if $1\leq i,j \leq m$ and
$i\ne j$; \item $h^{i\bar j}=O(1)$ and $h_{i\bar j}=O(1)$, if
$m+1\leq i,j \leq n$; \item $h^{i\bar j}=O(|t_i|)$ and $h_{i\bar
j}=O(\frac{u_i^3}{|t_i|})$ if $i\leq m < j$ or $j \leq m<i$.
\end{enumerate}
\end{cor}
{\bf Proof.} We need the following simple calculus results:
\begin{eqnarray}\label{calculus10}
\int_{c^{-1}e^{-\frac{2\pi^2}{l_j}}}^c \frac{1}{r_j}\sin^2\tau_j\
dr_j =u_j^{-1}(\frac{\pi}{2}+O(u_j)).
\end{eqnarray}
For any $k \geq 1$,
\begin{eqnarray}\label{calculus20}
\int_{c^{-1}e^{-\frac{2\pi^2}{l_j}}}^c r_j^{k-1}\sin^2\tau_j\ dr_j
=O(u_j^2)c^k
\end{eqnarray}
and for $k\leq -1$,
\begin{eqnarray}\label{calculus30}
\int_{c^{-1}e^{-\frac{2\pi^2}{l_j}}}^c r_j^{k-1}\sin^2\tau_j\ dr_j
=O(u_j^2)c^{-k}\bigg ( e^{-\frac{2\pi^2}{l_j}}\bigg )^k.
\end{eqnarray}
On the collar $\Omega_c^j$, the metric $\lambda$ is given by
\eqref{jmetric}. $h^{i\bar j}$ is given by the formula
\[
h^{i\bar j}=\int_{X} \psi_i\bar{\psi_j}\lambda^{-2}dv.
\]
By using the above calculus facts, we can compute the above
integral on the collars. The bound on $R_c$ was calculated in
\cite{ma1}. A simple computation shows that the first part of all
of the above claims hold. The second parts of these claims can be
obtained by inverting the matrix $(h^{i\bar j})$ together with
Masur's result on the nondegenerate extension of the submatrix
$(h^{i\bar j})_{i,j
>m}$.This finishes the proof.

\qed

Now we are ready to compute the harmonic Beltrami differentials
$B_i=A_i\pz\otimes d\bar z$.

\begin{lemma}\label{aj10}
 For $c$ small, on the genuine collar $\Omega_c^j$, the coefficient functions
$A_i$ of the harmonic Beltrami differentials have the form:
\begin{enumerate}
\item $A_i=\frac{z_j}{\bar{z_j}}\sin^2\tau_j \big (
\bar{p_i^j(z_j)}+\bar{b_i^j}\big )$ if $i\ne j$; \item
$A_j=\frac{z_j}{\bar{z_j}}\sin^2\tau_j(\bar{p_j(z_j)}+\bar{b_j})$
\end{enumerate}
where
\begin{enumerate}
\item $p_i^j(z_j)=\sum_{k\leq -1}a_{ik}^j\rho_j^{-k}z_j^k
+\sum_{k\geq 1}a_{ik}^jz_j^k$ if $i\ne j$; \item
$p_j(z_j)=\sum_{k\leq -1}a_{jk}\rho_j^{-k}z_j^k +\sum_{k\geq
1}a_{jk}z_j^k$.
\end{enumerate}
In the above expressions, $\rho_j=e^{-\frac{2\pi^2}{l_j}}$ and the
coefficients satisfy the following conditions:
\begin{enumerate}
\item $\sum_{k\leq -1}|a_{ik}^j|c^{-k}=O(u_j^{-2})$ and
$\sum_{k\geq 1}|a_{ik}^j|c^{k}=O(u_j^{-2})$ if $i\geq m+1$; \item
$\sum_{k\leq -1}|a_{ik}^j|c^{-k}=O(u_j^{-2}) O\big (
\frac{u_i^3}{|t_i|}\big )$ and $\sum_{k\geq
1}|a_{ik}^j|c^{k}=O(u_j^{-2}) O\big ( \frac{u_i^3}{|t_i|}\big )$
if $i\leq m$ and $i\ne j$; \item $\sum_{k\leq
-1}|a_{jk}|c^{-k}=O(\frac{u_j}{|t_j|})$ and $\sum_{k\geq
1}|a_{jk}|c^{k}=O(\frac{u_j}{|t_j|})$; \item $|b_i^j|=O(u_j)$ if
$i\geq m+1$; \item $|b_i^j|=O(u_j)O\big ( \frac{u_i^3}{|t_i|}\big
)$ if $i\leq m$ and $i\ne j$; \item
$b_j=-\frac{u_j}{\pi\bar{t_j}}(1+O(u_0))$.
\end{enumerate}
\end{lemma}
{\bf Proof.} The duality between the harmonic Beltrami
differentials and the holomorphic quadratic differentials is given
by
\begin{eqnarray}\label{dual100}
B_i=\lambda^{-1}\sum_{l=1}^n h_{i\bar l}\bar{\psi_l}
\end{eqnarray}
which implies $A_i=\lambda^{-1}\sum_{l=1}^n h_{i\bar
l}\bar{\phi_l}$. Now by Wolpert's estimate on the length of the
short geodesic $\gamma_j$ in \cite{wol1} we have
$l_j=-\frac{2\pi^2}{\log|t_j|} (1+O(u_j))$. This implies there is
a constant $0<\mu<1$ such that $\mu|t_j|<\rho_j<\mu^{-1}|t_j|$.
The lemma follows from equation \eqref{dual100} by replacing $c$
by $\mu c$, a simple computation together with Theorem \ref{imp}
and Corollary \ref{wpasymp}.

\qed

To estimate the curvature of the Ricci metric, we need to estimate
the asymptotics of the Ricci metric by using Theorem \ref{wpcurv}.
So we need the following estimates on the norms of the harmonic
Beltrami differentials.

\begin{lemma}\label{esafij1}
Let $\Vert \cdot \Vert_k$ be the norm as defined in Definition
\ref{norm}. We have
\begin{enumerate}
\item $\Vert A_i\Vert_{0,\Omega_c^i}=O\big ( \frac{u_i}{|t_i|}\big
)$ and $\Vert A_i\Vert_{0,X\setminus\Omega_c^i} =O\big (
\frac{u_i^3}{|t_i|}\big )$, if $i\leq m$; \item $\Vert
A_i\Vert_{0}=O(1)$, if $i \geq m+1$; \item $\Vert f_{i\bar
i}\Vert_{0,\Omega_c^i} =O\big ( \frac{u_i^2}{|t_i|^2}\big )$ and
$\Vert f_{i\bar i}\Vert_{0,X\setminus\Omega_c^i} =O\big (
\frac{u_i^6}{|t_i|^2}\big )$, if $i\leq m$; \item $\Vert f_{i\bar
j}\Vert_{0}=O(1)$, if $i,j \geq m+1$; \item  $\Vert f_{i\bar
j}\Vert_{0,\Omega_c^i} =O\big ( \frac{u_i u_j^3}{|t_i t_j|}\big )$
and $\Vert f_{i\bar j}\Vert_{0,\Omega_c^j} =O\big ( \frac{u_i^3
u_j}{|t_i t_j|}\big )$ and $\Vert f_{i\bar
j}\Vert_{0,X\setminus(\Omega_c^i\cup\Omega_c^j)}
=O\big ( \frac{u_i^3 u_j^3}{|t_i t_j|}\big )$\\
if $i,j\leq m$ and $i\ne j$; \item  $\Vert f_{i\bar
j}\Vert_{0,\Omega_c^i} =O\big ( \frac{u_i}{|t_i|}\big )$ and
$\Vert f_{i\bar j}\Vert_{0,X\setminus\Omega_c^i} =O\big (
\frac{u_i^3}{|t_i|}\big )$, if $i\leq m$ and $j\geq m+1$; \item
$|f_{i\bar j}|_{L^1}=O(1)$, if $i,j \geq m+1$; \item $|f_{i\bar
j}|_{L^1}=O(\frac{u_i^3}{|t_i|})$, if $i\leq m$ and $j \geq m+1$;
\item $|f_{i\bar j}|_{L^1}=O(\frac{u_i^3 u_j^3}{|t_i t_j|})$, if
$i,j\leq m$ and $i\ne j$.
\end{enumerate}
\end{lemma}
{\bf Proof.} We choose $c$ small enough such that for each $1\leq
j\leq m$,
\[
\tan(u_j\log c)<-10 u_j
\]
when $|(t,s)|<\delta$. A simple computation shows that, when
$1\leq p\leq 10$, on the collar $\Omega_c^j$ we have
\[
|r_j^k\sin^p\tau_j|\leq c^k|\log c|^p u_j^p
\]
if $k \geq 1$, and
\[
|r_j^k\sin^p\tau_j|\leq c^{-k}|\log c|^p\rho_j^k u_j^p
\]
if $k \leq -1$.

To prove the first claim, note that on $\Omega_c^i$ we have
\begin{align*}
\begin{split}
|A_i|=&\left | \frac{z_i}{\bar{z_i}}\right |
|\sin^2\tau_i(\bar{p_i}+\bar{b_i})| \leq \sum_{k\leq
-1}|a_{ik}|\rho_i^{-k}r_i^k\sin^2\tau_i
+\sum_{k\geq 1}|a_{ik}|r_i^k\sin^2\tau_i +|b_j|\\
\leq & (\log c)^2 u_i^2\big (\sum_{k\leq -1}|a_{ik}|c^{-k}
+\sum_{k\geq 1}|a_{ik}|c^k \big ) +|b_j|\\
=& O(u_i^2)O\big ( \frac{u_i}{|t_i|}\big )+ O(u_i^2)O\big (
\frac{u_i}{|t_i|}\big )+O\big ( \frac{u_i}{|t_i|}\big ) =O\big (
\frac{u_i}{|t_i|}\big ).
\end{split}
\end{align*}
Similarly, on $\Omega_c^j$ with $j\ne i$, we have $|A_i|=O\big (
\frac{u_i^3}{|t_i|}\big )$. Also, on $R_c$ we have $|A_i|=O\big (
\frac{u_i^3}{|t_i|}\big )$ by the work of Masur \cite{ma1},
equation \eqref{dual100} together with Theorem \ref{imp} and
Corollary \ref{wpasymp}. This finishes the proof of the first
claim.

The second claim can be proved in a similar way. Claim (3)-(6)
follow from the first and second claims by using the fact that
$f_{i\bar j}=A_i\bar{A_j}$. Claim (7) follows from claim (4) and
the fact that the area of $X$ is a fixed positive constant using
the Gauss-Bonnet theorem.

Now we prove claim (9). On $\Omega_c^i$, by using a similar
estimate as above, we have
\begin{align*}
\begin{split}
|f_{i\bar j}|=&|\sin^4\tau_i(\bar{p_i}+\bar{b_i})(p_j^i+b_j^i)|
\leq |\sin^4\tau_i\bar{p_i}p_j^i|+|\sin^4\tau_i\bar{b_i}p_j^i|
+|\sin^4\tau_i\bar{p_i}b_j^i|+|\sin^4\tau_i\bar{b_i}b_j^i|\\
\leq & O\big (\frac{u_i^3 u_j^3}{|t_i t_j|}\big )
+|\sin^4\tau_i\bar{b_i}b_j^i| =O\big (\frac{u_i^3 u_j^3}{|t_i
t_j|}\big ) +O\big (\frac{u_i^2 u_j^3}{|t_i t_j|}\big
)\sin^4\tau_i.
\end{split}
\end{align*}
So
\[
|f_{i\bar j}|_{L^1(\Omega_c^i)}\leq \int_{\Omega_c^i} \bigg (O\big
(\frac{u_i^3 u_j^3}{|t_i t_j|}\big ) +O\big (\frac{u_i^2
u_j^3}{|t_i t_j|}\big )\sin^4\tau_i\bigg ) dv =O\big (\frac{u_i^3
u_j^3}{|t_i t_j|}\big ).
\]
Similarly, $|f_{i\bar j}|_{L^1(\Omega_c^j)}\leq O\big (\frac{u_i^3
u_j^3}{|t_i t_j|}\big )$. The estimate $|f_{i\bar
j}|_{L^1(X\setminus(\Omega_c^i\cup\Omega_c^j))}= O\big
(\frac{u_i^3 u_j^3}{|t_i t_j|}\big )$ follows from claim (5). This
proves claim (9). Similarly we can prove claim (8).

\qed

In the following, we will denote the operator $(\Box+1)^{-1}$ by
$T$. We then have the following estimates about $L^2$ norms:
\begin{lemma}\label{spec10}
Let $f\in C^{\infty}(X,\mathbb C)$. Then we have
\begin{eqnarray}\label{lt10}
\int_X |Tf|^2\ dv \leq \int_X Tf\cdot\bar{f}\ dv \leq \int_X
|f|^2\ dv.
\end{eqnarray}
\end{lemma}
{\bf Proof.} This lemma is a simple application of the spectral
decomposition of the operator $(\Box+1)$ and the fact that all
eigenvalues of this operator are greater than or equal to $1$. One
can also prove it directly by using integration by part.

\qed

To estimate the Ricci metric, we also need to estimate the
functions $e_{i\bar j}$. We localize these functions on the
collars by constructing the following approximation functions.

Pick a positive constant $c_1<c$ and define the cut-off function
$\eta\in C^{\infty}(\mathbb R, [0,1])$ by
\begin{eqnarray}\label{cutoff}
\begin{cases}
\eta(x)=1, & x\leq \log c_1;\\
\eta(x)=0, & x\geq \log c;\\
0<\eta(x)<1, & \log c_1<x<\log c.
\end{cases}
\end{eqnarray}
It is clear that the derivatives of $\eta$ are bounded by
constants which only depend on $c$ and $c_1$. Let
$\widetilde{e_{i\bar j}}(z)$ be the function on $X$ defined in the
following way where $z$ is taken to be $z_i$ on the collar
$\Omega_c^i$:
\begin{enumerate}
\item if $i\leq m$ and $j \geq m+1$, then
\begin{eqnarray*}
\widetilde{e_{i\bar j}}(z)=
\begin{cases}
\frac{1}{2}\sin^2\tau_i\bar{b_i}b_j^i, & z \in \Omega_{c_1}^i;\\
(\frac{1}{2}\sin^2\tau_i\bar{b_i}b_j^i)\eta(\log r_i), & z \in
\Omega_c^i \text{ and }
c_1<r_i<c;\\
(\frac{1}{2}\sin^2\tau_i\bar{b_i}b_j^i)\eta(\log\rho_i-\log r_i),
& z \in \Omega_c^i \text{ and }
c^{-1}\rho_i<r_i<c_1^{-1}\rho_i;\\
0, & z \in X\setminus\Omega_c^i;\\
\end{cases}
\end{eqnarray*}
\item if $i,j\leq m$ and $i\ne j$, then
\begin{eqnarray*}
\widetilde{e_{i\bar j}}(z)=
\begin{cases}
\frac{1}{2}\sin^2\tau_i\bar{b_i}b_j^i, & z \in \Omega_{c_1}^i;\\
(\frac{1}{2}\sin^2\tau_i\bar{b_i}b_j^i)\eta(\log r_i), & z \in
\Omega_c^i \text{ and }
c_1<r_i<c;\\
(\frac{1}{2}\sin^2\tau_i\bar{b_i}b_j^i)\eta(\log\rho_i-\log r_i),
& z \in \Omega_c^i \text{ and }
c^{-1}\rho_i<r_i<c_1^{-1}\rho_i;\\
\frac{1}{2}\sin^2\tau_j\bar{b_i^j}b_j, & z \in \Omega_{c_1}^j;\\
(\frac{1}{2}\sin^2\tau_i\bar{b_i^j}b_j)\eta(\log r_j), & z \in
\Omega_c^j \text{ and }
c_1<r_j<c;\\
(\frac{1}{2}\sin^2\tau_i\bar{b_i^j}b_j)\eta(\log\rho_j-\log r_j),
& z \in \Omega_c^j \text{ and }
c^{-1}\rho_j<r_j<c_1^{-1}\rho_j;\\
0, & z \in X\setminus(\Omega_c^i\cup\Omega_c^j);\\
\end{cases}
\end{eqnarray*}
\item if $i\leq m$, then
\begin{eqnarray*}
\widetilde{e_{i\bar i}}(z)=
\begin{cases}
\frac{1}{2}\sin^2\tau_i |b_i|^2, & z \in \Omega_{c_1}^i;\\
(\frac{1}{2}\sin^2\tau_i |b_i|^2)\eta(\log r_i), & z \in
\Omega_c^i \text{ and }
c_1<r_i<c;\\
(\frac{1}{2}\sin^2\tau_i |b_i|^2)\eta(\log\rho_i-\log r_i), & z
\in \Omega_c^i \text{ and }
c^{-1}\rho_i<r_i<c_1^{-1}\rho_i;\\
0, & z \in X\setminus\Omega_c^i.\\
\end{cases}
\end{eqnarray*}
\end{enumerate}
Also, let $\widetilde{f_{i\bar j}}=(\Box+1)\widetilde{e_{i\bar
j}}$. It is clear that the supports of these approximation
functions are contained in the corresponding collars. We have the
following estimates:
\begin{lemma}\label{etildesti}
Let $\widetilde{e_{i\bar j}}$ be the functions constructed above.
Then
\begin{enumerate}
\item $e_{i\bar i}=\widetilde{e_{i\bar i}} +O\big
(\frac{u_i^4}{|t_i|^2}\big )$, if $i \leq m$; \item $e_{i\bar
j}=\widetilde{e_{i\bar j}} +O\big (\frac{u_i^3 u_j^3}{|t_i
t_j|}\big )$, if $i,j \leq m$ and $i\ne j$; \item $e_{i\bar
j}=\widetilde{e_{i\bar j}} +O\big (\frac{u_i^3}{|t_i|}\big )$, if
$i\leq m$ and $j\geq m+1$; \item $\Vert e_{i\bar j}\Vert_0=O(1)$,
if $i,j\geq m+1$.
\end{enumerate}
\end{lemma}
{\bf Proof.} The last claim follows from the maximum principle and
Lemma \ref{esafij1}. To prove the first claim, we note that the
maximum principle implies
\[
\Vert e_{i\bar i}-\widetilde{e_{i\bar i}}\Vert_0\leq \Vert
f_{i\bar i}-\widetilde{f_{i\bar i}}\Vert_0.
\]
Now we compute the right hand side of the above inequality. Since
$\widetilde{f_{i\bar i}}\mid_{X\setminus\Omega_c^i}=0$, by Lemma
\ref{esafij1} we know that $\Vert f_{i\bar i}-\widetilde{f_{i\bar
i}}\Vert_{0,X\setminus\Omega_c^i} =O\big
(\frac{u_i^6}{|t_i|^2}\big )$. On $\Omega_{c_1}^i$ we have
\[
|f_{i\bar i}-\widetilde{f_{i\bar i}}|\leq
|\sin^4\tau_i\bar{p_i}b_i|+|\sin^4\tau_i\bar{b_i}p_i|+
|\sin^4\tau_i\bar{p_i}p_i|=O\big (\frac{u_i^6}{|t_i|^2}\big )
\]
which implies $\Vert f_{i\bar i}-\widetilde{f_{i\bar
i}}\Vert_{0,\Omega_{c_1}^i} =O\big (\frac{u_i^6}{|t_i|^2}\big )$.
On $\Omega_c^i\setminus\Omega_{c_1}^i$ with $c_1\leq r_i \leq c$,
we have
\begin{align*}
\begin{split}
|f_{i\bar i}-\widetilde{f_{i\bar i}}|\leq &
(1-\eta)|b_i|^2\sin^4\tau_i+|\sin^4\tau_i\bar{p_i}b_i|
+|\sin^4\tau_i\bar{b_i}p_i|+|\sin^4\tau_i\bar{p_i}p_i|\\
&+\frac{|b_i|^2u_i^{-2}|\eta''|}{4}\sin^4\tau_i
+\frac{|b_i|^2u_i^{-1}|\eta'|}{2}\sin^2\tau_i|\sin 2\tau_i|\\
=& O\big (\frac{u_i^4}{|t_i|^2}\big ).
\end{split}
\end{align*}
Similarly, on $\Omega_c^i\setminus\Omega_{c_1}^i$ with
$c^{-1}\rho_i \leq r_i \leq c_1^{-1}\rho_i$, we have $|f_{i\bar
i}-\widetilde{f_{i\bar i}}|\leq O\big (\frac{u_i^4}{|t_i|^2}\big
)$. By combining the above estimate, we have $\Vert f_{i\bar
i}-\widetilde{f_{i\bar i}}\Vert_0 =O\big
(\frac{u_i^4}{|t_i|^2}\big )$ which implies the first claim. The
second and the third claims can be proved in a similar way.

\qed

As a corollary we prove the following estimates which are more
refined than those of Trapani's on the Ricci metric \cite{tr1}.
The precise constants of the leading terms will be used later to
compute the curvature of the Ricci metric.
\begin{cor}\label{ricciest}
Let $(t,s)$ be the pinching coordinates. Then we have
\begin{enumerate}
\item $\tau_{i\bar i}=\frac{3}{4\pi^2}\frac{u_i^2}{|t_i|^2}
(1+O(u_0))$ and $\tau^{i\bar i}
=\frac{4\pi^2}{3}\frac{|t_i|^2}{u_i^2} (1+O(u_0))$, if $i \leq m$;
\item $\tau_{i\bar j}=O\bigg
(\frac{u_i^2u_j^2}{|t_it_j|}(u_i+u_j)\bigg )$ and $\tau^{i\bar
j}=O(|t_it_j|)$, if $i,j \leq m$ and $i\ne j$; \item $\tau_{i\bar
j}=O\big (\frac{u_i^2}{|t_i|}\big )$ and $\tau^{i\bar
j}=O(|t_i|)$, if $i\leq m$ and $j\geq m+1$; \item $\tau_{i\bar
j}=O(1)$, if $i,j \geq m+1$.
\end{enumerate}
\end{cor}
\begin{rem}
The second part of the above corollary can be made sharper.
However, it will not be useful for our later estimates.
\end{rem}
{\bf Proof.} The second part of the corollary is obtained by
inverting the matrix $(\tau_{i\bar j})$ in the first part together
with the fact that the matrix $(h_{i\bar j})_{i,j\geq m+1}$ is
nondegenerate which was proved by Masur and the fact that the
matrix $(\tau_{i\bar j})_{i,j\geq m+1}$ is bounded from below by a
constant multiple of the matrix $(h_{i\bar j})_{i,j\geq m+1}$
which was proved by Wolpert.

Now we prove the first part. In the following, we use $C_0$ to
denote all universal constants which may change. Recall that
\begin{eqnarray}\label{use1111}
\tau_{i\bar j}=h^{\alpha\bar\beta}R_{i\bar j\alpha\bar\beta}.
\end{eqnarray}
To prove the last claim, let $i,j\geq m+1$. We first notice that
if $\alpha\ne\beta$ or $\alpha=\beta\geq m+1$, then
$|h^{\alpha\bar\beta}|\Vert A_\alpha\Vert_0\Vert
A_\beta\Vert_0=O(1)$ by Lemma \ref{esafij1} and Corollary
\ref{wpasymp}. In this case, we  have
\begin{align*}
\begin{split}
|R_{i\bar j\alpha\bar\beta}|\leq & \left |\int_X e_{i\bar
j}f_{\alpha\bar\beta}\ dv\right | +\left |\int_X
e_{i\bar\beta}f_{\alpha\bar j}\ dv\right | \leq C_0(\Vert e_{i\bar
j} \Vert_0\Vert f_{\alpha\bar\beta} \Vert_0
+\Vert e_{i\bar\beta} \Vert_0\Vert f_{\alpha\bar j} \Vert_0)\\
\leq & C_0(\Vert f_{i\bar j} \Vert_0\Vert f_{\alpha\bar\beta}
\Vert_0 +\Vert f_{i\bar\beta} \Vert_0\Vert f_{\alpha\bar j}
\Vert_0) =O(1)\Vert A_\alpha\Vert_0\Vert A_\beta\Vert_0
\end{split}
\end{align*}
which implies $|h^{\alpha\bar\beta}R_{i\bar
j\alpha\bar\beta}|=O(1)$. If $\alpha=\beta\leq m$ we have
\begin{align*}
\begin{split}
|R_{i\bar j\alpha\bar\alpha}|\leq & \left |\int_X e_{i\bar
j}f_{\alpha\bar\alpha}\ dv\right | +\left |\int_X
e_{i\bar\alpha}f_{\alpha\bar j}\ dv\right | \leq (\Vert e_{i\bar
j} \Vert_0 |f_{\alpha\bar\alpha}|_{L^1} +\bigg (\int_X
|e_{i\bar\alpha}|^2\ dv\int_X |f_{\alpha\bar j}|^2\ dv
\bigg )^{\frac{1}{2}}\\
\leq & O(1)O\big ( \frac{u_i^3}{|t_i|^2} \big ) +\bigg (\int_X
|f_{i\bar\alpha}|^2\ dv\int_X |f_{\alpha\bar j}|^2\ dv
\bigg )^{\frac{1}{2}}\\
=& O\big ( \frac{u_i^3}{|t_i|^2} \big ) +\bigg (\int_X f_{i\bar
i}f_{\alpha\bar\alpha} \ dv \int_X f_{\alpha\bar \alpha}f_{j\bar
j}\ dv \bigg )^{\frac{1}{2}} \leq O\big ( \frac{u_i^3}{|t_i|^2}
\big ) +\Vert A_i\Vert_0\Vert
A_j\Vert_0|f_{\alpha\bar\alpha}|_{L^1} =O\big (
\frac{u_i^3}{|t_i|^2} \big )
\end{split}
\end{align*}
which implies $|h^{\alpha\bar\alpha}R_{i\bar
j\alpha\bar\alpha}|=O(1)$. So we have proved that last claim.

To prove the third claim, let $i\leq m$ and $j\geq m+1$. If
$\alpha\ne\beta$ or $\alpha=\beta\geq m+1$ in formula
\eqref{use1111}, by using integration by part we have
\begin{align*}
\begin{split}
|R_{i\bar j\alpha\bar\beta}|\leq & \left |\int_X f_{i\bar
j}e_{\alpha\bar\beta}\ dv\right | +\left |\int_X
f_{i\bar\beta}e_{\alpha\bar j}\ dv\right | \leq C_0(\Vert
e_{\alpha\bar\beta} \Vert_0|f_{i\bar j}|_{L^1}
+\Vert e_{\alpha\bar j} \Vert_0|f_{i\bar\beta}|_{L^1})\\
\leq & C_0(\Vert f_{\alpha\bar\beta} \Vert_0|f_{i\bar j}|_{L^1}
+\Vert f_{\alpha\bar j} \Vert_0|f_{i\bar\beta}|_{L^1}) =O\big (
\frac{u_i^3}{|t_i|} \big ) \Vert A_\alpha\Vert_0\Vert
A_\beta\Vert_0 +O(1)\Vert A_\alpha\Vert_0|f_{i\bar\beta}|_{L^1}.
\end{split}
\end{align*}
By the above argument we have $|h^{\alpha\bar\beta}O\big (
\frac{u_i^3}{|t_i|} \big ) \Vert A_\alpha\Vert_0\Vert
A_\beta\Vert_0|=O\big ( \frac{u_i^3}{|t_i|} \big )$ and by Lemma
\ref{esafij1} we have $|h^{\alpha\bar\beta}\Vert
A_\alpha\Vert_0|f_{i\bar\beta}|_{L^1}|= O\big (
\frac{u_i^3}{|t_i|} \big )$. So the claim is true in this case.

If $\alpha=\beta\leq m$ and $\alpha\ne i$, we have
\[
|R_{i\bar j\alpha\bar\alpha}|\leq \left |\int_X f_{i\bar
j}e_{\alpha\bar\alpha}\ dv\right | +\left |\int_X
f_{i\bar\alpha}e_{\alpha\bar j}\ dv\right |.
\]
To estimate the second term in the above formula, we have
\[
\left |\int_X f_{i\bar\alpha}e_{\alpha\bar j}\ dv\right | \leq
\Vert e_{\alpha\bar j} \Vert_0|f_{i\bar\alpha}|_{L^1} \leq \Vert
f_{\alpha\bar j} \Vert_0|f_{i\bar\alpha}|_{L^1} =O\big (
\frac{u_\alpha}{|t_\alpha|} \big ) O\big (
\frac{u_i^3u_\alpha^3}{|t_it_\alpha|} \big ) =O\big (
\frac{u_i^3u_\alpha^4}{|t_i||t_\alpha|^2} \big ).
\]
To estimate the first term, we have
\begin{align*}
\begin{split}
\left |\int_X f_{i\bar j}e_{\alpha\bar\alpha}\ dv\right |\leq &
\left |\int_X f_{i\bar j}\widetilde{e}_{\alpha\bar\alpha}\
dv\right | +\left |\int_X f_{i\bar j}
(e_{\alpha\bar\alpha}-\widetilde{e}_{\alpha\bar\alpha})\ dv\right |\\
\leq & \left |\int_{\Omega_c^\alpha} f_{i\bar
j}\widetilde{e}_{\alpha\bar\alpha}\ dv\right | +\Vert
e_{\alpha\bar\alpha}-\widetilde{e}_{\alpha\bar\alpha}\Vert_0
|f_{i\bar j}|_{L^1}\\
\leq & \Vert f_{i\bar j} \Vert_{0,\Omega_c^\alpha}
|\widetilde{e}_{\alpha\bar\alpha}|_{L^1} +O\big (
\frac{u_\alpha^4}{|t_\alpha|^2} \big ) O\big ( \frac{u_i^3}{|t_i|}
\big ) =O\big ( \frac{u_i^3u_\alpha^3}{|t_i||t_\alpha|^2} \big )
\end{split}
\end{align*}
which implies $|h^{\alpha\bar\alpha}R_{i\bar j\alpha\bar\alpha}|
=O\big ( \frac{u_i^3}{|t_i|} \big )$.

Finally, if $\alpha=\beta=i$, we have
\[
|R_{i\bar j i\bar i}|=2 \left |\int_X f_{i\bar j}e_{i\bar i}\
dv\right | \leq 2 \Vert e_{i\bar i} \Vert_0|f_{i\bar j}|_{L^1}
\leq 2 \Vert f_{i\bar i} \Vert_0|f_{i\bar j}|_{L^1} =O\big (
\frac{u_i^2}{|t_i|^2} \big )O\big ( \frac{u_i^3}{|t_i|} \big )
=O\big ( \frac{u_i^5}{|t_i|^3} \big )
\]
which implies $|h^{i\bar i}R_{i\bar j i\bar i}| =O\big (
\frac{u_i^2}{|t_i|} \big )$. This proves the third claim.

The second claim can be proved in a similar way. Now we prove the
first claim. If $\alpha\ne\beta$ or $\alpha=\beta\geq m+1$ in
formula \eqref{use1111}, we have
\begin{align*}
\begin{split}
|R_{i\bar i\alpha\bar\beta}|\leq & \left |\int_X f_{i\bar
i}e_{\alpha\bar\beta}\ dv\right | +\left |\int_X
f_{i\bar\beta}e_{\alpha\bar i}\ dv\right | \leq \Vert
e_{\alpha\bar\beta} \Vert_0|f_{i\bar i}|_{L^1} +\bigg (\int_X
|e_{\alpha\bar i}|^2 \ dv
\int_X |f_{i\bar\beta}|^2\ dv \bigg )^{\frac{1}{2}}\\
\leq & \Vert f_{\alpha\bar\beta} \Vert_0|f_{i\bar i}|_{L^1} +\bigg
(\int_X |f_{\alpha\bar i}|^2 \ dv \int_X |f_{i\bar\beta}|^2\ dv
\bigg )^{\frac{1}{2}} \leq ( \Vert f_{\alpha\bar\beta} \Vert_0
+\Vert A_{\alpha} \Vert_0\Vert A_{\beta} \Vert_0)|f_{i\bar
i}|_{L^1}
\end{split}
\end{align*}
which implies  $|h^{\alpha\bar\beta}R_{i\bar i\alpha\bar\beta}|
=O\big ( \frac{u_i^3}{|t_i|^2} \big )$.

If $\alpha=\beta\leq m$ and $\alpha\ne i$, we have
\begin{align*}
\begin{split}
|R_{i\bar i\alpha\bar\alpha}|\leq & \left |\int_X e_{i\bar
i}f_{\alpha\bar\alpha}\ dv\right | +\left |\int_X
e_{i\bar\alpha}f_{\alpha\bar i}\ dv\right |.
\end{split}
\end{align*}
To estimate the second term in the above inequality, we have
\[
\left |\int_X e_{i\bar\alpha}f_{\alpha\bar i}\ dv\right | \leq
\Vert e_{i\bar\alpha} \Vert_0 |f_{\alpha\bar i}|_{L^1} \leq \Vert
f_{i\bar\alpha} \Vert_0 |f_{\alpha\bar i}|_{L^1} =O\big (
\frac{u_i u_\alpha}{|t_i t_\alpha|} \big ) O\big ( \frac{u_i^3
u_\alpha^3}{|t_i t_\alpha|} \big ) =O\big ( \frac{u_i^4
u_\alpha^4}{|t_i t_\alpha|^2} \big ).
\]
To estimate the first term in the above inequality, we have
\begin{align*}
\begin{split}
\left |\int_X e_{i\bar i}f_{\alpha\bar\alpha}\ dv\right |\leq &
\left |\int_X \widetilde{e}_{i\bar i}f_{\alpha\bar\alpha}\
dv\right | +\left |\int_X (e_{i\bar i}-\widetilde{e}_{i\bar i})
f_{\alpha\bar\alpha}\ dv\right |\\
\leq & \left |\int_{\Omega_c^i} \widetilde{e}_{i\bar
i}f_{\alpha\bar\alpha}\ dv\right | +\Vert e_{i\bar
i}-\widetilde{e}_{i\bar i}\Vert_0
|f_{\alpha\bar\alpha}|_{L^1}\\
\leq &\Vert f_{\alpha\bar\alpha} \Vert_{0,\Omega_c^i}
|\widetilde{e}_{i\bar i}|_{L^1} +\Vert e_{i\bar
i}-\widetilde{e}_{i\bar i}\Vert_0
|f_{\alpha\bar\alpha}|_{L^1}\\
=& O\big ( \frac{u_\alpha^6}{|t_\alpha|^2} \big ) O\big (
\frac{u_i^3}{|t_i|^2} \big )+ O\big (
\frac{u_\alpha^3}{|t_\alpha|^2} \big ) O\big (
\frac{u_i^4}{|t_i|^2} \big ) =O\big ( \frac{u_i^3 u_\alpha^3}{|t_i
t_\alpha|^2} \big ).
\end{split}
\end{align*}
These imply $|h^{\alpha\bar\alpha}R_{i\bar i\alpha\bar\alpha}|
=O\big ( \frac{u_i^3}{|t_i|^2} \big )$.

Finally, we compute $h^{i\bar i}R_{i\bar i i\bar i}$. Clearly
$R_{i\bar i i\bar i}=2\int_X e_{i\bar i}f_{i\bar i}\ dv$ and
\[
\int_X e_{i\bar i}f_{i\bar i}\ dv =\int_X \widetilde{e}_{i\bar
i}\widetilde{f}_{i\bar i}\ dv +\int_X \widetilde{e}_{i\bar
i}(f_{i\bar i}-\widetilde{f}_{i\bar i})\ dv +\int_X (e_{i\bar
i}-\widetilde{e}_{i\bar i})f_{i\bar i}\ dv.
\]
We also have
\[
\left | \int_X \widetilde{e}_{i\bar i}(f_{i\bar
i}-\widetilde{f}_{i\bar i})\ dv \right | \leq \Vert f_{i\bar
i}-\widetilde{f}_{i\bar i}\Vert_0 |\widetilde{e}_{i\bar
i}|_{L^1}=O\big ( \frac{u_i^7}{|t_i|^4} \big )
\]
and
\[
\left | \int_X f_{i\bar i}(e_{i\bar i}-\widetilde{e}_{i\bar i})\
dv \right | \leq \Vert e_{i\bar i}-\widetilde{e}_{i\bar i}\Vert_0
|f_{i\bar i}|_{L^1}=O\big ( \frac{u_i^7}{|t_i|^4} \big ).
\]
Also, we have $\Vert \widetilde{e}_{i\bar i}
\Vert_{0,\Omega_c^i\setminus\Omega_{c_1}^i} =O\big (
\frac{u_i^4}{|t_i|^2} \big )$ and
 $\Vert \widetilde{f}_{i\bar i}
\Vert_{0,\Omega_c^i\setminus\Omega_{c_1}^i} =O\big (
\frac{u_i^4}{|t_i|^2} \big )$. So
\begin{align*}
\begin{split}
\int_X \widetilde{e}_{i\bar i}\widetilde{f}_{i\bar i}\ dv=&
\int_{\Omega_{c_1}^i} \widetilde{e}_{i\bar i}\widetilde{f}_{i\bar
i}\ dv +\int_{\Omega_c^i\setminus\Omega_{c_1}^i}
\widetilde{e}_{i\bar i}\widetilde{f}_{i\bar i}\ dv
=\frac{3\pi^2}{16}|b_i|^4 u_i(1+O(u_0))+O\big (
\frac{u_i^8}{|t_i|^4} \big ).
\end{split}
\end{align*}
By using Corollary \ref{wpasymp} we have $h^{i\bar i}R_{i\bar i
i\bar i}=\frac{3}{4\pi^2}\frac{u_i^2}{|t_i|^2} (1+O(u_0))$. By
combining the above results we have proved this corollary.

\qed

It is well known that there is a complete asymptotic Poincar\'e
metric $\omega_p$ on $\mathcal{M}_g$. We briefly describe it here.
Please see \cite{ls2} for more details.

Let $\bar M$ be a compact K\"ahler manifold of dimension $m$. Let
$Y\subset \bar M$ be a divisor of normal crossings and let $M=\bar
M\setminus Y$. Cover $\bar M$ by coordinate charts
$U_1,\cdots,U_p,\cdots,U_q$ such that $(\bar
U_{p+1}\cup\cdots\cup\bar U_q)\cap Y=\Phi$. We also assume that,
for each $1\leq \alpha \leq p$, there is a constant $n_\alpha$
such that $U_\alpha\setminus
Y=(\Delta^\ast)^{n_\alpha}\times\Delta^{m-n_\alpha}$ and on
$U_\alpha$, $Y$ is given by $z_1^\alpha\cdots
z_{n_\alpha}^\alpha=0$. Here $\Delta$ is the disk of radius
$\frac{1}{2}$ and $\Delta^\ast$ is the punctured disk of radius
$\frac{1}{2}$. Let $\{\eta_i\}_{1\leq i\leq q}$ be the partition
of unity subordinate to the cover $\{U_i\}_{1\leq i\leq q}$. Let
$\omega$ be a K\"ahler metric on $\bar M$ and let $C$ be a
positive constant. Then for $C$ large, the K\"ahler form
\[
\omega_p=C\omega+\sum_{i=1}^p\sqrt{-1}\partial\bar\partial \bigg
(\eta_i\log\log\frac{1}{z_1^i\cdots z_{n_i}^i}\bigg )
\]
defines a complete metric on $M$ with finite volume since on each
$U_i$ with $1\leq i\leq p$, $\omega_p$ is bounded from above and
below by the local Poincar\'e metric on $U_i$. We call this metric
the asymptotic Poincar\'e metric.

As a direct application of the above corollary, we have
\begin{theorem}
The Ricci metric is equivalent to the asymptotic Poincar\'e
metric. More precisely, there is a positive constant $C$ such that
\[
C^{-1}\omega_p \leq \omega_\tau \leq C\omega_p.
\]
\end{theorem}

Now we estimate the holomorphic sectional curvature of the Ricci
metric. We will show that the holomorphic sectional curvature is
negative in the degeneration directions and is bounded in other
directions. We will need the following estimates on the norms to
estimate the error terms.

\begin{lemma}\label{bochner}
Let $f,g\in C^{\infty}(X,\mathbb C)$ be smooth functions such that
$(\Box+1)f=g$. Then there is a constant $C_0$ such that
\begin{enumerate}
\item $|K_0 f|_{L^2}\leq C_0 |K_0 g|_{L^2}$; \item $|K_1K_0
f|_{L^2}\leq C_0 |K_0 g|_{L^2}$;
\end{enumerate}
\end{lemma}
{\bf Proof.} Let $h=|K_0 f|^2$. By using Schwarz inequality, we
easily see that the lemma follows from the Bochner formula:
\[
\Box h+h+|K_1K_0 f|^2=K_0 f \bar{K_0 g}+\bar{K_0 f} K_0 g
-|f-g|^2.
\]

\qed

We also need the estimates on the sections $K_0 f_{i\bar j}$. We
have:
\begin{lemma}\label{esafij2}
Let $K_0$ and $K_1$ be the Maass operators defined in Section
\ref{sec3}. Then
\begin{enumerate}
\item $\Vert K_0 f_{i\bar i}\Vert_{0,\Omega_c^i} =O\big (
\frac{u_i^2}{|t_i|^2}\big )$ and $\Vert K_0 f_{i\bar
i}\Vert_{0,X\setminus\Omega_c^i} =O\big (
\frac{u_i^6}{|t_i|^2}\big )$, if $i\leq m$; \item $\Vert K_0
f_{i\bar j}\Vert_{0}=O(1)$, if $i,j \geq m+1$; \item  $\Vert K_0
f_{i\bar j}\Vert_{0,\Omega_c^i} =O\big ( \frac{u_i u_j^3}{|t_i
t_j|}\big )$ and $\Vert K_0 f_{i\bar j}\Vert_{0,\Omega_c^j} =O\big
( \frac{u_i^3 u_j}{|t_i t_j|}\big )$ and $\Vert K_0 f_{i\bar
j}\Vert_{0,X\setminus(\Omega_c^i\cup\Omega_c^j)} =O\big (
\frac{u_i^3 u_j^3}{|t_i t_j|}\big )$,\\ if $i,j\leq m$ and $i\ne
j$; \item  $\Vert K_0 f_{i\bar j}\Vert_{0,\Omega_c^i} =O\big (
\frac{u_i}{|t_i|}\big )$ and $\Vert K_0 f_{i\bar
j}\Vert_{0,X\setminus\Omega_c^i} =O\big ( \frac{u_i^3}{|t_i|}\big
)$, if $i\leq m$ and $j\geq m+1$; \item $\Vert f_{i\bar
i}-\widetilde f_{i\bar i}\Vert_{1} =O\big (
\frac{u_i^4}{|t_i|^2}\big )$, if $i\leq m$.
\end{enumerate}
\end{lemma}
This lemma can be proved by using similar methods as we used in
the proof of Lemma \ref{esafij1} together with direct
computations. So are the following $L^1$ and $L^2$ estimates:
\begin{lemma}\label{l1l2}
Let $P=K_1K_0$ be the operator defined Section \ref{sec3}. We have
\begin{enumerate}
\item $|f_{i\bar i}|_{L^2}^2=O\big ( \frac{u_i^5}{|t_i|^4}\big )$,
if $i \leq m$; \item $|K_0 f_{i\bar i}|_{L^2}^2=O\big (
\frac{u_i^5}{|t_i|^4}\big )$, if $i \leq m$; \item $|K_0f_{i\bar
j}|_{L^2}^2=O\big ( \frac{u_i^3u_j^3}{|t_it_j|^2}\big )$, if $i,j
\leq m$ and $i \ne j$; \item $|K_0f_{i\bar j}|_{L^2}^2=O\big (
\frac{u_i^3}{|t_i|^2}\big )$, if $i\leq m$ and $j \geq m+1$; \item
$|K_0f_{i\bar j}|_{L^2}^2=O(1)$, if $i,j \geq m+1$; \item
$|P(\widetilde{e}_{i\bar i})|_{L^1} =O\big (
\frac{u_i^3}{|t_i|^2}\big )$, if $i \leq m$.
\end{enumerate}
\end{lemma}

To estimate the curvature of the Ricci metric by using formula
\eqref{finalcurv}, we first expand the term $\int_X Q_{k\bar
l}(e_{i\bar j}) e_{\alpha\bar\beta}\ dv$. A simple computation
shows that
\begin{lemma}\label{qkl} We have
\begin{align}\nonumber
\begin{split}
\int_X Q_{k\bar l}(e_{i\bar j}) e_{\alpha\bar\beta}\ dv=&-\int_X
f_{k\bar l}(K_0e_{i\bar j}\bar K_0e_{\alpha\bar\beta}
+\bar K_0e_{i\bar j}K_0e_{\alpha\bar\beta})\ dv\\
&-\int_X (\Box e_{i\bar j} K_0e_{\alpha\bar\beta}\bar K_0 e_{k\bar
l} +\Box e_{\alpha\bar\beta}K_0e_{i\bar j} \bar K_0e_{k\bar l})\
dv.
\end{split}
\end{align}
\end{lemma}

To estimate the holomorphic sectional curvature, in formula
\eqref{finalcurv} we let $i=j=k=l$. We decompose
$\widetilde{R}_{i\bar ii\bar i}$ into two parts:
\[
\widetilde{R}_{i\bar ii\bar i}=G_1+G_2
\]
where $G_1$ consists of those terms in the right hand side of
\eqref{finalcurv} with all indices $\alpha$, $\beta$, $\gamma$,
$\delta$, $p$ and $q$ equal to $i$ and $G_2=\widetilde{R}_{i\bar
ii\bar i}-G_1$ consists of those terms in \eqref{finalcurv} where,
in each term, at least one of the indices $\alpha$, $\beta$,
$\gamma$, $\delta$, $p$ or $q$ is not $i$. If $i\leq m$, the
leading term is $G_1$ which is given by
\begin{align}\label{g1}
\begin{split}
G_1=&24h^{i\bar i}\int_{X}(\Box+1)^{-1}
(\xi_i(e_{i\bar i}))\bar\xi_i(e_{i\bar i})\ dv\\
&+6 h^{i\bar i}\int_{X} Q_{i\bar i}
(e_{i\bar i})e_{i\bar i}\ dv\\
&-36 \tau^{i\bar i}(h^{i\bar i})^2 \left | \int_{X}\xi_i(e_{i\bar
i})e_{i\bar i}
\ dv\right |^2\\
&+\tau_{i\bar i}h^{i\bar i}R_{i\bar ii\bar i}.
\end{split}
\end{align}

The main theorem of this section is the following estimate of the
holomorphic sectional curvature of the Ricci metric.
\begin{theorem}\label{mainholo}
Let $X_0\in\bar{\mathcal{M}_g}\setminus\mathcal{M}_g$ be a
codimension $m$ point and let
$(t_1,\cdots,t_m,s_{m+1},\cdots,s_n)$ be the pinching coordinates
at $X_0$ where $t_1,\cdots,t_m$ correspond to the degeneration
directions. Then the holomorphic sectional curvature is negative
in the degeneration directions and is bounded in the
non-degeneration directions. More precisely, there is a $\delta>0$
such that, if $|(t,s)|<\delta$, then
\begin{eqnarray}\label{important100}
\widetilde R_{i\bar ii\bar i}=
\frac{3u_i^4}{8\pi^4|t_i|^4}(1+O(u_0)) >0
\end{eqnarray}
if $i\leq m$ and
\begin{eqnarray}\label{important200}
\widetilde R_{i\bar ii\bar i}=O(1)
\end{eqnarray}
if $i\geq m+1$.

Furthermore, on $\mathcal M_g$, the holomorphic sectional
curvature, the bisectional curvature and the Ricci curvature of
the Ricci metric are bounded from above and below.
\end{theorem}
{\bf Proof.} We first compute the asymptotics of the holomorphic
sectional curvature. By Lemma \ref{qkl} we know that
\[
\int_X Q_{i\bar i}(e_{i\bar i})e_{i\bar i}\ dv =\int_X |K_0
e_{i\bar i}|^2(2e_{i\bar i}-4f_{i\bar i}) \ dv.
\]
By \eqref{g1} we have
\begin{align}\label{g1-10}
\begin{split}
G_1=&24h^{i\bar i}\int_{X}T (\xi_i(e_{i\bar i}))\bar\xi_i(e_{i\bar
i})\ dv +6 h^{i\bar i}\int_{X} |K_0
e_{i\bar i}|^2(2e_{i\bar i}-4f_{i\bar i})\ dv\\
&-36\tau^{i\bar i}(h^{i\bar i})^2 \left | \int_{X}\xi_i(e_{i\bar
i})e_{i\bar i} \ dv\right |^2 +\tau_{i\bar i}h^{i\bar i}R_{i\bar
ii\bar i}.
\end{split}
\end{align}
We first consider the degeneration directions. Assume $i \leq m$.
In this case $G_1$ is the leading term. We have the following
lemma.
\begin{lemma}\label{g2nopf}
If $i\leq m$, then $|G_2|=O\big ( \frac{u_i^5}{|t_i|^4}\big )$.
\end{lemma}
{\bf Proof.} The lemma follows from a case by case check. We will
prove it in the appendix.

\qed

Now we go back to the proof of Theorem \ref{mainholo}. We compute
each term of $G_1$. By the proof of Corollary \ref{ricciest} we
know that $h^{i\bar i}R_{i\bar i i\bar
i}=\frac{3}{4\pi^2}\frac{u_i^2}{|t_i|^2} (1+O(u_0))$. So we have
\begin{align}\label{g1-4}
\begin{split}
\tau_{i\bar i}h^{i\bar i}R_{i\bar ii\bar i}=& \bigg
(\frac{3u_i^2}{4\pi^2|t_i|^2}\bigg )^2(1+O(u_0))
=\frac{9u_i^4}{16\pi^4|t_i|^4}(1+O(u_0)).
\end{split}
\end{align}

Now we compute the second term. We have
\begin{align}\label{g1-2}
\begin{split}
&\int_X |K_0 e_{i\bar i}|^2
(2e_{i\bar i}-4f_{i\bar i})\ dv\\
=&\int_X |K_0 \widetilde e_{i\bar i}|^2 (2\widetilde e_{i\bar
i}-4\widetilde f_{i\bar i})\ dv +\int_X (|K_0 e_{i\bar i}|^2-|K_0
\widetilde e_{i\bar i}|^2)
(2\widetilde e_{i\bar i}-4\widetilde f_{i\bar i})\ dv\\
&+\int_X |K_0 e_{i\bar i}|^2 (2(e_{i\bar i}-\widetilde e_{i\bar
i}) -4(f_{i\bar i}-\widetilde f_{i\bar i}))\ dv.
\end{split}
\end{align}
For the second term in the above equation, we have
\begin{align}\nonumber
\begin{split}
&\left | \int_X (|K_0 e_{i\bar i}|^2-|K_0 \widetilde e_{i\bar
i}|^2) (2\widetilde e_{i\bar i}-4\widetilde f_{i\bar i})\ dv\right
| \leq \Vert |K_0 e_{i\bar i}|^2-|K_0 \widetilde e_{i\bar i}|^2
\Vert_0
\int_X (2|\widetilde e_{i\bar i}|+4|\widetilde f_{i\bar i}|)\ dv\\
\leq & \Vert |K_0 e_{i\bar i}|+|K_0 \widetilde e_{i\bar i}|\Vert_0
\Vert K_0(e_{i\bar i}-\widetilde e_{i\bar i})\Vert_0 \int_X
(2|\widetilde e_{i\bar i}|+4|\widetilde f_{i\bar i}|)\ dv =O\big
(\frac{u_i^2}{|t_i|^2}\big )O\big (\frac{u_i^4}{|t_i|^2}\big
)O\big (\frac{u_i^3}{|t_i|^2}\big )=O\big
(\frac{u_i^9}{|t_i|^6}\big ).
\end{split}
\end{align}
For the second term in the above equation, we have
\begin{align}\nonumber
\begin{split}
&\left | \int_X |K_0 e_{i\bar i}|^2 (2(e_{i\bar i}-\widetilde
e_{i\bar i}) -4(f_{i\bar i}-\widetilde f_{i\bar i}))\ dv\right |
\leq  C_0\Vert K_0 e_{i\bar i}\Vert_0^2 (2\Vert e_{i\bar
i}-\widetilde e_{i\bar i}\Vert_0
+4\Vert f_{i\bar i}-\widetilde f_{i\bar i}\Vert_0)\\
&=O\big (\frac{u_i^4}{|t_i|^4}\big )O\big
(\frac{u_i^4}{|t_i|^2}\big ) =O\big (\frac{u_i^8}{|t_i|^6}\big ).
\end{split}
\end{align}
So we get
\begin{align}\label{g1-20}
\begin{split}
&\int_X |K_0 e_{i\bar i}|^2 (2e_{i\bar i}-4f_{i\bar i})\ dv
=\int_X |K_0 \widetilde e_{i\bar i}|^2 (2\widetilde e_{i\bar
i}-4\widetilde f_{i\bar i})\ dv
+O\big (\frac{u_i^8}{|t_i|^6}\big )\\
=&\int_{\Omega_{c_1}^i} |K_0 \widetilde e_{i\bar i}|^2
(2\widetilde e_{i\bar i}-4\widetilde f_{i\bar i})\ dv
+\int_{\Omega_c^i\setminus\Omega_{c_1}^i} |K_0 \widetilde e_{i\bar
i}|^2 (2\widetilde e_{i\bar i}-4\widetilde f_{i\bar i})\ dv +O\big
(\frac{u_i^8}{|t_i|^6}\big ).
\end{split}
\end{align}
We also have the estimate
\[
\left | \int_{\Omega_c^i\setminus\Omega_{c_1}^i} |K_0 \widetilde
e_{i\bar i}|^2 (2\widetilde e_{i\bar i}-4\widetilde f_{i\bar i})\
dv\right |\leq C_0 \Vert K_0 \widetilde e_{i\bar i}\Vert_0^2(\Vert
\widetilde e_{i\bar
i}\Vert_{0,\Omega_c^i\setminus\Omega_{c_1}^i}+\Vert \widetilde
f_{i\bar i}\Vert_{0,\Omega_c^i\setminus\Omega_{c_1}^i})=O\big
(\frac{u_i^8}{|t_i|^6}\big ).
\]
A direct computation shows that
\[
\int_{\Omega_{c_1}^i} |K_0 \widetilde e_{i\bar i}|^2 (2\widetilde
e_{i\bar i}-4\widetilde f_{i\bar i})\
dv=-\frac{3u_i^7}{64\pi^4|t_i|^6}(1+O(u_0)).
\]
So
\begin{eqnarray}\label{g1-2000}
6 h^{i\bar i}\int_{X} |K_0 e_{i\bar i}|^2(2e_{i\bar i}-4f_{i\bar
i})\ dv =-\frac{9u_i^4}{16\pi^4|t_i|^4}(1+O(u_0)).
\end{eqnarray}

Now we compute the third term. We have
\begin{align}\label{g1-3}
\begin{split}
\int_{X}\xi_i(e_{i\bar i})e_{i\bar i}\ dv=&
\int_{X}\xi_i(\widetilde e_{i\bar i})\widetilde e_{i\bar i}\ dv
+\int_{X}\xi_i(\widetilde e_{i\bar i}) (e_{i\bar i}-\widetilde
e_{i\bar i})\ dv +\int_{X}\xi_i (e_{i\bar i}-\widetilde e_{i\bar
i})e_{i\bar i}\ dv.
\end{split}
\end{align}
By using the same method as above, we obtain
\begin{align}\nonumber
\begin{split}
&\left |\int_{X}\xi_i(\widetilde e_{i\bar i}) (e_{i\bar
i}-\widetilde e_{i\bar i})\ dv\right | \leq C_0\Vert
\xi_i(\widetilde e_{i\bar i})\Vert_0 \Vert e_{i\bar i}-\widetilde
e_{i\bar i}\Vert_0 \leq C_0\Vert A_i\Vert_0 \Vert
K_1K_0(\widetilde e_{i\bar i})\Vert_0\Vert
e_{i\bar i}-\widetilde e_{i\bar i}\Vert_0 \\
\leq  &C_0\Vert A_i\Vert_0 \Vert\widetilde e_{i\bar i}\Vert_2\Vert
e_{i\bar i}-\widetilde e_{i\bar i}\Vert_0 =O\big
(\frac{u_i}{|t_i|}\big )O\big (\frac{u_i^2}{|t_i|^2}\big ) O\big
(\frac{u_i^4}{|t_i|^2}\big )=O\big (\frac{u_i^7}{|t_i|^5}\big )
\end{split}
\end{align}
and
\begin{align}\nonumber
\begin{split}
&\left |\int_{X}\xi_i (e_{i\bar i}-\widetilde e_{i\bar i}
)e_{i\bar i}\ dv \right |\leq \Vert \xi_i (e_{i\bar i}-\widetilde
e_{i\bar i}) \Vert_0 \int_{X}e_{i\bar i}\ dv\leq \Vert A_i\Vert_0
\Vert e_{i\bar i}-\widetilde e_{i\bar i}\Vert_2 h_{i\bar i}\\
\leq &\Vert A_i\Vert_0 \Vert f_{i\bar i}-\widetilde f_{i\bar
i}\Vert_1 h_{i\bar i} =O\big (\frac{u_i}{|t_i|}\big )O\big
(\frac{u_i^4}{|t_i|^2}\big ) O\big (\frac{u_i^3}{|t_i|^2}\big
)=O\big (\frac{u_i^8}{|t_i|^5}\big )
\end{split}
\end{align}
and
\begin{align}\nonumber
\begin{split}
\left | \int_{\Omega_c^i\setminus\Omega_{c_1}^i}\xi_i(\widetilde
e_{i\bar i})\widetilde e_{i\bar i}\ dv \right | \leq C_0 \Vert
\xi_i (\widetilde e_{i\bar i}) \Vert_0 \Vert \widetilde e_{i\bar
i} \Vert_{0,\Omega_c^i\setminus\Omega_{c_1}^i}=O\big
(\frac{u_i^7}{|t_i|^5}\big ).
\end{split}
\end{align}
By putting the above results together, we get
\[
\int_{X}\xi_i(e_{i\bar i})e_{i\bar i}\ dv=
\int_{\Omega_{c_1}^i}\xi_i(\widetilde e_{i\bar i})\widetilde
e_{i\bar i}\ dv+O\big (\frac{u_i^7}{|t_i|^5}\big ).
\]
On $\Omega_{c_1}^i$ we have
\[
\xi_i(\widetilde{e}_{i\bar
i})=-\frac{z_i}{\bar{z_i}}\sin^2\tau_i\bar{b_i}P(\widetilde{e}_{i\bar
i})-\frac{z_i}{\bar{z_i}}\sin^2\tau_i\bar{p_i}P(\widetilde{e}_{i\bar
i}).
\]
However, we have $\Vert
\frac{z_i}{\bar{z_i}}\sin^2\tau_i\bar{p_i}P(\widetilde{e}_{i\bar
i})\Vert_{0,\Omega_{c_1}^i}=O\big (\frac{u_i^5}{|t_i|^3}\big )$
which implies
\[
\left |
\int_{\Omega_{c_1}^i}\frac{z_i}{\bar{z_i}}\sin^2\tau_i\bar{p_i}P(\widetilde{e}_{i\bar
i})\widetilde e_{i\bar i}\ dv\right |=O\big
(\frac{u_i^8}{|t_i|^5}\big ).
\]
A direct computation shows that
\[
\int_{\Omega_{c_1}^i}-\frac{z_i}{\bar{z_i}}\sin^2\tau_i\bar{b_i}P(\widetilde{e}_{i\bar
i}))\widetilde e_{i\bar i}\ dv=-\frac{u_i^6}{32\pi^3|t_i|^4
t_i}(1+O(u_0))
\]
which implies
\[
\int_{X}\xi_i(e_{i\bar i})e_{i\bar i}\ dv=
-\frac{u_i^6}{32\pi^3|t_i|^4 t_i}(1+O(u_0)).
\]
So we obtain
\begin{eqnarray}\label{g1-3000}
36\tau^{i\bar i}(h^{i\bar i})^2 \left | \int_{X}\xi_i(e_{i\bar
i})e_{i\bar i} \ dv\right
|^2=\frac{3u_i^4}{16\pi^4|t_i|^4}(1+O(u_0)).
\end{eqnarray}

Now we estimate the first term. We have
\begin{align*}
\begin{split}
\int_X T\xi_i(e_{i\bar i})\bar{\xi_i}(e_{i\bar i})\ dv =&\int_X
T\xi_i(\widetilde e_{i\bar i})\bar{\xi_i}(\widetilde e_{i\bar i})\
dv+\int_X T\xi_i(e_{i\bar i}-\widetilde e_{i\bar
i})\bar{\xi_i}(\widetilde e_{i\bar i})\ dv\\
&+\int_X T\xi_i(e_{i\bar i})\bar{\xi_i}(e_{i\bar i}-\widetilde
e_{i\bar i})\ dv.
\end{split}
\end{align*}
By using the same method we can get
\begin{align*}
\begin{split}
\left | \int_X T\xi_i(e_{i\bar i}-\widetilde e_{i\bar
i})\bar{\xi_i}(\widetilde e_{i\bar i})\ dv \right |\leq & C_0
\Vert T\xi_i(e_{i\bar i}-\widetilde e_{i\bar i})\Vert_0\Vert
\bar{\xi_i}(\widetilde e_{i\bar i})\Vert_0\leq C_0 \Vert
\xi_i(e_{i\bar i}-\widetilde e_{i\bar i})\Vert_0\Vert
\bar{\xi_i}(\widetilde e_{i\bar i})\Vert_0\\
=& O\big (\frac{u_i^5}{|t_i|^3}\big )O\big
(\frac{u_i^3}{|t_i|^3}\big )=O\big (\frac{u_i^8}{|t_i|^6}\big ).
\end{split}
\end{align*}
Similarly,
\[
\left | \int_X T\xi_i(e_{i\bar i})\bar{\xi_i}(e_{i\bar
i}-\widetilde e_{i\bar i})\ dv\right |=O\big
(\frac{u_i^8}{|t_i|^6}\big ).
\]
So we have
\[
\int_X T\xi_i(e_{i\bar i})\bar{\xi_i}(e_{i\bar i})\ dv =\int_X
T\xi_i(\widetilde e_{i\bar i})\bar{\xi_i}(\widetilde e_{i\bar i})\
dv+O\big (\frac{u_i^8}{|t_i|^6}\big ).
\]
To estimate $T\xi_i(\widetilde e_{i\bar i})$, we introduce another
approximation function. Pick $c_2<c_1$ and let $\eta_1\in
C^{\infty}(\mathbb R, [0,1])$ be the cut-off function defined by
\begin{eqnarray}\label{cutoff1}
\eta_1=
\begin{cases}
\eta_1(x)=1, & x\leq \log c_2;\\
\eta_1(x)=0, & x\geq \log c_1;\\
0<\eta_1(x)<1, & \log c_2<x<\log c_1.
\end{cases}
\end{eqnarray}
For $i\leq m$ define the function $d_i$ by
\begin{eqnarray*}
d_i(z)=
\begin{cases}
-\frac{1}{8}\sin^2\tau_i\cos 2\tau_i |b_i|^2\bar{b_i}, & z \in \Omega_{c_2}^i;\\
(-\frac{1}{8}\sin^2\tau_i\cos 2\tau_i |b_i|^2\bar{b_i})\eta_1(\log
r_i), & z \in \Omega_{c_1}^i \text{ and }
c_2<r_i<c_1;\\
(-\frac{1}{8}\sin^2\tau_i\cos 2\tau_i
|b_i|^2\bar{b_i})\eta_1(\log\rho_i-\log r_i), & z \in
\Omega_{c_1}^i \text{ and }
c_1^{-1}\rho_i<r_i<c_2^{-1}\rho_i;\\
0, & z \in X\setminus\Omega_{c_1}^i.\\
\end{cases}
\end{eqnarray*}
A simple computation shows that
\[
\Vert \xi_i(\widetilde e_{i\bar i})-(\Box+1)d_i\Vert_0 =O\big
(\frac{u_i^5}{|t_i|^3}\big )
\]
which implies
\[
\Vert T \xi_i(\widetilde e_{i\bar i})-d_i\Vert_0 =O\big
(\frac{u_i^5}{|t_i|^3}\big ).
\]
So
\[
\int_X T\xi_i(\widetilde e_{i\bar i})\bar{\xi_i}(\widetilde
e_{i\bar i})\ dv=\int_X d_i\bar{\xi_i}(\widetilde e_{i\bar i})\ dv
+\int_X (T\xi_i(\widetilde e_{i\bar i})-d_i)\bar{\xi_i}(\widetilde
e_{i\bar i})\ dv.
\]
We have the estimate
\[
\left | \int_X (T\xi_i(\widetilde e_{i\bar
i})-d_i)\bar{\xi_i}(\widetilde e_{i\bar i})\ dv \right | \leq
C_0\Vert T \xi_i(\widetilde e_{i\bar i})-d_i\Vert_0 \Vert
\bar{\xi_i}(\widetilde e_{i\bar i})\Vert_0=O\big
(\frac{u_i^8}{|t_i|^6}\big )
\]
which implies
\[
\int_X T\xi_i(e_{i\bar i})\bar{\xi_i}(e_{i\bar i})\ dv=\int_X
d_i\bar{\xi_i}(\widetilde e_{i\bar i})\ dv+O\big
(\frac{u_i^8}{|t_i|^6}\big ).
\]
We also have
\[
d_i\bar{\xi_i}(\widetilde e_{i\bar
i})=-d_i\frac{\bar{z_i}}{z_i}\sin^2\tau_i b_i\bar
P(\widetilde{e}_{i\bar i})-d_i\frac{\bar{z_i}}{z_i}\sin^2\tau_i
p_i\bar P(\widetilde{e}_{i\bar i}).
\]
Since $\Vert d_i\frac{\bar{z_i}}{z_i}\sin^2\tau_i p_i\bar
P(\widetilde{e}_{i\bar i})\Vert_0=O\big (\frac{u_i^8}{|t_i|^6}\big
)$ and $\Vert d_i\frac{\bar{z_i}}{z_i}\sin^2\tau_i b_i\bar
P(\widetilde{e}_{i\bar
i})\Vert_{0,\Omega_{c_1}^i\setminus\Omega_{c_2}^i}=O\big
(\frac{u_i^8}{|t_i|^6}\big )$, we get
\[
\int_X T\xi_i(e_{i\bar i})\bar{\xi_i}(e_{i\bar i})\
dv=\int_{\Omega_{c_2}^i} d_i\bar{\xi_i}(\widetilde e_{i\bar i})\
dv+O\big (\frac{u_i^8}{|t_i|^6}\big ).
\]
A direct computation shows that
\[
\int_X T\xi_i(e_{i\bar i})\bar{\xi_i}(e_{i\bar i})\ dv=
\frac{3u_i^7}{256\pi^4|t_i|^6}(1+O(u_0))
\]
which implies
\begin{eqnarray}\label{g1-1000}
24h^{i\bar i}\int_{X}T (\xi_i(e_{i\bar i}))\bar\xi_i(e_{i\bar i})\
dv =\frac{9u_i^4}{16\pi^4|t_i|^4}(1+O(u_0)).
\end{eqnarray}

By combining formulas \eqref{g1-1000}, \eqref{g1-2000},
\eqref{g1-3000} and \eqref{g1-4} we obtain
\[
G_1=\frac{3u_i^4}{8\pi^4|t_i|^4}(1+O(u_0)).
\]
Together with Lemma \ref{g2nopf} we proved formula
\eqref{important100}. The formula \eqref{important200} can be
proved using similar method with a case by case like the proof of
Lemma \ref{g2nopf}.

Now we give a weak estimate on the full curvature of the Ricci
metric. Let
\begin{enumerate}
\item $\Lambda_i=\frac{u_i}{|t_i|}$ if $i\leq m$; \item
$\Lambda_i=1$ if $i\geq m+1$.
\end{enumerate}
We can check the following estimates by using the methods in the
proof of Lemma \ref{g2nopf}. We have
\begin{eqnarray}\label{gener111}
\widetilde R_{i\bar jk\bar l}=O(1)
\end{eqnarray}
if $i,j,k,l\geq m+1$ and
\begin{eqnarray}\label{gener222}
\widetilde R_{i\bar jk\bar
l}=O(\Lambda_i\Lambda_j\Lambda_k\Lambda_l)O(u_0)
\end{eqnarray}
if at least one of these indices $i,j,k,l$ is less than or equal
to $m$ and they are not all equal to each other.

Now we prove the boundedness of the curvatures. For the
holomorphic sectional curvature, from \eqref{important100} and
\eqref{important200} and Corollary \ref{ricciest}, it is clear
that there is a constant $C_0>1$ depending on $X_0$ and $\delta$
such that if $|(t,s)|\leq\delta$, then
\begin{enumerate}
\item $C_0^{-1}\tau_{i\bar i}^2\leq \widetilde R_{i\bar ii\bar i}
\leq C_0\tau_{i\bar i}^2$, if $i\leq m$; \item $|\widetilde
R_{i\bar ii\bar i}| \leq C_0\tau_{i\bar i}^2$, if $i\geq m+1$.
\end{enumerate}
We cover the divisor $Y=\bar{\mathcal M}_g\setminus \mathcal M_g$
by such open coordinate charts. Since $Y$ is compact, we can pick
finitely many such coordinate charts $\Xi_1,\cdots,\Xi_q$ such
that $Y \subset \bigcup_{s=1}^q \Xi_s$. Clearly  there is an open
neighborhood $N$ of $Y$ such that $\bar N \subset \bigcup_{s=1}^q
\Xi_s$. From formulas \eqref{gener111}, \eqref{gener222} and the
above argument, we know that the holomorphic sectional curvature
of $\tau$ is bounded from above and below on $N$. However,
$\mathcal M_g\setminus N$ is a compact set of $\mathcal M_g$, so
the holomorphic sectional curvature is also bounded on $\mathcal
M_g\setminus N$ which implies the holomorphic sectional curvature
is bounded on $\mathcal M_g$.

The bisectional curvature and the Ricci curvature of the Ricci
metric can be proved to be bounded by using \eqref{gener111},
\eqref{gener222} and a similar argument as above, together with
the covering and compactness argument. This finishes the proof.

\qed

\begin{rem}
The estimates of the bisectional curvature and the Ricci curvature
are not optimal. A sharper estimate will be given in our next
paper \cite{lsy2}.
\end{rem}

\section{The perturbed Ricci metric and its curvatures}
In this section we introduce another new metric, the perturbed
Ricci metric. This metric is obtained by adding a constant
multiple of the Weil-Petersson metric to the Ricci metric. By
doing this we construct a natural complete metric whose
holomorphic sectional curvature is negatively bounded. We will see
that the holomorphic sectional curvature of the perturbed Ricci
metric near an interior point of the moduli space is dominated by
the curvature of the large constant multiple of the Weil-Petersson
metric. Similar argument holds for the holomorphic sectional
curvature of the perturbed Ricci metric in the non-degenerate
directions near a boundary point.

\begin{definition}
For any constant $C>0$, we call the metric
\[
\widetilde\tau_{i\bar j}=\tau_{i\bar j}+Ch_{i\bar j}
\]
the perturbed Ricci metric with constant $C$.
\end{definition}

We first give the curvature formula of the perturbed Ricci metric.
We use $P_{i\bar jk\bar l}$ to denote the curvature tensor of the
perturbed Ricci metric.

\begin{theorem}\label{perriccicurv}
Let $s_1,\cdots,s_n$ be local holomorphic coordinates at $s \in
M_g$. Then at $s$, we have
\begin{align}\label{finalpercurv}
\begin{split}
P_{i\bar j k\bar l}=&h^{\alpha\bar\beta}
\left\{\sigma_1\sigma_2\int_{X_s}
\left\{(\Box+1)^{-1}(\xi_k(e_{i\bar j}))
\bar{\xi}_l(e_{\alpha\bar\beta})+ (\Box+1)^{-1} (\xi_k(e_{i\bar
j})) \bar{\xi}_\beta(e_{\alpha\bar l})
\right\}\ dv\right\}\\
&+h^{\alpha\bar\beta} \left\{\sigma_1\int_{X_s}Q_{k\bar
l}(e_{i\bar j}) e_{\alpha\bar\beta}\ dv
\right\}\\
&-\widetilde\tau^{p\bar q}h^{\alpha\bar\beta}h^{\gamma\bar\delta}
\left\{\sigma_1\int_{X_s}\xi_k(e_{i\bar q}) e_{\alpha\bar\beta}\
dv\right\}\left\{ \widetilde\sigma_1\int_{X_s}\bar{\xi}_l(e_{p\bar
j})
e_{\gamma\bar\delta})\ dv\right\}\\
&+\tau_{p\bar j}h^{p\bar q}R_{i\bar q k\bar l} +CR_{i\bar j k\bar
l}.
\end{split}
\end{align}
\end{theorem}
{\bf Proof.} Let $s_1,\cdots,s_n$ be normal coordinates at a point
$s\in \M_g$ with respect to the Weil-Petersson metric. By formula
\eqref{1st}, at the point $s$ we have
\begin{align}\label{ppp1st}
\begin{split}
\partial_k\widetilde\tau_{i\bar j}&=\partial_k\tau_{i\bar j}
+C\partial_k h_{i\bar j}=h^{\alpha\bar\beta}\left\{ \sigma_1
\int_{X_s} (\xi_k(e_{i\bar j})e_{\alpha\bar\beta})\ dv \right\}
+\tau_{p\bar j}\Gamma_{ik}^p+C\partial_k h_{i\bar j}\\
&=h^{\alpha\bar\beta}\left\{ \sigma_1 \int_{X_s} (\xi_k(e_{i\bar
j})e_{\alpha\bar\beta})\ dv \right\}
\end{split}
\end{align}
since $\Gamma_{ik}^p=\partial_k h_{i\bar j}=0$ at this point. Now
at $s$ the curvature of the Weil-Petersson metric is
\[
R_{i\bar jk\bar l}=\partial_{\bar l}\partial_k h_{i\bar j}.
\]
The theorem follows from formulas \eqref{notation1},
\eqref{ppp1st} and \eqref{key140}.

\qed

Now we estimate the curvature of the perturbed Ricci metric using
formula \eqref{finalpercurv}. The following two linear algebra
lemmas will be used to handle the inverse matrix
$\widetilde\tau^{i\bar j}$ near an interior point and a boundary
point.

\begin{lemma}\label{interior}
Let $D$ be a neighborhood of $0$ in $\mathbb C^n$ and let $A$ and
$B$ be two positive definite $n\times n$ Hermitian matrix
functions on $D$ such that they are bounded from above and below
on $D$ and each entry of them are bounded. Then each entry of the
inverse matrix $(A+CB)^{-1}=O(C^{-1})$ when $C$ is very large.
\end{lemma}

{\bf Proof.} Consider the determinant $\det(A+CB)$. It is a
polynomial of $C$ of degree $n$ and the coefficient of the leading
term is $\det(B)$ which is bounded from below. All other
coefficients are bounded since they only depend on the entries of
$A$ and $B$. So we can pick $C$ large such that $\det(A+CB)\geq
\frac{1}{2}\det(B)C^n$. Now the determinant of the $(i,j)$-minor
of $A+CB$ is a polynomial of $C$ of degree at most $n-1$ and the
coefficients are bounded since they only depend on the entries of
$A$ and $B$. From the fact that the $(i,j)$-entry is the quotient
of the determinant of the $(i,j)$-minor and the determinant of the
matrix $A+CB$, the lemma follows directly.

\qed

\begin{lemma}\label{boundary1000}
Let $X_0\in \bar{\mathcal M}_g$ be a codimension $m$ boundary
point and let $(t_1,\cdots,s_n)$ be the pinching coordinates near
$X_0$. Then for $|(t,s)|<\delta$ with $\delta$ small, we have
that, for any $C>0$,
\begin{enumerate}
\item $0<\widetilde\tau^{i\bar i} < \tau^{i\bar i}$ for all $i$;
\item $\widetilde\tau^{i\bar j}=O(|t_it_j|)$, if $i,j\leq m$ and
$i\ne j$; \item $\widetilde\tau^{i\bar j}=O(|t_i|)$, if $i\leq m$
and $j \geq m+1 $; \item $\widetilde\tau^{i\bar j}=O(1)$, if $i,j
\geq m+1$.
\end{enumerate}
Furthermore, the bounds in the last three claims are independent
of the choice of $C$.
\end{lemma}

{\bf Proof.} The first claim is a general fact of linear algebra.
To prove the last three claims, we denote the submatrices
$(\tau_{i\bar j})_{i,j\geq m+1}$ and $(h_{i\bar j})_{i,j\geq m+1}$
by $A$ and $B$. These two matrices represent the non-degenerate
directions of the Ricci metric and the Weil-Petersson metric
respectively. By the work of Masur, we know that the matrix $B$
can by extended to the boundary non-degenerately. This implies
that $B$ has a positive lower bound. By Corollary \eqref{wpasymp}
we know that $B$ is bounded from above. Now by the work of
Wolpert, since $\omega_\tau\geq \widetilde C\omega_{WP}$ where
$\widetilde C$ only depend on the genus of the Riemann surface, we
know that $A$ has a positive lower bound. By Corollary
\ref{ricciest} we know that $A$ is bounded from above. So both
matrices $A$ and $B$ are bounded from above and below and all
their entries are bounded as long as $|(t,s)|\leq \delta$.

By Corollary \ref{wpasymp} and Corollary \ref{ricciest} we know
that
\begin{displaymath}
(\widetilde\tau_{i\bar j})= \left (
\begin{array}{cc}
\Upsilon & \Psi\\
\bar\Psi^T & A+CB
\end{array}
\right )
\end{displaymath}
where $\Upsilon$ is an $m\times m$ matrix given by
\begin{displaymath}
\Upsilon= \left (
\begin{array}{ccc}
\frac{u_1^2}{|t_1|^2}(\frac{3}{4\pi^2}+\frac{Cu_1}{2})(1+O(u_0)) & \ldots & \frac{u_1^2u_m^2}{|t_1t_m|}(O(u_0)+CO(u_1u_m))\\
\vdots & \vdots & \vdots\\
\frac{u_1^2u_m^2}{|t_1t_m|}(O(u_0)+CO(u_1u_m)) & \ldots &
\frac{u_m^2}{|t_m|^2}(\frac{3}{4\pi^2}+\frac{Cu_m}{2})(1+O(u_0))
\end{array}
\right )
\end{displaymath}
which represent the degenerate directions of the perturbed Ricci
metric and $\Psi$ is an $m\times(n-m)$ matrix given by
\begin{displaymath}
\Psi= \left (
\begin{array}{ccc}
\frac{u_1^2}{|t_1|}(O(1)+CO(u_1)) & \ldots & \frac{u_1^2}{|t_1|}(O(1)+CO(u_1))\\
\vdots & \vdots & \vdots\\
\frac{u_m^2}{|t_m|}(O(1)+CO(u_m)) & \ldots &
\frac{u_m^2}{|t_m|}(O(1)+CO(u_m))
\end{array}
\right )
\end{displaymath}
which represents the mixed directions of the perturbed Ricci
metric.

A direct computation shows that
\[
\det\widetilde\tau=\big\{\prod_{i=1}^m\frac{u_i^2}{|t_i|^2}
(\frac{3}{4\pi^2}+\frac{Cu_i}{2})\big\}\det(A+CB)(1+O(u_0))
\]
where the $O(u_0)$ term is independent of $C$. Let $\Phi_{ij}$ be
the $(i,j)$-minor of $(\widetilde\tau_{i\bar j})$ obtained by
deleting the $i$-th row and $j$-th column of
$(\widetilde\tau_{i\bar j})$. By using the fact that
\[
|\widetilde\tau^{i\bar j}|=\left | \frac{\det
\Phi_{ij}}{\det\widetilde\tau}\right |
\]
the lemma follows from a direct computation of the determinant of
$\Phi_{ij}$.

\qed

Now we prove the main theorem of this section.

\begin{theorem}\label{perholocurv}
For a suitable choice of positive constant $C$, the perturbed
Ricci metric $\widetilde\tau_{i\bar j}=\tau_{i\bar j}+Ch_{i\bar
j}$ is complete and its holomorphic sectional curvatures are
negative and bounded from above and below by negative constants.
Furthermore, the Ricci curvature of the perturbed Ricci metric is
bounded from above and below.
\end{theorem}
{\bf Proof.} It is clear that the metric $\widetilde\tau_{i\bar
j}$ is complete as long as $C\geq 0$ since it is greater than the
Ricci metric which is complete.

Now we estimate the holomorphic sectional curvature. We first show
that, for any codimension $m$ point $X_0 \in \bar{\mathcal
M}_g\setminus \mathcal{M}_g$, there are constants $C_0, \delta>0$
such that, if $(t,s)=(t_1,\cdots,t_m,s_{m+1}, \cdots,s_n)$ is the
pinching coordinates at $p$ with $|(t,s)|<\delta$ and $C\geq C_0$,
the holomorphic sectional curvature of the metric $\widetilde\tau$
is negative. We first consider the degeneration directions. Let
$i=j=k=l\leq m$. As in the proof of Theorem \ref{mainholo}, we let
\begin{align}\label{g1-1000}
\begin{split}
\widetilde G_1=&24h^{i\bar i}\int_{X}T (\xi_i(e_{i\bar
i}))\bar\xi_i(e_{i\bar i})\ dv +6 h^{i\bar i}\int_{X} |K_0
e_{i\bar i}|^2(2e_{i\bar i}-4f_{i\bar i})\ dv\\
&-36\widetilde\tau^{i\bar i}(h^{i\bar i})^2 \left |
\int_{X}\xi_i(e_{i\bar i})e_{i\bar i} \ dv\right |^2 +\tau_{i\bar
i}h^{i\bar i}R_{i\bar ii\bar i}
\end{split}
\end{align}
and $\widetilde G_2$ be the summation of those terms in
\eqref{finalpercurv} in which at least one of the indices
$p,q,\alpha,\beta,\gamma,\delta$ is not $i$. We have $P_{i\bar
ii\bar i}=\widetilde G_1+\widetilde G_2 +CR_{i\bar ii\bar i}$. We
notice here that we can use Lemma \ref{boundary1000} instead of
Corollary \ref{ricciest} in the proof of Lemma \ref{g2nopf}
without changing any estimate. This implies that $|\widetilde G_2|
=O\big (\frac{u_i^5}{|t_i|^4}\big )$. By the proof of Theorem
\ref{mainholo} we have
\begin{eqnarray}\label{tg1}
\widetilde G_1=\bigg ( \frac{9}{16\pi^4}- \frac{3}{16\pi^4} \big(
1+\frac{2\pi^2Cu_i}{3} \big )^{-1}\bigg
)\frac{u_i^4}{|t_i|^4}(1+O(u_0))
\end{eqnarray}
which implies
\begin{eqnarray}\label{tg10}
P_{i\bar ii\bar i}=\left (\bigg ( \frac{9}{16\pi^4}-
\frac{3}{16\pi^4} \big( 1+\frac{2\pi^2Cu_i}{3} \big )^{-1}\bigg
)\frac{u_i^4}{|t_i|^4}
+\frac{3C}{8\pi^2}\frac{u_i^5}{|t_i|^4}\right )(1+O(u_0))>0
\end{eqnarray}
as long as $\delta$ is small enough. Furthermore, $P_{i\bar ii\bar
i}$ is bounded above and below by constant multiple of
$\widetilde\tau_{i\bar i}^2$ where the constants may depend on
$C$. However, when $C$ is fixed, the constants are universal if
$\delta$ is small enough.

Now we let $i=j=k=l\geq m+1$. By the proof of Theorem
\ref{mainholo} and Lemma \ref{boundary1000} we know that $P_{i\bar
ii\bar i}=O(1)+CR_{i\bar ii\bar i}$. We also know that $R_{i\bar
ii\bar i}>0$ has a positive lower bound. Again, by using the
extension theorem of Masur, we can choose $C_0$ large enough such
that, when $C\geq C_0$, we have $P_{i\bar ii\bar i}>0$.
Furthermore,$P_{i\bar ii\bar i}$ is bounded from above and below
by constant multiple of $\widetilde\tau_{i\bar i}^2$ where the
constants may depend on $C$,$m$,$n$, $X_0$ and the choice of
$\nu_{m+1},\cdots,\nu_n$ if $\delta$ is small enough. We also have
estimates similar to \eqref{gener111} and \eqref{gener222}:
\begin{eqnarray}\label{gener333}
P_{i\bar jk\bar l}=O(1)+C R_{i\bar jk\bar l}
\end{eqnarray}
if $i,j,k,l\geq m+1$ and
\begin{eqnarray}\label{gener444}
P_{i\bar jk\bar l}=O(\Lambda_i\Lambda_j\Lambda_k\Lambda_l)O(u_0)+C
R_{i\bar jk\bar l}
\end{eqnarray}
if at least one of these indices $i,j,k,l$ is less than or equal
to $m$ and they are not all equal to each other. So we can choose
$\delta$ small such that, if $|(t,s)|\leq\delta$, then the
holomorphic sectional curvature is bounded from above and below by
negative constants which may depend on $C$.

Now we consider the interior points. Fix a point $p\in\mathcal
M_g$ and a small neighborhood $D$ of $p$ such that  $\bar
D\subset\mathcal M_g$. Since the Ricci metric and Weil-Petersson
metric are uniformly bounded in $\bar D$, we have $P_{i\bar ii\bar
i}=O(1)+CR_{i\bar ii\bar i}$. Using a similar argument as above,
we can choose a $C_0$ such that, when $C>C_0$,  the holomorphic
sectional curvature is bounded from above and below by negative
constants which may depend on $C$.

Since the divisor $\bar{\mathcal M}_g\setminus \mathcal{M}_g$ is
compact, we can find finitely many boundary charts  of $\mathcal
M_g$ described above such that the holomorphic sectional curvature
of $\widetilde\tau$ is pinched by two negative constants which
depend on $C$ on these charts. Furthermore, there is a
neighborhood $N$ of $\bar{\mathcal M}_g\setminus \mathcal{M}_g$ in
$\bar{\mathcal M}_g$ such that $\bar N$ is contained in the union
of these charts. It is clear that we can find a constant $C_1$
such that on $N$, the holomorphic sectional curvature of
$\widetilde\tau$ is pinched by negative constants when $C \geq
C_1$.

Also, since the set $\mathcal M_g\setminus N$ is compact, by the
above argument, we can find finitely many interior charts
described above such that their union covers $\mathcal
M_g\setminus N$ and a constant $C_2$, such that the holomorphic
sectional curvature of $\widetilde\tau$ is pinched by negative
constants when $C>C_2$. Again, the bounds may depend on $C$.By
taking a constant $C>\max\{C_1,C_2\}$, we have proved the first
part of the theorem. The Ricci curvature can be estimated in a
similar way as we did in the proof of Theorem \ref{mainholo}
together with Lemma \ref{interior} and \ref{boundary1000}.

\qed

\begin{rem}
By using the negativity of the Ricci curvature of the
Weil-Petersson metric and estimates \eqref{tg10}, \eqref{gener333}
and \eqref{gener444}, we can actually show that the Ricci
curvature of the perturbed Ricci metric is pinched between two
negative constants. The detail will be given in our next paper.
\end{rem}

\section{Equivalent metrics on the moduli space}
In this section, we prove the equivalence among the Ricci metric,
perturbed Ricci metric, K\"ahler-Einstein metric and the McMullen
metric. These equivalences imply that the Teichm\"uller metric is
equivalent to the K\"ahler-Einstein metric which gives a positive
answer to Yau's Conjecture. The main tool we use is the
Schwarz-Yau Lemma. Also, to control the McMullen metric, we give a
simple formula of the first derivative of the geodesic length
functions.

\begin{lemma}\label{51}
The Weil-Petersson metric is bounded above by a constant multiple
of the Ricci metric. Namely, there is a constant $\alpha>0$ such
that $\omega_{WP}\leq \alpha\omega_\tau$.
\end{lemma}
{\bf Proof.} This lemma follows from Corollary \ref{wpasymp} and
Corollary \ref{ricciest}. It also follows directly from
Schwarz-Yau Lemma.

\qed

By using this simple result, we have
\begin{theorem}\label{52}
The Ricci metric and the perturbed Ricci metric are equivalent.
\end{theorem}
{\bf Proof.} Since $\widetilde\tau_{i\bar j}=\tau_{i\bar
j}+Ch_{i\bar j}$ and $C>0$, we know that the Ricci metric is
bounded above by the perturbed Ricci metric. By using the above
lemma, we also have the bound of the other side.

\qed

By the work of Cheng and Yau \cite{cy1} and Mok and Yau
\cite{mokyau1}, there is a unique complete K\"ahler-Einstein
metric on the moduli space whose Ricci curvature is $-1$. One of
the main results of this section is the equivalence of the
K\"ahler-Einstein metric and the Ricci metric. To prove this
result, we need the following simple fact of linear algebra.

\begin{lemma}\label{lin100}
Let $A$ and $B$ be positive definite $n\times n$ Hermitian
matrices and let $\alpha,\beta$ be positive constants such that
$B\geq \alpha A$ and $\det(B) \leq \beta\det(A)$. Then there is a
constant $\gamma>0$ depending on $\alpha,\beta$ and $n$ such that
$B\leq \gamma A$.
\end{lemma}

\begin{theorem}\label{main100}
The Ricci metric is equivalent to the K\"ahler-Einstein metric
$g_{KE}$.
\end{theorem}
{\bf Proof.} Consider the identity map $i:(\mathcal M_g,g_{KE})\to
(\mathcal M_g, \widetilde\tau)$. We know that the
K\"ahler-Einstein metric is complete and its Ricci curvature is
$-1$. By Theorem \ref{perholocurv} we know that the holomorphic
sectional curvatures of the perturbed Ricci metric is bounded
above by a negative constant. From the Schwarz-Yau Lemma, there is
a constant $c_0>0$ such that
\[
g_{KE} \geq c_0\widetilde\tau.
\]
From Theorem \ref{52} we know that the K\"ahler-Einstein metric is
bounded below by a constant multiple of the Ricci metric
\begin{eqnarray}\label{bw}
g_{KE} \geq \tilde c_0\tau.
\end{eqnarray}
Now we consider the identity map $j:(\mathcal M_g,\tau)\to
(\mathcal M_g, g_{KE})$. By Theorem \ref{mainholo} we know that
the Ricci curvature of the Ricci metric is bounded from below.
Also, the Ricci curvature of the K\"ahler-Einstein metric is $-1$.
From the Schwarz-Yau Lemma for volume forms, there is a constant
$c_1>0$ such that
\begin{eqnarray}\label{abo}
\det(g_{KE})\leq c_1 \det(\tau).
\end{eqnarray}
By combining formula \eqref{bw}, \eqref{abo} and Lemma
\ref{lin100} we have proved the theorem.

\qed

Now we consider the McMullen metric. In \cite{mc} McMullen
constructed a new metric $g_{1/l}$ on $\mathcal M_g$ which is
equivalent to the Teichm\"hller metric and is K\"ahler hyperbolic.
More precisely, let $Log:\ \mathbb R_{+}\to [0,\infty)$ be a
smooth function such that
\begin{enumerate}
\item $Log(x)=\log x$ if $x \geq 2$; \item $Log(x)=0$ if $x \leq
1$.
\end{enumerate}
For suitable choices of small constants $\delta,\epsilon>0$, the
K\"ahler form of the McMullen metric $g_{1/l}$ is
\[
\omega_{1/l}=\omega_{WP}-i\delta\sum_{l_{\gamma}(X)<\epsilon}\partial\bar\partial
Log\frac{\epsilon}{l_\gamma}
\]
where the sum is taken over primitive short geodesics $\gamma$ on
$X$. We will also write this as $\omega_{M}$.

To compare the Ricci metric and the McMullen metric, we compute
the first order derivative of the short geodesics.

\begin{lemma}\label{shortgeo}
Let $X_0\in \bar{\mathcal M}_g$ be a codimension $m$ boundary
point and let $(t_1,\cdots,s_n)$ be the pinching coordinates near
$X_0$. Let $l_j$ be the length of the geodesic on the collar
$\Omega_c^j$. Then
\[
\partial_i l_j=-\pi u_j\bar{b_i^j}
\]
if $i\ne j$ and
\[
\partial_i l_j=-\pi u_j\bar{b_i}
\]
if $i= j$. Here $b_i^j$ and $b_i$ are defined in Lemma \ref{aj10}.
\end{lemma}

{\bf Proof.} It is clear that on the genuine collar $\Omega_c^j$,
$\lambda A_i$ is an anti-holomorphic quadratic differential. By
using the rs-coordinate $z$ on $\Omega_c^j$, we can denote
$\lambda A_i$ by $\kappa_i(\bar z)d\bar z^2$. We consider the
coefficient of the term $\bar z^{-2}$ in the expansion of
$\kappa_i$ and denote it by $C_{-2}(\kappa_i)$. From formula
\eqref{jmetric} and Lemma \ref{aj10} we know that
\begin{eqnarray}\label{geoesti}
C_{-2}(\kappa_i)=\frac{1}{2}u_j^2\bar{b_i^j}.
\end{eqnarray}
Now we use a different way to compute $C_{-2}(\kappa_i)$. Fix
$(t_0,s_0)$ with small norm and let $X=X_{t_0,s_0}$. Let $w$ be
the rs-coordinates on the $j$-th collar of $X_{t,s}$ and let $z$
be the rs-coordinate on the $j$-th collar of $X$. Clearly
$w=w(z,t,s)$ is holomorphic with respect to $z$ and when
$(t,s)=(t_0,s_0)$, we have $w=z$. We pull-back the metric on the
$j$-th collar of $X_{t,s}$ to $X$. We have
\[
\Lambda=\frac{1}{2}u_j^2 |w|^{-2}\csc^2(u_j\log |w|)\left
|\frac{\partial w}{\partial z} \right |^2
\]
is the K\"ahler-Einstein metric on the $j$-th collar of $X_{t,s}$.
Now from formulas \eqref{smalla} and \eqref{biga}, at point
$(t_0,s_0)$, a simple computation shows that
\begin{eqnarray}\label{geoesti10}
\kappa_i(\bar z)=-\frac{u_j\partial_i u_j}{\bar
z^2}+\frac{u_j^2+1}{\bar z^3}\partial_i\bar w
\mid_{(t_0,s_0)}-\frac{u_j^2+1}{\bar z^2}\partial_i\pzb\bar
w\mid_{(t_0,s_0)} -\partial_i\pzb\pzb\pzb\bar w\mid_{(t_0,s_0)}.
\end{eqnarray}
From the above formula we can see that
$C_{-2}(\kappa_i)=-u_j\partial_i u_j$ since the contribution of
the last three terms in the above formula to $C_{-2}(\kappa_i)$ is
$0$. By comparing equations \eqref{geoesti} and \eqref{geoesti10}
we have
\[
\partial_i u_j=-\frac{1}{2}u_j\bar{b_i^j}.
\]
The lemma follows from the fact that $l_j=2\pi u_j$. Again, the
above argument also works when $i=j$. In this case, we replace
$b_i^j$ by $b_i$.

\qed

Now we can prove another main theorem of this section.

\begin{theorem}
The Ricci metric is equivalent to the McMullen metric, the
Teichm\"uller metric and the Kobayashi metric.
\end{theorem}

{\bf Proof.} Royden proved that the Teichm\"uller metric is the
same as the Kobayashi metric. Also, the equivalence of the
McMullen metric and the Teichm\"uller metric was proved by
McMullen \cite{mc}. We only need to show the equivalence between
the Ricci metric and the McMullen $g_{1/l}$ metric.

Since the Ricci curvature of the $g_{1/l}$ metric is bounded from
below and it is complete, by the Schwarz-Yau lemma we know that
\[
\tau<\widetilde\tau\leq C_0 g_{1/l}
\]
for some constant $C_0$. Now we prove the other bound. Fix a
boundary point $X_0$ and the pinching coordinates near $X_0$. By
Theorem 1.1 and Theorem 1.7 of \cite{mc} we know that there are
constants $c_1,c_2$ such that, when $i \leq m$,
\begin{align}\label{mc1000}
\begin{split}
(g_{1/l})_{i\bar i}=&\left\Vert \frac{\partial}{\partial t_i}
\right\Vert^2_{g_{1/l}} <c_1 \left\Vert \frac{\partial}{\partial
t_i}\right\Vert^2_T \leq c_2 \bigg (\left\Vert
\frac{\partial}{\partial
t_i}\right\Vert^2_{WP}+\sum_{l_\gamma<\epsilon}
\left |(\partial\log l_\gamma)\frac{\partial}{\partial t_i}\right|^2\bigg )\\
=&  c_2 \bigg (\left\Vert \frac{\partial}{\partial
t_i}\right\Vert^2_{WP} +\sum_{j=1}^m \left |\partial_i\log
l_j\right |^2\bigg ).
\end{split}
\end{align}
By Lemma \ref{shortgeo} we know that
\[
\left |\partial_i\log l_j\right |^2=\left |\frac{-\pi
u_j\bar{b_i^j}}{l_j}\right |^2 =\frac{1}{4}\left|b_i^j\right |^2.
\]
From Lemma \ref{aj10} we have
\[
\sum_{j=1}^m \left |\partial_i\log l_j\right
|^2=\frac{1}{4}\frac{u_i^2}{\pi^2|t_i|^2} (1+O(u_0)).
\]
From the above formulas and Corollary \ref{wpasymp} and Corollary
\ref{ricciest} we know that there is a constant $c_3$ such that
\[
\left\Vert \frac{\partial}{\partial t_i}\right\Vert^2_{WP}
+\sum_{j=1}^m \left |\partial_i\log l_j\right |^2 \leq
c_3\tau_{i\bar i}
\]
which implies
\begin{eqnarray}\label{mc2000}
(g_{1/l})_{i\bar i}\leq c_4 \tau_{i\bar i}
\end{eqnarray}
where $c_4$ is another constant. The same argument works when
$i\geq m+1$. So formula \eqref{mc2000} holds for all $i$. Since
the McMullen metric is bounded from below by a constant multiple
of the Ricci metric and the diagonal terms of its metric matrix is
bounded from above by a constant multiple of the diagonal terms of
matrix of the Ricci metric, a simple linear algebra fact shows
that there is a constant $c_5$ such that
\[
\tau\geq c_5 g_{1/l}.
\]
The theorem follows from a compactness argument as we have used in
previous sections.

\qed

\section{Appendix: the proof of Lemma 4.10}

We will prove Lemma \ref{g2nopf} in this appendix which consists
of some computational details. We fix a nodal surface $X_0$ which
corresponding to a codimension $m$ boundary point in $\mathcal
M_g$. Let $(t,s)$ be the pinching coordinates near $X_0$ such that
$X_{0,0}=X_0$. Fix $(t,s)$ with small norm, we denote $X_{t,s}$ by
$X$. In the curvature formula \eqref{finalcurv}, we let
$i=j=k=l\leq m$. The term $G_2$ is a summation of the following
four types of terms:
\begin{enumerate}
\item $I=h^{\alpha\bar\beta} \left\{\sigma_1\sigma_2\int_{X}
\left\{T(\xi_k(e_{i\bar j})) \bar{\xi}_l(e_{\alpha\bar\beta})+T
(\xi_k(e_{i\bar j})) \bar{\xi}_\beta(e_{\alpha\bar l}) \right\}\
dv\right\}$ with $(\alpha,\beta)\ne (i,i)$; \item
$II=h^{\alpha\bar\beta} \left\{\sigma_1\int_{X_s}Q_{k\bar
l}(e_{i\bar j}) e_{\alpha\bar\beta}\ dv \right\}$ with
$(\alpha,\beta)\ne (i,i)$; \item $III=\tau^{p\bar
q}h^{\alpha\bar\beta}h^{\gamma\bar\delta}
\left\{\sigma_1\int_{X_s}\xi_k(e_{i\bar q}) e_{\alpha\bar\beta}\
dv\right\}\left\{ \widetilde\sigma_1\int_{X_s}\bar{\xi}_l(e_{p\bar
j})
e_{\gamma\bar\delta})\ dv\right\}$ \\
with $(p,q,\alpha,\beta,\gamma,\delta)\ne (i,i,i,i,i,i)$; \item
$IV=\tau_{p\bar j}h^{p\bar q}R_{i\bar q k\bar l}$ with
$(p,q)\ne(i,i)$
\end{enumerate}
where $T=(\Box+1)^{-1}$. Now we check that the norm of each type
is bounded by $O\big (\frac{u_i^5}{|t_i|^4}\big )$. In the
following, $C_0$ will be a unversal constant which may change but
is independent of the Riemann surface as long as $(t,s)$ has small
norm.

{\bf Case 1.} We check that each term in the sum $IV$ has the
desired bound. By Corollary \ref{ricciest} and its proof we have
\begin{eqnarray*}
R_{i\bar qi\bar i}=
\begin{cases}
O\big (\frac{u_i^5}{|t_i|^3}\big ), \text{ if } q \geq m+1;\\
O\big (\frac{u_i^5u_q^3}{|t_i|^3|t_q|}\big ), \text{ if } q \leq
m, \text{ and }
q\ne i;\\
O\big (\frac{u_i^5}{|t_i|^4}\big ), \text{ if } q= i.
\end{cases}
\end{eqnarray*}
By using the above formula and Corollary \ref{wpasymp} and
\ref{ricciest}, and by a case by case check we have
\[
|\tau_{p\bar i}h^{p\bar q}R_{i\bar qi\bar i}| =O\big
(\frac{u_i^7}{|t_i|^4}\big ).
\]
This proves that the norm of the last term is bounded by $=O\big
(\frac{u_i^5}{|t_i|^4}\big )$.

{\bf Case 2.} We check that each term in the sum $I$ has the
desired bound. Firstly, when $i=j=k=l$, we have
\begin{align}\label{app5}
\begin{split}
&\sigma_1\sigma_2 \left \{ T(\xi_k(e_{i\bar j}))
\bar{\xi}_l(e_{\alpha\bar\beta})+T (\xi_k(e_{i\bar j}))
\bar{\xi}_\beta(e_{\alpha\bar l})\right \}\\
=& 2\left \{  T(\xi_i(e_{i\bar i}))
\bar{\xi}_i(e_{\alpha\bar\beta}) + 2T(\xi_i(e_{i\bar\beta}))
\bar{\xi}_i(e_{\alpha\bar i}) + T(\xi_i(e_{i\bar i}))
\bar{\xi}_\beta(e_{\alpha\bar i})\right \}\\
&+2\left\{ T(\xi_i(e_{\alpha\bar i})) \bar{\xi}_i(e_{i\bar\beta})
+2 T(\xi_i(e_{\alpha\bar\beta})) \bar{\xi}_i(e_{i\bar i}) +
T(\xi_i(e_{\alpha\bar i})) \bar{\xi}_\beta(e_{i\bar i})
\right \}\\
&+2\left\{ T(\xi_\alpha(e_{i\bar i})) \bar{\xi}_i(e_{i\bar \beta})
+T(\xi_\alpha(e_{i\bar\beta})) \bar{\xi}_i(e_{i\bar i})
+T(\xi_\alpha(e_{i\bar i})) \bar{\xi}_\beta(e_{i\bar i})
\right \}\\
&+2T(\xi_\alpha(e_{i\bar\beta})) \bar{\xi}_i(e_{i\bar i}).
\end{split}
\end{align}
We estimate the integration of each term in the above summation.
To estimate these terms, we note that, if $\alpha\ne\beta$ or
$\alpha=\beta\geq m+1$, then
\begin{eqnarray}\label{app10}
\left |h^{\alpha\bar\beta}\Vert f_{\alpha\bar\beta}\Vert_1\right
|=O(1).
\end{eqnarray}
Also, we have
\begin{eqnarray}\label{app20}
\Vert P(e_{\alpha\bar\beta})\Vert_0\leq \Vert
e_{\alpha\bar\beta}\Vert_2 \leq C_0 \Vert
f_{\alpha\bar\beta}\Vert_1.
\end{eqnarray}
These formulae can be checked easily by using Theorem
\ref{laesti}, Corollary \ref{wpasymp}, Lemma \ref{esafij1} and
Lemma \ref{esafij2}.

Now we estimate $\left |h^{\alpha\bar\beta}\int_X T(\xi_i(e_{i\bar
i})) \bar{\xi}_i(e_{\alpha\bar\beta})\ dv\right |$. If
$\alpha\ne\beta$ or $\alpha=\beta\geq m+1$, we have
\begin{align*}
\begin{split}
&\left |\int_X T(\xi_i(e_{i\bar i}))
\bar{\xi}_i(e_{\alpha\bar\beta})\ dv\right |\leq \left ( \int_X
|T(\xi_i(e_{i\bar i}))|^2\ dv
\int_X |\bar\xi_i(e_{\alpha\bar\beta})|^2\ dv\right )^{\frac{1}{2}}\\
\leq&\left ( \int_X |\xi_i(e_{i\bar i})|^2\ dv \int_X
|\bar\xi_i(e_{\alpha\bar\beta})|^2\ dv\right )^{\frac{1}{2}}
=\left ( \int_X f_{i\bar i}|P(e_{i\bar i})|^2\ dv
\int_X f_{i\bar i}|P(e_{\alpha\bar\beta})|^2\ dv\right )^{\frac{1}{2}}\\
\leq & \Vert P(e_{i\bar i})\Vert_0\Vert
P(e_{\alpha\bar\beta})\Vert_0 h_{i\bar i}\leq C_0 \Vert f_{i\bar
i}\Vert_1 \Vert f_{\alpha\bar\beta}\Vert_1 h_{i\bar i} =O\big
(\frac{u_i^5}{|t_i|^4}\big )\Vert f_{\alpha\bar\beta}\Vert_1
\end{split}
\end{align*}
since $\Vert f_{i\bar i}\Vert_1=O\big (\frac{u_i^2}{|t_i|^2}\big
)$. Together with formula \eqref{app10} we have
\[
\left |h^{\alpha\bar\beta}\int_X T(\xi_i(e_{i\bar i}))
\bar{\xi}_i(e_{\alpha\bar\beta})\ dv\right |=O\big
(\frac{u_i^5}{|t_i|^4}\big ).
\]
If $\alpha=\beta\leq m$ and $\alpha\ne i$, we have
\begin{align}\label{app30}
\begin{split}
\left |\int_X T(\xi_i(e_{i\bar i}))
\bar{\xi}_i(e_{\alpha\bar\alpha})\ dv\right |\leq & \left |\int_X
T(\xi_i(e_{i\bar i}))
\bar{\xi}_i(\widetilde{e_{\alpha\bar\alpha}})\ dv\right |\\
&+\left |\int_X T(\xi_i(e_{i\bar i}))
\bar{\xi}_i(e_{\alpha\bar\alpha}-\widetilde{e_{\alpha\bar\alpha}})\
dv\right |.
\end{split}
\end{align}
From Lemma \ref{esafij2} we have
\[
\Vert
P(e_{\alpha\bar\alpha}-\widetilde{e_{\alpha\bar\alpha}})\Vert_0
\leq \Vert
e_{\alpha\bar\alpha}-\widetilde{e_{\alpha\bar\alpha}}\Vert_2 \leq
\Vert f_{\alpha\bar\alpha}-\widetilde{e_{\alpha\bar\alpha}}\Vert_1
=O\big (\frac{u_\alpha^4}{|t_\alpha|^2}\big ).
\]
So
\begin{align}\label{app40}
\begin{split}
&\left |\int_X T(\xi_i(e_{i\bar i}))
\bar{\xi}_i(e_{\alpha\bar\alpha}-\widetilde{e_{\alpha\bar\alpha}})\
dv\right | \leq \Vert
P(e_{\alpha\bar\alpha}-\widetilde{e_{\alpha\bar\alpha}})\Vert_0
\left |\int_X |T(\xi_i(e_{i\bar i}))||A_i|\ dv\right |\\
\leq & \Vert
P(e_{\alpha\bar\alpha}-\widetilde{e_{\alpha\bar\alpha}})\Vert_0
\left ( \int_X |T(\xi_i(e_{i\bar i}))|^2\ dv\int_X f_{i\bar i}\ dv\right )^{\frac{1}{2}}\\
\leq & \Vert
P(e_{\alpha\bar\alpha}-\widetilde{e_{\alpha\bar\alpha}})\Vert_0
\left ( \int_X |\xi_i(e_{i\bar i})|^2\ dv\int_X f_{i\bar i}\ dv\right )^{\frac{1}{2}}\\
=&\Vert
P(e_{\alpha\bar\alpha}-\widetilde{e_{\alpha\bar\alpha}})\Vert_0
\left ( \int_X f_{i\bar i}|P(e_{i\bar i})|^2\ dv\int_X f_{i\bar i}\ dv\right )^{\frac{1}{2}}\\
\leq & \Vert
P(e_{\alpha\bar\alpha}-\widetilde{e_{\alpha\bar\alpha}})\Vert_0
\Vert e_{i\bar i}\Vert_2 h_{i\bar i} =O\big
(\frac{u_\alpha^4}{|t_\alpha|^2}\big )O\big
(\frac{u_i^5}{|t_i|^4}\big ).
\end{split}
\end{align}
Since the support of $\widetilde{e_{\alpha\bar\alpha}}$ is inside
$\Omega_c^\alpha$, we know the support of
$P(\widetilde{e_{\alpha\bar\alpha}})$ is inside $\Omega_c^\alpha$.
From Lemma \ref{l1l2} we have
\begin{align}\label{app50}
\begin{split}
&\left |\int_X T(\xi_i(e_{i\bar i}))
\bar{\xi}_i(\widetilde{e_{\alpha\bar\alpha}})\ dv\right | =\left
|\int_{\Omega_c^\alpha} T(\xi_i(e_{i\bar i}))
\bar{\xi}_i(\widetilde{e_{\alpha\bar\alpha}})\ dv\right |\\
\leq & \Vert A_i\Vert_{0,\Omega_c^\alpha}\Vert T(\xi_i(e_{i\bar
i}))\Vert_0 |P(\widetilde{e_{\alpha\bar\alpha}})|_{L^1} \leq \Vert
A_i\Vert_{0,\Omega_c^\alpha}\Vert \xi_i(e_{i\bar i})\Vert_0
|P(\widetilde{e_{\alpha\bar\alpha}})|_{L^1}\\
=& \Vert A_i\Vert_{0,\Omega_c^\alpha}\Vert A_i\Vert_{0} \Vert
P(e_{i\bar i})\Vert_0 |P(\widetilde{e_{\alpha\bar\alpha}})|_{L^1}
=O\big (\frac{u_i^3}{|t_i|}\big )O\big (\frac{u_i}{|t_i|}\big )
O\big (\frac{u_i^2}{|t_i|^2}\big )O\big (\frac{u_\alpha^3}{|t_\alpha|^2}\big )\\
=& O\big (\frac{u_i^6}{|t_i|^4}\big )O\big
(\frac{u_\alpha^3}{|t_\alpha|^2}\big ).
\end{split}
\end{align}
By combining the inequalities \eqref{app40} and \eqref{app50} we
know that
\[
\left |\int_X T(\xi_i(e_{i\bar i}))
\bar{\xi}_i(e_{\alpha\bar\alpha})\ dv\right |= O\big
(\frac{u_i^5}{|t_i|^4}\big )O\big
(\frac{u_\alpha^3}{|t_\alpha|^2}\big ).
\]
From  Lemma \ref{wpasymp} we have
\[
\left |h^{\alpha\bar\alpha}\int_X T(\xi_i(e_{i\bar i}))
\bar{\xi}_i(e_{\alpha\bar\alpha})\ dv\right |= O\big
(\frac{u_i^5}{|t_i|^4}\big ).
\]
We finish the estimate of the first term in the sum \eqref{app5}.
The integration of other terms in this sum can be estimated in a
similar way.

{\bf Case 3.} We check that each term in the sum $III$ has the
desired bound. By Lemma \ref{ricciest} we first prove that when
$q\ne i$ and $k=i$,
\begin{eqnarray}\label{app60}
\left | h^{\alpha\bar\beta} \left\{\sigma_1\int_{X}\xi_k(e_{i\bar
q}) e_{\alpha\bar\beta}\ dv\right\}\right |=
\begin{cases}
O\big (\frac{u_i^{\frac{5}{2}}}{|t_i|^2}\big )O\big
(\frac{u_q}{|t_q|}\big )
\text{ if } q \leq m\\
O\big (\frac{u_i^{\frac{5}{2}}}{|t_i|^2}\big ) \text{ if } q \geq m+1\\
\end{cases}
\end{eqnarray}
Again, we do a case by base check. First we estimate $\left |
h^{\alpha\bar\beta}\int_{X}\xi_i(e_{i\bar q}) e_{\alpha\bar\beta}\
dv\right |$. If $\alpha\ne\beta$ or $\alpha=\beta\geq m+1$, we
have
\begin{align}\label{app70}
\begin{split}
&\left | \int_{X}\xi_i(e_{i\bar q}) e_{\alpha\bar\beta}\ dv\right
| =\left | \int_{X}e_{i\bar q}\xi_i( e_{\alpha\bar\beta})\
dv\right | \leq \left (\int_{X}|\xi_i( e_{\alpha\bar\beta})|^2\ dv
\int_X |e_{i\bar q}|^2\ dv\right )^{\frac{1}{2}}\\
\leq & \left (\int_{X}f_{i\bar i}|P( e_{\alpha\bar\beta})|^2\ dv
\int_X |f_{i\bar q}|^2\ dv\right )^{\frac{1}{2}} \leq \Vert P(
e_{\alpha\bar\beta})\Vert_0  \left (\int_{X}f_{i\bar i}\ dv
\int_{X}f_{i\bar i}f_{q\bar q}\ dv\right )^{\frac{1}{2}}\\
\leq & \Vert P( e_{\alpha\bar\beta})\Vert_0 \Vert A_q\Vert_0
h_{i\bar i} =O\big (\frac{u_i^3}{|t_i|^2}\big )\Vert
f_{\alpha\bar\beta}\Vert_1\Vert A_q\Vert_0.
\end{split}
\end{align}
This implies
\[
\left | h^{\alpha\bar\beta}\int_{X}\xi_i(e_{i\bar q})
e_{\alpha\bar\beta}\ dv\right |=O\big (\frac{u_i^3}{|t_i|^2}\big
)\Vert A_q\Vert_0.
\]
If $\alpha=\beta\leq m$ and $\alpha\ne i$, we have
\[
\left | \int_{X}\xi_i(e_{i\bar q}) e_{\alpha\bar\alpha}\ dv\right
| \leq \left | \int_{X}\xi_i(e_{i\bar
q})\widetilde{e_{\alpha\bar\alpha}}\ dv\right | +\left |
\int_{X}\xi_i(e_{i\bar q})
(e_{\alpha\bar\alpha}-\widetilde{e_{\alpha\bar\alpha}})\ dv\right
|.
\]
For the second term in the above formula, we have
\begin{align*}
\begin{split}
&\left | \int_{X}\xi_i(e_{i\bar q})
(e_{\alpha\bar\alpha}-\widetilde{e_{\alpha\bar\alpha}})\ dv\right
| =\left | \int_{X}e_{i\bar q}
\xi_i(e_{\alpha\bar\alpha}-\widetilde{e_{\alpha\bar\alpha}})\ dv\right |\\
\leq & \left ( \int_X |e_{i\bar q}|^2\ dv\int_X
|\xi_i(e_{\alpha\bar\alpha}-\widetilde{e_{\alpha\bar\alpha}})|^2\
dv\right )^{\frac{1}{2}} \leq \left ( \int_X |f_{i\bar q}|^2\
dv\int_X f_{i\bar i}
|P(e_{\alpha\bar\alpha}-\widetilde{e_{\alpha\bar\alpha}})|^2\ dv\right )^{\frac{1}{2}}\\
\leq & \Vert
P(e_{\alpha\bar\alpha}-\widetilde{e_{\alpha\bar\alpha}})\Vert_0
 \left ( \int_X f_{i\bar i}f_{q\bar q}\ dv \int_X f_{i\bar i}\ dv\right )^{\frac{1}{2}}
\leq \Vert
e_{\alpha\bar\alpha}-\widetilde{e_{\alpha\bar\alpha}}\Vert_2
\Vert A_q\Vert_0 h_{i\bar i}\\
\leq & \Vert
f_{\alpha\bar\alpha}-\widetilde{e_{\alpha\bar\alpha}}\Vert_2 \Vert
A_q\Vert_0 h_{i\bar i} =O\big (\frac{u_\alpha^4}{|t_\alpha|^2}\big
)O\big (\frac{u_i^3}{|t_i|^2}\big ) \Vert A_q\Vert_0.
\end{split}
\end{align*}
For the first term in the above formula, we have
\begin{align*}
\begin{split}
&\left | \int_{X}\xi_i(e_{i\bar
q})\widetilde{e_{\alpha\bar\alpha}}\ dv\right | =\left |
\int_{\Omega_c^\alpha}\xi_i(e_{i\bar
q})\widetilde{e_{\alpha\bar\alpha}}\ dv\right | \leq \Vert
A_i\Vert_{0,\Omega_c^\alpha} \Vert P(e_{i\bar q})\Vert_0
\int_{\Omega_c^\alpha}\widetilde{e_{\alpha\bar\alpha}}\ dv\\
&\leq \Vert A_i\Vert_{0,\Omega_c^\alpha} \Vert e_{i\bar q}\Vert_2
\int_{\Omega_c^\alpha}\widetilde{e_{\alpha\bar\alpha}}\ dv \leq
O\big (\frac{u_\alpha^3}{|t_\alpha|^2}\big )O\big
(\frac{u_i^3}{|t_i|}\big ) \Vert f_{i\bar q}\Vert_1.
\end{split}
\end{align*}
By combining the above two formulas we have the desired bound for
$\left | h^{\alpha\bar\alpha}\int_{X}\xi_i(e_{i\bar q})
e_{\alpha\bar\alpha}\ dv\right |$.

When $\alpha=\beta=i$, by using a similar method we can show that
$\left | h^{i\bar i}\int_{X}\xi_i(e_{i\bar q}) e_{i\bar i}\
dv\right |=O\big (\frac{u_i^3}{|t_i|^2}\big )\Vert A_q\Vert_0$.
From the above estimates we have proved that the term $\left |
h^{\alpha\bar\beta}\int_{X}\xi_i(e_{i\bar q}) e_{\alpha\bar\beta}\
dv\right |$ in formula \eqref{app60} has the desired estimate. By
 using similar method we can show that the other terms in
\eqref{app60} have the desired estimate. This proves formula
\eqref{app60}.

In a similar way, in the case $q=i$ we can prove that, when $k=i$,
\begin{eqnarray}\label{app80}
\left | h^{\alpha\bar\beta} \left\{\sigma_1\int_{X}\xi_k(e_{i\bar
q}) e_{\alpha\bar\beta}\ dv\right\}\right |=
\begin{cases}
O\big (\frac{u_i^3}{|t_i|^3}\big ), \text{ if } \alpha=\beta=i;\\
O\big (\frac{u_i^4}{|t_i|^3}\big ), \text{ if } \alpha\ne i \text{ or } \beta\ne i.\\
\end{cases}
\end{eqnarray}
By combining formulas \eqref{app70} and \eqref{app80} we conclude
that each term in the sum $III$ is of order $O\big
(\frac{u_i^5}{|t_i|^4}\big )$.

{\bf Case 4.} We need to show that each term in the sum $II$ is of
order $O\big (\frac{u_i^5}{|t_i|^4}\big )$. This case can be
 proved by a case by case check by using the similar estimates as in the third case
together with Lemma \ref{qkl}. This finishes the proof.

\qed

\begin{rem}
The method we estimate these terms can be directly applied to the
computations of the full curvature tensor and we can get certain
bounds for the bisectional curvature and the Ricci curvature of
the Ricci metric as well as the perturbed Ricci metric.
\end{rem}

\end{document}